\documentclass[11pt]{article}
\renewcommand{\baselinestretch}{1.5}
\newcounter{constsub}
\newcommand\crcount[1]{\newcounter{#1}\setcounter{#1}%
{\value{constsub}}\stepcounter{constsub}}

\usepackage{makeidx}
\usepackage{amsmath,amssymb}
\usepackage{verbatim}

\usepackage{exscale}
\usepackage{tabularx}
\usepackage{amsthm}
\usepackage{latexsym}
\usepackage[usenames,dvipsnames]{color}
\usepackage{graphicx,epsfig}
\usepackage{hyperref}
\usepackage{seceqn}
\usepackage{url}

\newcommand{\mb}[1]{\mbox{\boldmath $#1$}}

\newenvironment{pf}{\noindent {\bf Proof. }}{\hfill $\square$}

\newcommand{\bea}{\begin{eqnarray}}
\newcommand{\eea}{\end{eqnarray}}

\newcommand{\req}[1]{$(\ref{#1})$}

\newcommand{\vect}[1]{{\boldsymbol #1 }}

\newcommand{\inprod}[2]{\langle #1 , #2 \rangle }
\newcommand{\card}[1]{\left\lvert #1 \right\rvert}
\newcommand{\bc}{\begin{center}}
\newcommand{\ec}{\end{center}}

\newcommand{\by}{\vect{y}}
\newcommand{\bY}{\vect{Y}}
\newcommand{\refs}[1]{$(\ref{#1})$}

\newcommand{\bX}{\vect{X}}
\newcommand{\bx}{\vect{x}}

\newcommand{\bZ}{\boldsymbol Z}
\newcommand{\R}{\mathbb R}

\newcommand{\m}[1]{\mathcal{#1}}

\newcommand{\be}{\begin{equation}}
\newcommand{\ee}{\end{equation}}
\newcommand{\beaa}{\begin{eqnarray*}}
\newcommand{\eeaa}{\end{eqnarray*}}
\newcommand{\ben}{\begin{enumerate}}
\newcommand{\een}{\end{enumerate}}
\newcommand{\db}{\hspace*{\fill}{\zapf o}}

\newcommand{\bpn}{\begin{proposition}\twlsf}
\newcommand{\epn}{\db\end{proposition}}
\newcommand{\bdm}{\begin{displaymath}}
\newcommand{\edm}{\end{displaymath}}
\newcommand{\ba}{\begin{array}}
\newcommand{\ea}{\end{array}}

\newcommand{\st}{\mathop{\rm s.t.}}

\newtheorem{lemma}{Lemma}
\newtheorem{proposition}{Proposition}
\newtheorem{corollary}{Corollary}

\newtheorem{theorem}{Theorem}

\newtheorem{claim}{Claim}
\newcommand{\norm}[1]{\left\lVert#1\right\rVert}
\newcommand{\tnorm}[1]{\lVert\mkern-2mu |#1|\mkern-2mu\rVert}
\title{Finding the largest low-rank clusters with Ky Fan $2$-$k$-norm and $\ell_1$-norm\thanks{Supported in part by the U.~S.~Air Force 
Office of Scientific Research, a Discovery Grant from the 
Natural
Sciences and Engineering Research Council of Canada, and a grant
from MITACS.}}
\author{Xuan Vinh Doan\thanks{DIMAP and ORMS Group, Warwick Business School, University of Warwick, Coventry, CV4 7AL, United Kingdom, xuan.doan@wbs.ac.uk. The research was partly done when the author was at the Department of Combinatorics and Optimization,
University of Waterloo, Canada.} \and Stephen Vavasis\thanks{Department of
    Combinatorics and Optimization, University of Waterloo, 200
    University Avenue West, Waterloo, ON N2L 3G1, Canada,
    vavasis@uwaterloo.ca.}}

\date{November 2015}

\begin{document}
\maketitle

\begin{abstract}
We propose a convex optimization formulation with the Ky Fan $2$-$k$-norm and $\ell_1$-norm to
find $k$ largest approximately rank-one submatrix blocks of a given nonnegative matrix that has low-rank block diagonal structure with noise. We analyze low-rank and sparsity structures of the optimal solutions using properties of these two matrix norms. We show that, under certain hypotheses, with high probability, the approach can recover rank-one submatrix blocks even when they are corrupted with random noise and inserted into a
much larger matrix with other random noise blocks.
\end{abstract}

\section{Introduction}
\label{sec:intro}
Given a matrix $\mb{A}\in\R^{m\times n}$ that has low-rank block diagonal structure with noise, we would like to find that low-rank block structure of $\mb{A}$. Doan and Vavasis \cite{Doan10} have proposed a convex optimization formulation to find a large approximately rank-one submatrix of $\mb{A}$ with the nuclear norm and $\ell_1$-norm. The proposed LAROS problem (for ``large approximately rank-one submatrix'') in \cite{Doan10} can be used to sequentially extract features in data. For example, given a corpus of documents in some language, it can be used to co-cluster (or bicluster) both terms and documents, i.e., to identify simultaneously both subsets of terms and subsets of documents strongly related to each other from the \emph{term-document matrix} $\mb{A}\in\R^{m\times n}$ of the underlying corpus of $n$ documents with $m$ defined terms (see, for example, Dhillon \cite{Dhillon01}). Here, ``term'' means a word in the language, excluding common words such as articles and prepositions. The $(i,j)$ entry of $\mb{A}$ is the number of occurrences of term $i$ in document $j$, perhaps normalized. Another example is the biclustering of gene expression data to discover expression patterns of gene clusters with respect to different sets of experimental conditions (see the survey by Madeira and Oliveira \cite{Madeira04} for more details). Gene expression data can be represented by a matrix $\mb{A}$ whose rows are in correspondence with different genes and columns are in corresponence with different experimental conditions. The value $a_{ij}$ is the measurement of the expression level of gene $i$ under the experimental condition $j$.    

If the selected terms in a bicluster occur proportionally in the selected documents, we can intuitively assign a topic to that particular term-document bicluster. Similarly, if the expression levels of selected genes are proportional in all selected experiments of a bicluster in the second example, we can identify a expression pattern for the given gene-experimental condition bicluster. Mathematically, for each bicluster $i$, we obtain a subset ${\cal I}_i\subset\{1,\ldots,m\}$ and ${\cal J}_i\subset\{1,\ldots,n\}$ and the submatrix block $\mb{A}({\cal I}_i,{\cal J}_i)$ is approximately rank-one, i.e., $\mb{A}({\cal I}_i,{\cal J}_i)\approx \mb{w}_i\mb{h}_i^T$. Assuming there are $k$ biclusters and ${\cal I}_i\cap{\cal I}_j=\emptyset$ and ${\cal J}_i\cap{\cal J}_j=\emptyset$ for all $i\neq j$, we then have the following approximation:
\be
\label{eq:nmfapprox}
\mb{A}\approx [\bar{\mb{w}}_1,\ldots,\bar{\mb{w}}_k][\bar{\mb{h}}_1,\ldots,\bar{\mb{h}}_k]^T,
\ee 
where $\bar{\mb{w}}_i$ and $\bar{\mb{h}}_i$ are the zero-padded extensions of $\mb{w}_i$ and $\mb{h}_i$ to vectors of length $m$ and $n$ respectively. If the matrix $\mb{A}$ is nonnegative and consists of these $k$ (row- and column-exclusive) biclusters, we may assume that $\mb{w}_i,\mb{h}_i\geq \mb{0}$ for all $i$ (a consequence of Perron-Frobenius theorem, see, for example, Golub and Van Loan \cite{GVL} for more details). Thus $\mb{A}\approx\mb{W}\mb{H}^T$, where $\mb{W},\mb{H}\geq\mb{0}$, which is an approximate \emph{nonnegative matrix factorization} (NMF) of the matrix $\mb{A}$. In this paper, we shall follow the NMF representation to find row- and column-exclusive biclusters. Note that there are different frameworks for biclustering problems such as the graph partitioning models used in Dhillon \cite{Dhillon01}, Tanay et al. \cite{tanay2002discovering}, and Ames \cite{Ames13}, among other models (see, for example, the survey by Nan et al. \cite{fan2010recent}).

Approximate and exact NMF problems are difficult to solve. The LAROS problem proposed by Doan and Vavasis \cite{Doan10} can be used as a subroutine for a greedy algorithm with which columns of $\mb{W}$ and $\mb{H}$ are constructed sequentially. Each pair of columns corresponds to a feature (or pattern) in the original data matrix $\mb{A}$. Given the properties of LAROS problem, the most significant feature (in size and magnitude) will be constructed first with the appropriate parameter. 

The iterated use of the LAROS algorithm of \cite{Doan10} to extract
blocks one at a time, however, will not succeed
in the case that there are two or more hidden blocks of roughly the same magnitude.
In order to avoid this issue, we propose a new convex formulation that allows us to extract several (non-overlapping) features simultaneously. In Section \ref{sec:normopt}, we study the proposed convex relaxation and the properties of its optimal solutions. In Section \ref{sec:recovery}, we provide conditions to recover low-rank block structure of the block diagonal data matrix $\mb{A}$ in the presence of random noise. {Finally, we demonstrate our results with some numerical examples in Section \ref{sec:ex}, including a synthetic biclustering example and a synthetic gene expression example from the previous literature.}

\noindent  
{\bf Notation.} $\inprod{\mb{A}}{\bX}=\mbox{trace}(\mb{A}^T\bX)$ is used to denote the inner product of two matrices $\mb{A}$ and $\bX$ in $\R^{m\times n}$. $\norm{\bX}_1$ means the sum of the absolute values of all entries of $\bX$, i.e., the $\ell_1$-norm of $\mbox{vec}(\bX)$, the long vector constructed by the concatenation of all columns of $\bX$. Similarly, $\norm{\bX}_{\infty}$ is the maximum absolute value of entries of $\bX$, i.e, the $\ell_\infty$-norm of $\mbox{vec}(\bX)$.
\section{Matrix norm minimization}
\label{sec:normopt}
We start with the following general norm minimization problem, which has been considered in \cite{Doan10}.
\be
\label{eq:gnorm}
\ba{rl}
\min & \tnorm{\bX}\\
\st & \inprod{\mb{A}}{\bX}\geq 1,
\ea
\ee
where $\tnorm{\,\cdot\,}$ is an arbitrary norm function on $\R^{m\times n}$. The associated dual norm $\tnorm{\,\cdot\,}^{\star}$ is defined as
\be
\label{eq:dgnorm}
\ba{rl}
\tnorm{\mb{A}}^{\star}=\max & \inprod{\mb{A}}{\bY}\\
\st & \tnorm{\bY}\leq 1.
\ea
\ee 
These two optimization problems are closely related and their relationship is captured in the following lemmas and theorem discussed in Doan and Vavasis \cite{Doan10}.
{
\begin{lemma}
\label{lem:dgnorm}
Matrix $\bX^*$ is an optimal solution of Problem \refs{eq:gnorm} if and only if $\bY^*=\left(\tnorm{\mb{A}}^{\star}\right)\bX^*$ is an optimal solution of Problem \ref{eq:dgnorm}.
\end{lemma}
\begin{lemma}
\label{lem:sgd}
The set of all optimal solutions of Problem \refs{eq:dgnorm} is the subdifferential of the dual norm function $\tnorm{\,\cdot\,}^{\star}$ at $\mb{A}$, $\partial\tnorm{\mb{A}}^{\star}$.
\end{lemma}
}
\begin{theorem}[Doan and Vavasis \cite{Doan10}]
\label{thm:optsol}
The following statements are true:
\ben
\item[(i)] The set of optimal solutions of Problem \refs{eq:gnorm} is $(\tnorm{\mb{A}}^{\star})^{-1}\partial\tnorm{\mb{A}}^{\star}$, where $\partial\tnorm{\,\cdot\,}^{\star}$ is the subdifferential of the dual norm function $\tnorm{\,\cdot\,}^{\star}$.
\item[(ii)] Problem \refs{eq:gnorm} has a unique optimal solution if and only if the dual norm function $\tnorm{\,\cdot\,}^{\star}$ is differentiable at $\mb{A}$.
\een
\end{theorem}
The LAROS problem in \cite{Doan10} belongs to a special class of \refs{eq:gnorm} with parametric matrix norms of the form $\tnorm{\bX}_\theta=\tnorm{\mb{X}}+\theta\Vert\mb{X}\Vert_1$ where $\tnorm{\,\cdot\,}$ is a \emph{unitarily invariant norm} and $\theta$ is a nonnegative parameter, $\theta\geq 0$:
\be
\label{eq:tnorm}
\ba{rl}
\min & \tnorm{\mb{X}}+\theta\Vert\mb{X}\Vert_1\\
\st & \inprod{\mb{A}}{\bX}\geq 1.
\ea
\ee
A norm $\tnorm{\,\cdot\,}$ is unitarily invariant if $\tnorm{\mb{U}\bX\mb{V}}=\tnorm{\bX}$ for all pairs of unitary matrices $\mb{U}$ and $\mb{V}$ (see, for example, Lewis \cite{Lewis95} for more details). For the LAROS problem, $\tnorm{\bX}$ is the nuclear norm, $\tnorm{\bX}=\norm{\bX}_*$, which is the sum of singular values of $\bX$. In order to characterize the optimal solutions of \refs{eq:tnorm}, we need to compute the dual norm $\tnorm{\,\cdot\,}_\theta^\star$:
\be
\label{eq:dtnorm}
\ba{rl}
\tnorm{\mb{A}}_\theta^{\star}=\max & \inprod{\mb{A}}{\bY}\\
\st & \tnorm{\mb{Y}}+\theta\Vert\mb{Y}\Vert_1\leq 1.
\ea
\ee 
The following proposition, which is a straightforward generalization of Proposition 7 in \cite{Doan10}, provides a dual formulation to compute $\tnorm{\,\cdot\,}_\theta^\star$.
\begin{proposition}
\label{prop:dnorm}
The dual norm $\tnorm{\mb{A}}_{\theta}^{\star}$ with $\theta>0$ is the optimal value of the following optimization problem:
\be
\label{eq:tfunc}
\ba{rl}
\tnorm{\mb{A}}_{\theta}^{\star}=\min & \max\left\{\tnorm{\bY}^\star,\theta^{-1}\norm{\bZ}_{\infty}\right\}\\
\st & \bY + \bZ = \mb{A}.
\ea
\ee
\end{proposition}
The optimality conditions of \refs{eq:tnorm} are described in the following proposition, which is again a generalization of Proposition 9 in \cite{Doan10}.
\begin{proposition}
\label{prop:optimality}
Consider a feasible solution $\bX$ of Problem \refs{eq:tnorm}. If there exists $(\bY,\bZ)$ that satisfies the conditions below,
\ben
\item[(i)] $\bY+\bZ=\mb{A}$ and $\tnorm{\bY}^\star=\theta^{-1}\norm{\bZ}_{\infty}$,
\item[(ii)] $\bX\in\alpha\partial\tnorm{\bY}^\star$, $\alpha\geq 0$,
\item[(iii)] $\bX\in\beta\partial\norm{\bZ}_{\infty}$, $\beta\geq 0$,
\item[(iv)] $\alpha+\theta\beta=\left(\norm{\mb{A}}^*_\theta\right)^{-1}$,
\een
then $\bX$ is an optimal solution of Problem \refs{eq:tnorm}. In addition, if
\ben
\item[(v)] $\tnorm{\,\cdot\,}^\star$ is differentiable at $\bY$ or $\norm{\,\cdot\,}_{\infty}$ is differentiable at $\bZ$,
\een
then $\bX$ is the unique optimal solution.
\end{proposition}
The low-rank structure of solutions obtained from the LAROS problem comes from the fact that the dual norm of the nuclear norm is the spectral norm (or $2$-norm), $\norm{\bX}=\sigma_1(\bX)$, the largest singular value of $\bX$. More exactly, it is due to the structure of the subdifferential $\partial\norm{\,\cdot\,}$. According to Zi\c{e}tak \cite{Zietak93}, if $\mb{Y}=\mb{U}\mb{\Sigma}\mb{V}^T$ is a singular value decomposition of $\mb{Y}$ and $s$ is the multiplicity of the largest singular value of $\mb{Y}$, the subdifferential $\partial\norm{\mb{Y}}$ is written as follows:
$$
\partial\norm{\mb{Y}}=\left\{\mb{U}
\begin{bmatrix}
\mb{S} &\mb{0}\cr
\mb{0} & \mb{0}
\end{bmatrix}
\mb{V}^T: \mb{S}\in{\cal S}^s_+, \norm{\mb{S}}_*=1\right\},
$$
where ${\cal S}^s_+$ is the set of positive semidefinite matrices of size $s$. The description of the subdifferential shows that the maximum possible rank of $\bX\in\alpha\partial\norm{\bY}$ is the multiplicity
of the largest singular value of $\mb{Y}$ and if $s=1$, we achieve rank-one solutions. This structural property of the subdifferential $\partial\norm{\,\cdot\,}$ motivates the norm optimization formulation for the LAROS problem, which aims to find a \emph{single} approximately rank-one submatrix of the data matrix $\mb{A}$. We now propose a new pair of norms that would allow us to handle several approximately rank-one submatrices simultaneously instead of individual ones. Let consider the following norm, which we call Ky Fan $2$-$k$-norm given its similar formulation to that of the classical Ky Fan $k$-norm:
\be
\label{eq:kyfank2}
\tnorm{\mb{A}}_{k,2}=\left(\sum_{i=1}^k\sigma_i^2(\mb{A})\right)^{\frac{1}{2}},
\ee
where $\sigma_1\geq\ldots\sigma_k\geq 0$ are the first $k$ largest singular values of $\mb{A}$, $k\leq k_0=\mbox{rank}(\mb{A})$. The dual norm of the Ky Fan $2$-$k$-norm is denoted by $\tnorm{\,\cdot\,}_{k,2}^\star$. According to Bhatia \cite{Bhatia97}, Ky Fan $2$-$k$-norm is a \emph{Q-norm}, which is unitarily invariant (Definition IV.2.9 \cite{Bhatia97}). Since Ky Fan $2$-$k$-norm is unitarily invariant, we can define its corresponding symmetric gauge function, $\norm{\,\cdot\,}_{k,2}:\mathbb{R}^n\rightarrow\R$, as follows:
\be
\label{eq:vkyfank2}
\norm{\bx}_{k,2}=\left(\sum_{i=1}^k\card{x}_{(i)}^2\right)^{\frac{1}{2}},
\ee
where $\card{x}_{(i)}$ is the $(n-i+1)$-st order statistic of $\card{\bx}$. The dual norm of this gauge function (or more exactly, its square), has been used in Argyriou et al. \cite{Argyriou12} as a regularizer in sparse prediction problems. More recently, its matrix counterpart is considered in McDonald et al. \cite{McDonald14} as a special case of the matrix cluster norm defined in \cite{Jacob09}, whose square is used for multi-task learning regularization. On the other hand, the square Ky Fan $2$-$k$-norm is considered as a penalty in low-rank regression analysis in Giraud \cite{Giraud11}. In this paper, we are going to use dual Ky Fan $2$-$k$-norm, not its square, in our formulation given its structural properties, which will be explained later.

When $k=1$, the Ky Fan $2$-$k$-norm becomes the spectral norm, whose subdifferential has been used to characterize the low-rank structure of the optimal solutions of the LAROS problem. We now propose the following optimization problem, of which the LAROS problem is a special instance with $k=1$:
\be
\label{eq:nprob}
\ba{rl}
\min & \tnorm{\bX}_{k,2}^\star+\theta\norm{\bX}_1\\
\st & \inprod{\mb{A}}{\bX}\geq 1,
\ea
\ee
where $\theta$ is a nonnegative parameter, $\theta\geq 0$. The proposed formulation is an instance of the parametric problem \refs{eq:tnorm} and we can use results obtained in Proposition \ref{prop:dnorm} and \ref{prop:optimality} to characterize its optimal solutions. Before doing so, we first provide an equivalent semidefinite optimization formulation for \refs{eq:nprob} in the following proposition.  
\begin{proposition}
\label{prop:sdpequiv}
Assuming $m\geq n$, the optimization problem \refs{eq:nprob} is then equivalent to the following semidefinite optimization problem:
\be
\label{eq:sdpprob}
\ba{rl}
\displaystyle\min_{p,\vect{P},\vect{Q},\vect{R},\bX} & p + \mbox{trace}(\mb{R}) + \theta\inprod{\mb{E}}{\mb{Q}}\\
\st & kp - \mbox{trace}(\mb{P}) = 0,\\
\quad & p\mb{I} -\mb{P}\succeq 0,\\
\quad & \begin{pmatrix}
\mb{P} & -\frac{1}{2}\mb{X}^T\\
-\frac{1}{2}\mb{X} & \mb{R}
\end{pmatrix}\succeq 0,\\
&\mb{Q}\geq\bX,\,\mb{Q}\geq-\bX,\\
&\inprod{\mb{A}}{\bX}\geq 1,
\ea
\ee
where $\mb{E}$ is the matrix of all ones.
\end{proposition}

\begin{pf}
We first consider the dual norm $\tnorm{\bX}_{k,2}^\star$. We have:
\be
\label{eq:mk2dual}
\ba{rl}
\tnorm{\mb{X}}^\star_{k,2}=\displaystyle\max & \inprod{\mb{X}}{\bY}\\
\st & \tnorm{\bY}_{k,2}\leq 1.
\ea
\ee
Since $m\geq n$, we have: $\left(\tnorm{\bY}_{k,2}\right)^2 = \tnorm{\bY^T\bY}_k$, where $\tnorm{\,\cdot\m}_k$ is the Ky Fan $k$-norm, i.e., the sum of $k$ largest singular values. Since $\bY^T\bY$ is symmetric, $\tnorm{\bY^T\bY}_k$ is actually the sum of $k$ largest eigenvalues of $\bY^T\bY$. Similar to $\norm{\bx}_k$, which is the sum of $k$ largest elements of $\bx$, we obtain the following (dual) optimization formulation for $\tnorm{\bY^T\bY}_k$ (for example, see Laurent and Rendl \cite{Laurent05}):
$$
\ba{rl}
\tnorm{\bY^T\bY}_k = \min & kz + \mbox{trace}(\mb{U})\\
\st & z\mb{I} + \mb{U}\succeq \bY^T\bY,\\
\quad & \mb{U}\succeq 0.
\ea
$$
Applying the Schur complement, we have:
$$
\ba{rl}
\tnorm{\bY^T\bY}_k = \min & kz + \mbox{trace}(\mb{U})\\
\st & \begin{pmatrix}
z\mb{I} + \mb{U} & \bY^T\\
\bY & \mb{I}
\end{pmatrix}\succeq 0,\\
\quad & \mb{U}\succeq 0.
\ea
$$
Thus, the dual norm $\tnorm{\,\cdot\,}_{k,2}^\star$ can be computed as follows:
$$
\ba{rl}
\tnorm{\mb{X}}_{k,2}^\star = \max & \inprod{\mb{X}}{\bY}\\
\st & kz + \mbox{trace}(\mb{U})\leq 1,\\
\quad & \begin{pmatrix}
z\mb{I} + \mb{U} & \bY^T\\
\bY & \mb{I}
\end{pmatrix}\succeq 0,\\
\quad & \mb{U}\succeq 0.
\ea
$$
Applying strong duality theory under Slater's condition, we have:
\be
\label{eq:msk2dual}
\ba{rl}
\tnorm{\mb{X}}_{k,2}^\star = \min & p + \mbox{trace}(\mb{R})\\
\st & kp - \mbox{trace}(\mb{P}) = 0,\\
\quad & p\mb{I} -\mb{P}\succeq 0,\\
\quad & \begin{pmatrix}
\mb{P} & -\frac{1}{2}\mb{X}^T\\
-\frac{1}{2}\mb{X} & \mb{R}
\end{pmatrix}\succeq 0.
\ea
\ee
The reformulation of $\norm{\bX}_1$ is straightforward with the new decision variable $\mb{Q}$ and additional constraints $\mb{Q}\geq\bX$ and $\mb{Q}\geq-\bX$, given the fact that the main problem is a minimization problem.
\end{pf}

Proposition \ref{prop:sdpequiv} indicates that in general, we can solve \refs{eq:nprob} by solving its equivalent semidefinite optimization formulation \refs{eq:sdpprob} with any SDP solver. We are now ready to study some properties of optimal solutions of \refs{eq:nprob}. We have: $\tnorm{\bX}_{k,2}^\star+\theta\norm{\bX}_1$ is a norm for $\theta\geq 0$ and we denote it by $\tnorm{\bX}_{k,2,\theta}$. According to Proposition \ref{prop:dnorm}, the dual norm $\tnorm{\bX}_{k,2,\theta}^\star$,
\be
\label{eq:npk2dual}
\ba{rl}
\tnorm{\mb{A}}^\star_{k,2,\theta}=\displaystyle\max & \inprod{\mb{A}}{\bX}\\
\st & \norm{\bX}_{k,2,\theta}\leq 1,
\ea
\ee
can be calculated by solving the following optimization problem given $\theta>0$:
\be
\label{eq:ndk2dual}
\ba{rl}
\tnorm{\mb{A}}_{k,2,\theta}^\star=\min & \max\left\{\tnorm{\bY}_{k,2},\theta^{-1}\norm{\bZ}_{\infty}\right\}\\
\st & \bY + \bZ = \mb{A}.
\ea
\ee
Similar to Proposition \ref{prop:optimality}, we can provide the optimality conditions for \refs{eq:nprob} in the following proposition.
\begin{proposition}
\label{prop:noptimality}
Consider a feasible solution $\bX$ of Problem \refs{eq:nprob}. If there exists $(\bY,\bZ)$ that satisfies the conditions below,
\ben
\item[(i)] $\bY+\bZ=\mb{A}$ and $\tnorm{\bY}_{k,2}=\theta^{-1}\norm{\bZ}_{\infty}$,
\item[(ii)] $\bX\in\alpha\partial\tnorm{\bY}_{k,2}$, $\alpha\geq 0$,
\item[(iii)] $\bX\in\beta\partial\norm{\bZ}_{\infty}$, $\beta\geq 0$,
\item[(iv)] $\alpha+\theta\beta=\left(\norm{\mb{A}}^*_{k,2,\theta}\right)^{-1}$,
\een
then $\bX$ is an optimal solution of Problem \refs{eq:tnorm}. In addition, if
\ben
\item[(v)] $\tnorm{\,\cdot\,}_{k,2}$ is differentiable at $\bY$ or $\norm{\,\cdot\,}_{\infty}$ is differentiable at $\bZ$,
\een
then $\bX$ is the unique optimal solution.
\end{proposition}
The optimality conditions presented in Proposition \ref{prop:noptimality} indicate that some properties of optimal solutions of \refs{eq:nprob} can be derived from the structure of $\partial\tnorm{\,\cdot\,}_{k,2}$. We shall characterize the subdifferential $\partial\tnorm{\,\cdot\,}_{k,2}$ next. According to Watson \cite{Watson93}, since $\tnorm{\,\cdot\,}_{k,2}$ is a unitarily invariant norm, $\partial\tnorm{\mb{A}}_{k,2}$ is related to $\partial\norm{\sigma(\mb{A})}_{k,2}$, where $\sigma(\mb{A})$ is the vector of singular values of $\mb{A}$. Let $\mb{A}\neq\mb{0}$ be a matrix with singular values that satisfy
$$
\sigma_1\geq\ldots>\sigma_{k-t+1}=\ldots=\sigma_k=\ldots=\sigma_{k+s}>\ldots\geq\sigma_p,
$$
where $p=\min\{m,n\}$, so that the multiplicity of $\sigma_k$ is $s+t$. The subdifferential $\partial\norm{\mb{\sigma}}_{k,2}$ is characterized in the following lemma.
\begin{lemma}
\label{lem:sk2vnorm}
$\mb{v}\in\partial\norm{\mb{\sigma}}_{k,2}$ if and only $\mb{v}$ satisfies the following conditions:
\ben
\item[(i)] $\displaystyle v_i=\frac{\sigma_i}{\norm{\mb{\sigma}}_{k,2}}$ for all $i=1,\ldots,k-t$.
\item[(ii)] $\displaystyle v_i = \tau_i\frac{\sigma_k}{\norm{\mb{\sigma}}_{k,2}}$, $0\leq \tau_i\leq 1$ for all $i=k-t+1,\ldots,k+s$, and $\displaystyle\sum_{i=k-t+1}^{k+s}\tau_i=t$.
\item[(iii)] $v_i=0$ for all $i=k+s+1,\ldots,p$.
\een
\end{lemma}

\begin{pf}
Let ${\cal N}_k$ be the collection of all subsets with $k$ elements of $\{1,\ldots,p\}$, we have:
$$
\norm{\mb{\sigma}}_{k,2}=\max_{N\in{\cal N}_k}f_N(\mb{\sigma}),
$$
where $\displaystyle f_N(\mb{\sigma})=\left(\sum_{i\in N}\sigma_i^2\right)^{\frac{1}{2}}$ for all $N\in{\cal N}_k$. According to Dubovitski-Milyutin's theorem (see, for example, Tikhomirov \cite{tikhomirov1998principles}), the subdifferential of $\norm{\,\cdot\,}_{k,2}$ is computed as follows:
$$
\partial\norm{\mb{\sigma}}_{k,2}=\mbox{conv}\left\{\partial f_N(\mb{\sigma})\,:\,N\in{\cal N}_k,\,f_N(\mb{\sigma})=\norm{\mb{\sigma}}_{k,2}\right\}.
$$
With the structure of $\mb{\sigma}$, clearly $\{1,\ldots,k-t\}\in N$ for all $N\in{\cal N}_k$ such that $f_N(\mb{\sigma})=\norm{\mb{\sigma}}_{k,2}$. The remaining $t$ elements of $N$ are chosen from $s+t$ values from $\{k-t+1,\ldots,k+s\}$. Since $\mb{\sigma}\neq\mb{0}$, all $f_N$ that satisfy $f_N(\mb{\sigma})=\norm{\mb{\sigma}}_{k,2}$ is differentiable at $\mb{\sigma}$ (even in the case $\sigma_k=0$) and
$$
\frac{\partial f_N(\mb{\sigma})}{\partial\sigma_i}=\frac{\sigma_i}{\norm{\mb{\sigma}}_{k,2}},\,\forall\,i\in N,\quad \frac{\partial f_N(\mb{\sigma})}{\partial\sigma_i}=0,\,\quad\,i\not\in N.
$$
Thus if $\mb{v}\in\partial\norm{\mb{\sigma}}_{k,2}$, for all $i=1,\ldots,k-t$, we have: $\displaystyle v_i=\frac{\sigma_i}{\norm{\mb{\sigma}}_{k,2}}$ and $v_i=0$ for all $i=k+s+1,\ldots,p$. 

We now have: counting arguments for the appearance of each index in $\{k-t+1,\ldots,k+s\}$ with respect to all subsets $N\in{\cal N}_k$ that satisfy $f_N(\mb{\sigma})=\norm{\mb{\sigma}}_{k,2}$ allow us to characterize $v_i$ for $i=k-t+1,\ldots,k+s$ as $\displaystyle v_i = \tau_i\frac{\sigma_k}{\norm{\mb{\sigma}}_{k,2}}$, $0\leq \tau_i\leq 1$ and $\displaystyle\sum_{i=k-t+1}^{k+s}\tau_i=t$.
\end{pf}

We are ready to characterize the subdifferential of $\tnorm{\,\cdot\,}_{k,2}$ with the following proposition.
\begin{proposition}
\label{prop:sk2mnorm}
Consider $\mb{A}\neq\mb{0}$. Let $\mb{A}=\mb{U}\mb{\Sigma}\mb{V}^T$ be a particular singular value decomposition of $\mb{A}$ and assume that $\sigma(\mb{A})$ satisfies $\sigma_1\geq\ldots>\sigma_{k-t+1}=\ldots=\sigma_k=\ldots=\sigma_{k+s}>\ldots\geq\sigma_p$. Then, $\mb{G}\in\partial\tnorm{\mb{A}}_{k,2}$ if and only if there exists $\mb{T}\in\R^{(s+t)\times (s+t)}$ such that
$$
\mb{G}=\frac{1}{\tnorm{\mb{A}}_{k,2}}\left(\mb{U}_{[:,1:k-t]}\mb{\Sigma}_{[1:k-t,1:k-t]}\mb{V}_{[:,1:k-t]}^T + \sigma_k\mb{U}_{[:,k-t+1:k+s]}\mb{T}\mb{V}_{[:,k-t+1:k+s]}^T\right),
$$
where $\mb{T}$ is symmetric positive semidefinite, $\norm{\mb{T}}\leq 1$ and $\norm{\mb{T}}_*=t$.
\end{proposition}

\begin{pf}
According to Watson \cite{Watson93}, we have:
$$
\partial\tnorm{\mb{A}}_{k,2}=\left\{\mb{U}\mbox{Diag}(\mb{g})\mb{V}^T:\,\mb{A}=\mb{U}\mb{\Sigma}\mb{V}^T\,\mbox{is any SVD of}\,\mb{A},\,\mb{g}\in\partial\norm{\sigma(\mb{A})}_{k,2}\right\}.
$$
Let $\mb{A}=\mb{U}\mb{\Sigma}\mb{V}^T$ be a particular singular value decomposition of $\mb{A}$ and assume that a singular value $\sigma_i>0$ has the multiplicity of $r$ with corresponding singular vectors $\mb{U}_i\in\R^{m\times r}$ and $\mb{V}_i\in\R^{n\times r}$. Then for any singular value decomposition of $\mb{A}$, $\mb{A}=\bar{\mb{U}}\mb{\Sigma}\bar{\mb{V}}^T$, there exists an orthonormal matrix $\mb{W}\in\R^{r\times r}$, $\mb{W}\mb{W}^T=\mb{I}$, such that $\bar{\mb{U}}_i=\mb{U}_i\mb{W}$ and $\bar{\mb{V}}_i=\mb{V}_i\mb{W}$ (for example, see Zi\c{e}tak \cite{Zietak93}).

Combining these results with Lemma \ref{lem:sk2vnorm}, the proof is straightforward with a singular value (or eigenvalue) decomposition of matrix $\mb{T}$. 
\end{pf}

\begin{corollary}
\label{col:diff}
$\tnorm{\,\cdot\,}_{k,2}$ is differentiable at any $\mb{A}\neq\mb{0}$ such that $\sigma_k>\sigma_{k+1}$ ($\sigma_{p+1}=0$) or $\sigma_k=0$.
\end{corollary}

\begin{pf}
If $\sigma_k=0$, then, according to Proposition \ref{prop:sk2mnorm},
$$
\mb{G}\in\partial\tnorm{\mb{A}}_{k,2}\,\Leftrightarrow\,\mb{G}=\frac{1}{\tnorm{\mb{A}}_{k,2}}\mb{U}_{[:,1:k-t]}\mb{\Sigma}_{[1:k-t,1:k-t]}\mb{V}_{[:,1:k-t]}^T.
$$ 
 Now, if $\sigma_k>\sigma_{k+1}$, we have: $s=0$, thus $\mb{T}=\mb{I}$ is unique since $\mb{T}\in{\cal S}^t$, $\norm{\mb{T}}_*=t$, and $\norm{\mb{T}}\leq 1$. Thus $\partial\tnorm{\mb{A}}_{k,2}$ is a singleton, which implies $\tnorm{\,\cdot\,}_{k,2}$ is differentiable at $\mb{A}$. 
\end{pf}

Proposition \ref{prop:sk2mnorm} shows that the problem \refs{eq:nprob} with $\theta=0$ is a convex optimization problem that finds $k$-approximation of a matrix $\mb{A}$. It also shows that intuitively, the problem \refs{eq:nprob} can be used to recover $k$ largest approximately rank-one submatrices with $\theta>0$. {Note that for Ky Fan $k$-norm, if $\sigma_k(\mb{A})>\sigma_{k+1}(\mb{A})$, its subdifferential at $\mb{A}$ is a singleton with a unique subgradient:
$$
\partial\tnorm{\mb{A}}_k=\left\{\mb{U}
\begin{bmatrix}
\mb{I}_k &\mb{0}\cr
\mb{0} & \mb{0}
\end{bmatrix}
\mb{V}^T\right\},
$$
where $\mb{A}=\mb{U}\Sigma\mb{V}^T$ is a singular value decomposition of $\mb{A}$ and $\mb{I}_k$ is the identity matrix in $\R^{k\times k}$ (see for example, Watson \cite{Watson93}). In this particular case, the unique subgradient of the Ky Fan $k$-norm provides the information of singular vectors corresponding to the $k$ largest singular values. Having said that, it does not preserve the information of singular values. When $\theta=0$, the proposed formulation with the Ky Fan $k$-norm will not return the rank-$k$ approximation of the matrix $\mb{A}$ as the Ky Fan $2$-$k$-norm does.} In the next section, we shall study the recovery of these submatrices under the presence of random noise.
\section{Recovery with Block Diagonal Matrices and Random Noise}
\label{sec:recovery}
We consider $\mb{A}=\mb{B}+\mb{R}$, where $\mb{B}$ is a block
diagonal matrix,  each block having rank one,
while $\mb{R}$ is a noise
matrix.  The main theorem shows that under certain assumptions concerning
the noise, the positions of the blocks can be recovered
from the solution of
\refs{eq:nprob}.
As mentioned in the introduction, this corresponds to solving a special
case of the approximate NMF problem, that is, a factorization $\mb{A}\approx \mb{W}\mb{H}^T$,
where $\mb{W}$ and $\mb{H}$ are nonnegative matrices.  The special case solved is that 
$\mb{W}$ and $\mb{H}$ each
consist of nonnegative columns with nonzeros in disjoint positions (so that
$\mb{A}$ is approximately a matrix with disjoint blocks each of rank one).
Even this special case of NMF is NP-hard unless further
restrictions are placed on the data model given the fact that the (exact) LAROS problem is NP-hard (see \cite{Doan10} for details).

Before starting the proof of the theorem, we need to consider some properties of subgaussian random variables. A random variable
$x$ is {\em $b$-subgaussian} if $\mathbb{E}[x]=0$ and there exists a $b>0$ such that for all $t\in\R$,
\be
\mathbb{E}\left[e^{tx}\right]\le e^{\frac{b^2t^2}{2}}.
\label{eq:subgau}
\ee
We can apply the Markov inequality for the $b$-subgaussian random variable $x$ and obtain the following inequalities:
\begin{equation}
\mathbb{P}(x\ge t)\le \exp(-t^2/(2b^2))\mbox{ and } \mathbb{P}(x\le -t)\le \exp(-t^2/(2b^2)),\quad\forall\,t>0.
\label{eq:tsubgau}
\end{equation}
The next three lemmas, which show several properties of random matrices and vectors with independent subgaussian entries, are adopted from Doan and Vavasis \cite{Doan10} and references therein.
\begin{lemma}
\label{lem:indsub}
Let $x_1,\ldots,x_k$ be independent $b$-subgaussian random variables and let $a_1,\ldots,a_k$ be scalars that satisfy $\displaystyle\sum_{i=1}^ka_i^2=1$. Then $x = \displaystyle\sum_{i=1}^ka_ix_i$ is a $b$-subgaussian random variable.
\end{lemma}
 
\begin{lemma}
\label{lem:2norm}
Let $\mb{B}\in\R^{m\times n}$ be a random matrix, where $b_{ij}$ are independent $b$-subgaussian random variables for all $i=1,\ldots,m$, and $j=1,\ldots,n$. 
Then for any $u>0$,
$$
\mathbb{P}\left(\norm{\mb{B}}\geq u\right)\leq \exp\left(-\left(\frac{8u^2}{81b^2}-(\log 7)(m+n)\right)\right).
$$
\begin{comment}
\item[(ii)]
 $\displaystyle\mathbb{P}\left(\norm{\mb{C}\mb{B}}\geq u\right)\leq
\exp\left(-\left(\frac{8u^2}{81b^2\norm{\vect{C}}^2}-(\log 7)
(m+n)\right)\right)$, where $\mb{C}$ is a deterministic matrix.
\item[(iii)] $\displaystyle\mathbb{P}\left(\norm{\mb{B}}_2\geq u\right)\leq \exp\left(-\left(\frac{2u^2}{9b^2}-(\log 7)m\right)\right)$ if $n=1$,i.e., $\mb{B}=\mb{b}$ is a random vector in $\R^m$.
\een
\end{comment}
\end{lemma}

\begin{lemma}
\label{lem:subgauinnerprod}
Let $\mb{x},\mb{y}$ be two vectors in $\R^n$ with i.i.d.~$b$-subgaussian
entries.  Then for any $t>0$,
$$
\mathbb{P}\left(\bx^T\by\geq t\right)\leq\exp\left(-\min\left\{\frac{t^2}{(4eb^2)^2n},\frac{t}{4eb^2}\right\}\right),\mbox{and }\mathbb{P}\left(\bx^T\by\leq -t\right)\leq\exp\left(-\min\left\{\frac{t^2}{(4eb^2)^2n},\frac{t}{4eb^2}\right\}\right).
$$
\end{lemma}

With these properties of subgaussian variables presented, we are now
able to state and prove the main theorem, which gives sufficient conditions
for 
optimization problem \refs{eq:nprob} to recover $k$ blocks in the presence of
noise.

\begin{theorem} 
\label{thm:recovery}
Suppose $\mb{A}=\mb{B}+\mb{R}$, where $\mb{B}$ is a block
diagonal matrix with $k_0$ blocks, that is,
$\mb{B}=\rm{diag}(\mb{B}_1,\ldots,\mb{B}_{k_0})$, where
$\mb{B}_i=\bar{\sigma}_i\bar{\mb{u}}_i\bar{\mb{v}}_i^T$,
$\bar{\mb{u}}_i\in\R^{m_i}$, $\bar{\mb{v}}_i\in\R^{n_i}$,
$\norm{\bar{\mb{u}}_i}_2=\norm{\bar{\mb{v}}_i}_2=1$,
$\bar{\mb{u}}_i>\mb{0}$, $\bar{\mb{v}}_i>\mb{0}$ for all $i=1,\ldots,k_0$.
Assume the blocks are ordered so that
$\bar{\sigma_1}\ge\bar{\sigma_2}\ge \cdots\ge\bar{\sigma}_{k_0}> 0$. 
Matrix $\mb{R}$ is a random matrix composed
of blocks in which each entry is a translated $b$-subgaussian
variable, i.e., there exists $\mu_{ij}\ge 0$ such that elements of the
matrix block $\mb{R}_{ij}/(\phi_i\phi_j)^{1/2}-\mu_{ij}\mb{e}_{m_i}\mb{e}_{n_j}^T$ are
independent $b$-subgaussian random variables for all
$i,j=1,\ldots,k_0$. Here $\phi_i=\bar{\sigma}_i/\sqrt{m_in_i}$, $i=1,\ldots,k_0$, is a scaling
factor to match the scale of $\mb{R}_{ij}$ with that of $\mb{B}_i$
and $\mb{B}_j$, and $\mb{e}_m$ denotes the $m$-vector of all $1$'s.

We define the following positive scalars that control the degree of
heterogeneity among the first $k$ blocks:
\begin{eqnarray}
\delta_u &\le& \min_{i=1,\ldots,k}\norm{\bar{\mb{u}}_i}_1/\sqrt{m_i}, 
\label{eq:deltau_def}\\
\delta_v &\le& \min_{i=1,\ldots,k}\norm{\bar{\mb{v}}_i}_1/\sqrt{n_i}, 
\label{eq:deltav_def}\\
\xi_u &\le& \min_{i=1,\ldots,k}\left(\min_{j=1,\ldots,m_i}\bar{u}_{i,j}\right)\sqrt{m_i}, 
\label{eq:xiu_def}\\
\xi_v &\le& \min_{i=1,\ldots,k}\left(\min_{j=1,\ldots,n_i}\bar{v}_{i,j}\right)\sqrt{n_i},
\label{eq:xiv_def}\\
\pi_u &\ge& \max_{i=1,\ldots,k}\left(\max_{j=1,\ldots,m_i}\bar{u}_{i,j}\right)\sqrt{m_i}, 
\label{eq:piu_def}\\
\pi_v &\ge& \max_{i=1,\ldots,k}\left(\max_{j=1,\ldots,n_i}\bar{v}_{i,j}\right)\sqrt{n_i}, 
\label{eq:piv_def}\\
\rho_m & \ge & \max_{i,j=1,\ldots,k}m_i/m_j, 
\label{eq:rhom_def}\\
\rho_n & \ge & \max_{i,j=1,\ldots,k}n_i/n_j,
\label{eq:rhon_def}\\
\rho_\sigma &\ge & \bar\sigma_1 / \bar\sigma_k. 
\label{eq:rhosig_def}
\end{eqnarray}

We also assume that the blocks do not diverge much from being square;
more precisely we assume that $m_i\le O(n_j^2)$ and $n_i\le O(m_j^2)$ for
$i,j=1,\ldots,k$.
Let $\mb{p}=(k,\delta_u,\delta_v,\pi_u,\pi_v,\xi_u,\xi_v,\rho_\sigma,\rho_m,\rho_n)$
denote the vector of parameters controlling the heterogeneity.

For the remaining noise blocks $i=k+1,\ldots,k_0$, we assume that 
their dominant singular values are substantially smaller than those
of the first $k$:
\begin{equation}
\bar\sigma_{k+1}\le \frac{0.23 \bar\sigma_k}{k+1},
\label{eq:sigmanoisebd}
\end{equation}
that
their scale is bounded:
\crcount{noisephi}
\begin{equation}
\phi_i\le c_{\thenoisephi} \phi_j,
\label{eq:noisephi}
\end{equation}
for all $i=k+1,\ldots,k_0$ and $j=1,\ldots,k$,
where $c_{\thenoisephi}$ is a constant,
and that their size is bounded:
\crcount{miniNoise}
\begin{equation}
\sum_{i=k+1}^{k_0} (m_i+n_i)
\le 
c_{\theminiNoise}(\mb{p},c_{\thenoisephi},b) \min_{i=1,\ldots k} m_in_i,
\label{eq:mnp1bd}
\end{equation}
where $c_{\theminiNoise}(\mb{p},c_{\thenoisephi},b)$ is given by \refs{eq:c0def} below.
Assume that 
\crcount{mubd}
\begin{equation}
\mu_{ij}\le c_{\themubd}(\mb{p},c_{\thenoisephi}),
\label{eq:mubd}
\end{equation}
for all $i,j=1,\ldots,k_0$, where $c_{\themubd}(\mb{p},c_{\thenoisephi})$ 
is given by \refs{eq:c1def} below.
Then provided that
\crcount{thetarange}
\begin{equation}
c_{\thethetarange}(\mb{p})\left(\sum_{i=1}^k m_in_i\right)^{-1/2} \le \theta \le
2c_{\thethetarange}(\mb{p})\left(\sum_{i=1}^k m_in_i\right)^{-1/2},
\label{eq:thetarange}
\end{equation}
where $c_{\thethetarange}(\mb{p})$ is given by \refs{eq:cthetarangedef} below,
the optimization problem \refs{eq:nprob} will return $\mb{X}$ with nonzero
entries precisely in the positions of $\mb{B_1},\ldots,\mb{B_k}$ with
probability exponentially close to $1$ as $m_i,n_i\rightarrow \infty$ for all
$i=1,\ldots,k_0$.
\end{theorem}

\noindent
{\bf Remarks.}
\begin{enumerate}
\item
Note that the theorem does not recover the exact values of $(\bar\sigma_i,
\bar{\mb{u}}_i,\bar{\mb{v}}_i)$; it is clear that this is impossible in
general under the assumptions made.
\item
The theorem is valid under arbitrary permutation of the rows and
columns (i.e., the block structure may be `concealed') since
\refs{eq:nprob} is invariant under such transformations.
\item
Given the fact that for all $i=1,\ldots,k$,
$$
0<\left(\min_{j=1,\ldots,m_i}\bar{u}_{i,j}\right)\sqrt{m_i}\leq \norm{\bar{\mb{u}}_i}_1/\sqrt{m_i}\leq 1\leq \left(\max_{j=1,\ldots,m_i}\bar{u}_{i,j}\right)\sqrt{m_i},
$$
we can always choose $\xi_u$, $\delta_u$, and $\pi_u$ such that $0<\xi_{u}\leq \delta_{u}\leq 1\leq \pi_{u}$. Similarly, we assume $0<\xi_{v}\leq \delta_{v}\leq 1\leq \pi_{v}$. 
These parameters measure how much $\bar{\mb{u}}_i$
and $\bar{\mb{v}}_i$ diverge from $\mb{e}_{m_i}$ and $\mb{e}_{n_i}$ after
normalization respectively.  
The best case for our theory (i.e., the least restrictive
values of parameters) occurs when all of these scalars are equal to 1.
Similarly
$\rho_\sigma,\rho_m,\rho_n\geq 1$, and the best case for the
theory is when they are all equal to 1.

\item
It is an implicit assumption of
the theorem that the scalars contained in $\mb{p}$ as well as $b$, which
controls the subgaussian random variables, stay fixed
as $m_i,n_i\rightarrow\infty$.

{
\item As compared to the recovery result in Ames \cite{Ames13} for the planted $k$-biclique problem, our result for the general bicluster problem is in general weaker in terms of noise magnitude (as compared to data magnitude) but stronger in terms of block sizes. Ames \cite{Ames13} requires $m_i=\tau_i^2 n_i$, where $\tau_i$ are scalars for all $i$, $i=1,\ldots,k+1,$ whereas we only need $m_i\leq O(n_j^2)$ and $n_i\leq O(m_j^2)$ for $i,j=1,\ldots,k$. More importantly, the noise block size, $n_{k+1}$, is more restricted as compared to data block sizes, $n_i$, for $i=1,\ldots,k$, in Ames \cite{Ames13} with the condition
$$
c_1\left(\sqrt{k}+\sqrt{n_{k+1}}+1\right)\sqrt{\sum_{i=1}^{k+1}n_i} + \beta\tau_{k+1}n_{k+1}\leq c_2\gamma\min_{i=1,\ldots,k}n_i.
$$
In contrast, for our recovery result, \refs{eq:mnp1bd} means that the total size of the noise
blocks can be much larger (approximately the square) than the size of the
data blocks.  Thus, the theorem shows that the $k$ blocks can be found
even though they are hidden in a much larger matrix. In the special case when $k_0=k+1$, $m_i=n_i=n$, $\bar{\sigma}_i=\bar{\sigma}$ for all $i=1,\ldots,k$, and $m_{k+1}=n_{k+1}$, combining \refs{eq:noisephi} and \refs{eq:mnp1bd}, we will obtain the following condition, which clearly shows the relationship between block sizes:
$$
\frac{\bar{\sigma}_{k+1}}{c_0\bar{\sigma}}n\leq n_{k+1}\leq \frac{c_1(\mb{p},c_0,b)}{2}n^2.
$$}
{
\item As compared to the recovery result in Doan and Vavasis \cite{Doan10} when $k=1$, our recovery result is for a more general setting with $\bar{\sigma}_2>0$ instead of $\bar{\sigma}_2=0$ as in Doan and Vavasis \cite{Doan10}. We therefore need additional conditions on $\bar{\sigma}_i$, $i=1,2$. In addition, we need to consider the off-diagonal blocks $(i,j)$ for $i,j=1,\ldots,k$, which is not needed when $k=1$. This leads to more (stringent) conditions on the noise magnitudes. Having said that, the conditions on the parameter $\theta$ and block sizes remain similar. We still require $\theta$ to be in the order of $(m_1n_1)^{-1/2}$ as in Doan and Vavasis \cite{Doan10}. The conditions $m_1\leq O(n_1^2)$ and $n_1\leq O(m_1^2)$ are similar to the condition $m_1n_1\geq\Omega((m_1+n_1)^{4/3})$ in Doan and Vavasis \cite{Doan10}. Finally, the condition $m_2+n_2\leq c_1(\mb{p},c_0,b)m_1n_1$ is close to the condition $m_1n_1\geq \Omega(m_1+m_2+n_1+n_2)$, which again shows the similarity of these recovery results in terms of block sizes.}
\end{enumerate}

In order to simplify the proof, we first consolidate all blocks $i=k+1,\ldots,k_0$ into a single block and call it block $(k+1)$ of size $\bar{m}_{k+1}\times\bar{n}_{k+1}$ where
$\displaystyle\bar{m}_{k+1}=\sum_{i=k+1}^{k_0}m_i$ and 
$\displaystyle \bar{n}_{k+1}=\sum_{i=k+1}^{k_0}n_i$
The only difference is that the new block 
$\bar{\mb{B}}_{k+1,k+1}\in\mathbb{R}^{\bar{m}_{k+1}\times\bar{n}_{k+1}}$ is now 
a block diagonal matrix with $k_0-k$ blocks instead of a rank-one block.
Similarly, 
new blocks $\bar{\mb{R}}_{i,k+1}$ and $\bar{\mb{R}}_{k+1,i}$, $i=1,\ldots,k_0$, now have more than one subblock with different parameters $\mu$ instead of a single one. This new block structure helps us derive the optimality conditions more concisely. Clearly, we would like to achieve the optimal solution $\bX$ with the following structure
$$
\bX= \begin{pmatrix}\sigma_1\mb{u}_1\mb{v}_1^T & \mb{0} & \cdots  & \cdots& \mb{0}\\
\mb{0} & \ddots & \ddots & & \mb{0} \\
\vdots & \ddots & \ddots & \mb{0}& \vdots \\
\vdots & & \mb{0}& \sigma_k\mb{u}_k\mb{v}_k^T & \mb{0}\\
\mb{0} & \cdots & \cdots & \mb{0} & \mb{0} \end{pmatrix},
$$
where $\norm{\mb{u}_i}_2=\norm{\mb{v}_i}_2=1$ for $i=1,\ldots,k$. Padding appropriate zeros to $\mb{u}_i$ and $\mb{v}_i$ to construct $\mb{u}_i^0\in\R^m_+$ and $\mb{v}_i^0\in\R^n_+$ for $i=1,\ldots,k$, we obtain sufficient optimality conditions based on Proposition \ref{prop:noptimality} as follows:
\begin{quote}
There exist $\bY$ and $\bZ$ such that $\bY+\bZ=\mb{A}$ and
$$
\bY=\tnorm{\mb{A}}_{k,2,\theta}^\star\left[\sum_{i=1}^k\sigma_i\mb{u}_i^0(\mb{v}_i^0)^T+\mb{W}\right],\quad \bZ=\theta\tnorm{\mb{A}}_{k,2,\theta}^\star\mb{V},
$$		
where $\sigma_i>0$ for $i=1,\ldots,k$, $\displaystyle\sum_{i=1}^k\sigma_i^2=1$, $\displaystyle\norm{\mb{W}}\leq\min_{i=1,\ldots,k}\{\sigma_i\}$, $\mb{W}\mb{v}_i^0=\mb{0}$, $\mb{W}^T\mb{u}_i^0=\mb{0}$, for $i=1,\ldots,k$, and $\norm{\mb{V}}_{\infty}\leq 1$, $\mb{V}_{ii}=\mb{e}_{m_i}\mb{e}_{n_i}^T$, for $i=1,\dots,k$.
\end{quote}
Since $\mb{A}$ has the block structure, we can break these optimality conditions into appropriate conditions for each block. Starting with diagonal $(i,i)$ blocks, $i=1,\ldots,k$, the detailed conditions are:
\begin{eqnarray}
\sigma_i\mb{u}_i\mb{v}_i^T + \mb{W}_{ii}&=&\lambda(\bar{\sigma}_i\bar{\mb{u}}_i\bar{\mb{v}}_i^T+\mb{R}_{ii})-\theta\mb{e}_{m_i}\mb{e}_{n_i}^T,\label{eq:biicond1}\\
\mb{W}_{ii}^T\mb{u}_i&=&\mb{0},\label{eq:biicond2}\\
\mb{W}_{ii}\mb{v}_i&=&\mb{0},\label{eq:biicond3}
\end{eqnarray}
where $\lambda=1/\tnorm{\mb{A}}_{k,2,\theta}^\star$. For non-diagonal $(i,j)$ blocks, $i\neq j$ and $i,j=1,\ldots,k$, we obtain the following conditions:
\begin{eqnarray}
\mb{W}_{ij}+\theta\mb{V}_{ij}&=&\lambda\mb{R}_{ij},\label{eq:bijcond1}\\
\mb{W}_{ij}^T\mb{u}_i&=&\mb{0},\label{eq:bijcond2}\\
\mb{W}_{ij}\mb{v}_j&=&\mb{0},\label{eq:bijcond3}\\
\norm{\mb{V}_{ij}}_{\infty}&\leq &1.\label{eq:bijcond4}
\end{eqnarray}
For $(i,k+1)$ blocks, $i=1,\ldots,k$, we have:
\begin{eqnarray}
\mb{W}_{i,k+1}+\theta\mb{V}_{i,k+1}&=&\lambda\bar{\mb{R}}_{i,k+1},\label{eq:bik1cond1}\\
\mb{W}_{i,k+1}^T\mb{u}_i&=&\mb{0},\label{eq:bik1cond2}\\
\norm{\mb{V}_{i,k+1}}_{\infty}&\leq &1.\label{eq:bik1cond3}
\end{eqnarray}
Similarly, for $(k+1,j)$ blocks, $j=1,\ldots,k$, the conditions are:
\begin{eqnarray}
\mb{W}_{k+1,j}+\theta\mb{V}_{k+1,j}&=&\lambda\bar{\mb{R}}_{k+1,j},\label{eq:bk1jcond1}\\
\mb{W}_{k+1,j}\mb{v}_j&=&\mb{0},\label{eq:bk1jcond2}\\
\norm{\mb{V}_{k+1,j}}_{\infty}&\leq &1.\label{eq:bk1jcond3}
\end{eqnarray}
Finally, the $(k+1,k+1)$ block needs the following conditions:
\begin{eqnarray}
\mb{W}_{k+1,k+1}+\theta\mb{V}_{k+1,k+1}&=&\lambda\left(\bar{\mb{B}}_{k+1,k+1}+\bar{\mb{R}}_{k+1,k+1}\right),\label{eq:bk1k1cond1}\\
\norm{\mb{V}_{k+1,k+1}}_{\infty}&\leq &1.\label{eq:bk1k1cond2}
\end{eqnarray}
The remaining conditions are not block separable. We still need $\sigma_i>0$, $i=1,\ldots,k$, and $\displaystyle\sum_{i=1}^k\sigma_i^2=1$. The last condition, which is $\displaystyle\norm{\mb{W}}\leq\min_{i=1,\ldots,k}\{\sigma_i\}$, can be replaced by the following sufficient conditions that are block separable by applying the fact that $\norm{\mb{W}}^2\leq\sum_{i,j}\norm{\mb{W}_{ij}}^2$:
\be
\label{eq:bijnormcond}
\norm{\mb{W}_{ij}}\leq\frac{1}{k+1}\min_{i=1,\ldots,k}\{\sigma_i\},\quad\,i,j=1,\ldots,k+1.
\ee

With these sufficient block separable conditions, in order to construct $(\mb{V} ,\mb{W})$, we now need to construct $(\mb{V}_{ij},\mb{W}_{ij})$ for different pairs $(i,j)$ block by block. The block by block details are shown
in the following analysis.

{In the following proof, we assume that the random matrix $\mb{R}$ is
chosen in stages: the diagonal blocks $\mb{R}_{ii}$, $i=1,\ldots,k$, are selected
before the off-diagonal blocks.  This allows us to treat the diagonal blocks as
deterministic during the analysis of the off-diagonal blocks.  This technique of
staging independent random variables is by now standard in the literature;  see
e.g., the ``golfing'' analysis of the matrix completion problem by Gross \cite{gross2011recovering}.}

\begin{subsection}{Analysis for block $(i,i)$, $i=1,\ldots,k$}
We begin with the proof of the existence of a $\lambda>0$ that satisfies the optimality conditions. We then show the sufficient condition \refs{eq:bijnormcond} for block $(i,i)$, $i=1,\ldots,k$. The final condition that needs to be proved for these blocks is the positivity of $\mb{u}_i$ and $\mb{v}_i$, $i=1,\ldots,k$.
\begin{subsubsection}{Existence of $\mathbf{\lambda^*}$}
The conditions for $(i,i)$ block, $i=1,\ldots,k$, namely,
\refs{eq:biicond1}--\refs{eq:biicond3},
indicate that $(\sigma_i,\mb{u}_i,\mb{v}_i)$ is the dominant singular triple of $\mb{L}_{i}=\lambda(\bar{\sigma}_i\bar{\mb{u}}_i\bar{\mb{v}}_i^T+\mb{R}_{ii})-\theta\mb{e}_{m_i}\mb{e}_{n_i}^T$. They also indicate that 
\begin{equation}
\norm{\mb{W}_{ii}} = \sigma_2(\mb{L}_i)
\label{eq:normWii}
\end{equation}
since \refs{eq:biicond1}--\refs{eq:biicond3} are equivalent to the first
step of a singular value decomposition of $\mb{L}_{i}$.

For the rest of this analysis, it is more convenient notationally 
work with $\tau=\lambda/\theta$ rather than with $\lambda$ directly.
The condition $\displaystyle\sum_{i=1}^k\sigma_i^2=1$ becomes
\begin{equation}
f(\tau)=\sum_{i=1}^k\norm{\tau(\bar{\sigma}_i\bar{\mb{u}}_i\bar{\mb{v}}_i^T+\mb{R}_{ii})-\mb{e}_{m_i}\mb{e}_{n_i}^T}^2-\theta^{-2}=0.
\label{eq:oftaudef}
\end{equation}
We will prove that there exists $\tau^*>0$ such that $f(\tau^*)=0$. More precisely, we will focus our analysis of $f(\tau)$ for $\tau\in[\tau_\ell,\tau_u]$, where
$\tau_\ell$ is given by \refs{eq:tauldef} and $\tau_u$ is given by \refs{eq:tauudef}
below
and prove that
there exists $\tau^*\in[\tau_\ell,\tau_u]$ such that $f(\tau^*)=0$.

Letting $\mb{Q}_{ij}=\mb{R}_{ij}/\sqrt{\phi_i\phi_j}-\mu_{ij}\mb{e}_{m_i}\mb{e}_{n_j}^T$ for 
$i,j=1,\ldots,k$, we have: 
$\mb{Q}_{ij}$ are $b$-subgaussian random matrices with independent elements. 
The function $f$ can be rewritten as follows:
\begin{eqnarray}
f(\tau) &=&\sum_{i=1}^k\norm{\tau\bar{\sigma}_i\bar{\mb{u}}_i\bar{\mb{v}}_i^T-(1-\tau\phi_i\mu_{ii})\mb{e}_{m_i}\mb{e}_{n_i}^T+\tau\phi_i\mb{Q}_{ii}}^2-\theta^{-2}\nonumber\\
&=&\sum_{i=1}^k\norm{\mb{P}_i(\tau)+\tau\phi_i\mb{Q}_{ii}}^2-\theta^{-2}, \label{eq:ftaudef}
\end{eqnarray}
where $\mb{P}_i(\tau)=\tau\bar{\sigma}_i\bar{\mb{u}}_i\bar{\mb{v}}_i^T-(1-\tau\phi_i\mu_{ii})\mb{e}_{m_i}\mb{e}_{n_i}^T$. Applying triangle inequality, we have:
\begin{equation}
\norm{\mb{P}_i(\tau)}-\tau\phi_i\norm{\mb{Q}_{ii}}\leq\norm{\mb{P}_i(\tau)+\tau\phi_i\mb{Q}_{ii}}\leq\norm{\mb{P}_i(\tau)}+\tau\phi_i\norm{\mb{Q}_{ii}}.
\label{eq:triPQ}
\end{equation}

We start the analysis with $\norm{\mb{P}_i(\tau)}$. We first define the following function
\be
\label{eq:gfunc}
g_i(\tau;a)=\phi_i^2\tau^2 - 2a\phi_i\tau (1-\mu_{ii}\phi_i\tau)+ (1-\mu_{ii}\phi_i\tau)^2,
\ee
which is a quadratic function in $\tau$ with any fixed parameter $a$.  
Note by \refs{eq:tauudef1} below that $\tau_u\le 0.3/(\phi_i\mu_{ii})$ for all
$i=1,\ldots, k$, so $1-\mu_{ii}\phi_i\tau\ge 0$ and $\tau\ge 0$ for $\tau\in[\tau_\ell,\tau_u]$.  Therefore, provided $a\le 1$ and $\tau\in[\tau_\ell,\tau_u]$, 
\begin{eqnarray}
g_i(\tau;a)&=&\left(\phi_i\tau - (1-\mu_{ii}\phi_i\tau)\right)^2 +
2(1-a)\phi_i\tau(1-\mu_{ii}\phi_i\tau) \nonumber \\
& \ge & \left(\phi_i\tau - (1-\mu_{ii}\phi_i\tau)\right)^2 \nonumber\\
&= & \left(\phi_i(1+\mu_{ii})\tau - 1\right)^2
\label{eq:gibd1} \\
& \ge & 0 \label{eq:gibd2}.
\end{eqnarray}

We now analyze the dominant singular triple of
$\mb{P}_i(\tau)=\tau\bar{\sigma}_i\bar{\mb{u}}_i\bar{\mb{v}}_i^T
-(1-\tau\phi_i\mu_{ii})\mb{e}_{m_i}\mb{e}_{n_i}^T$
for a fixed $\tau\in[\tau_\ell,\tau_u]$. It is clear that
dominant right singular vector lies in
$\mbox{span}\{\bar{\mb{v}}_i,\mb{e}_{n_i}\}$ since this is the range
of $(\mb{P}_i(\tau))^T$. Letting $\zeta_i=\norm{\mb{P}_i(\tau)}^2$
be the square of the dominant singular value, we have: $\zeta_i$ is
a solution of the following eigenvector problem:
$$
\left(\mb{P}_i(\tau)\right)^T\mb{P}_i(\tau)(\alpha\bar{\mb{v}}_i+\beta\mb{e}_{n_i})=\zeta_i(\alpha\bar{\mb{v}}_i+\beta\mb{e}_{n_i}).
$$
Expanding and gathering multiples of $\bar{\mb{v}}_i$ and $\mb{e}_{n_i}$, we obtain the following $2\times 2$ eigenvalue problem 
\begin{equation}
\mb{M}_i\left(
\begin{array}{c}
\alpha \\
\beta
\end{array}
\right)
=
\zeta_i
\left(
\begin{array}{c}
\alpha \\
\beta
\end{array}
\right),
\label{eq:Peig}
\end{equation}
where
\begin{equation}
\mb{M}_i=
\left(
\begin{array}{cc}
\tau^2\bar{\sigma}_i^2-\tau\bar{\sigma}_ih_i(\tau)\norm{\bar{\mb{u}}_i}_1\norm{\bar{\mb{v}}_i}_1\ & \tau^2\bar{\sigma}_i^2-\tau\bar{\sigma}_ih_i(\tau)\norm{\bar{\mb{u}}_i}_1n_i \\
 (h_i(\tau))^2\norm{\bar{\mb{v}}_i}_1m_i - \tau\bar{\sigma}_ih_i(\tau)\norm{\bar{\mb{u}}_i}_1 & (h_i(\tau))^2m_in_i - \tau\bar{\sigma}_ih_i(\tau)\norm{\bar{\mb{u}}_i}_1\norm{\bar{\mb{v}}_i}_1
\end{array}
\right),
\label{eq:Mdef}
\end{equation}
and $h_i(\tau)=1-\tau\phi_i\mu_{ii}$, $i=1,\ldots,k$.
Thus, $\zeta_i$ is a root of the equation 
\begin{equation}
\zeta_i^2-\mbox{trace}(\mb{M}_i)\zeta_i+\det(\mb{M}_i)=0,
\label{eq:muquad}
\end{equation}
where 
\begin{eqnarray}
\mbox{trace}(\mb{M}_i)&=&\tau^2\bar{\sigma}_i^2-2\tau\bar{\sigma}_i(1-\tau\phi_i\mu_{ii})\norm{\bar{\mb{u}}_i}_1\norm{\bar{\mb{v}}_i}_1+(1-\tau\phi_i\mu_{ii})^2m_in_i\nonumber\\
&=&m_in_i\left[\tau^2\phi_i^2-2\tau\phi_i(1-\tau\phi_i\mu_{ii})\delta_{u,i}\delta_{v,i}+(1-\tau\phi_i\mu_{ii})^2\right]\nonumber \\
& = &
m_in_ig_i(\tau;\delta_{u,i}\delta_{v,i}).
\label{eq:trMi}
\end{eqnarray}
and
\begin{eqnarray}
\det(\mb{M}_i) &=& \tau^2\bar{\sigma}_i^2(1-\tau\phi_i\mu_{ii})^2(m_i-\norm{\bar{\mb{u}}_i}_1^2)(n_i-\norm{\bar{\mb{v}}_i}_1^2)\nonumber\\
&=&
m_i^2n_i^2\tau^2\phi_i^2(1-\tau\phi_i\mu_{ii})^2(1-\delta_{u,i}^2)(1-\delta_{v,i}^2)
\label{eq:detMformula} \\
&\geq &
0. \nonumber
\end{eqnarray}
Here, we have introduced notation
\begin{eqnarray*}
\delta_{u,i}&=&\Vert\bar{\mb{u}}_i\Vert_1/\sqrt{m_i}, \\
\delta_{v,i}&=&\Vert\bar{\mb{v}}_i\Vert_1/\sqrt{n_i},
\end{eqnarray*}
that we will continue to use for the remainder of the proof.
It is apparent that $\delta_{u,i}\in[\delta_u,1]$ by \refs{eq:deltau_def}
and similarly $\delta_{v,i}\in[\delta_v,1]$.

Let $\Delta$ be the discriminant of the quadratic equation \refs{eq:muquad}, that is, 
\begin{equation}
\Delta=\mbox{trace}(\mb{M}_i)^2-4\det(\mb{M}_i).
\label{eq:Deltadef}
\end{equation}
 We have:
\begin{eqnarray}
\Delta& = &m_in_i\left[\tau^2\phi_i^2-2\tau\phi_i(1-\tau\phi_i\mu_{ii})\left(\delta_{u,i}\delta_{v,i}+\sqrt{(1-\delta_{u,i}^2)(1-\delta_{v,i}^2)}\right)+(1-\tau\phi_i\mu_{ii})^2\right]
\nonumber\\
& & \,\mbox{}\cdot m_in_i\left[\tau^2\phi_i^2-2\tau\phi_i(1-\tau\phi_i\mu_{ii})\left(\delta_{u,i}\delta_{v,i}-\sqrt{(1-\delta_{u,i}^2)(1-\delta_{v,i}^2)}\right)+(1-\tau\phi_i\mu_{ii})^2\right] 
\nonumber \\
& = & (m_in_i)^2g_i\left(\tau;\delta_{u,i}\delta_{v,i}+\sqrt{(1-\delta_{u,i}^2)(1-\delta_{v,i}^2)}\right)\cdot g_i\left(\tau;\delta_{u,i}\delta_{v,i}-\sqrt{(1-\delta_{u,i}^2)(1-\delta_{v,i}^2)}\right).  \label{eq:Deltafac}
%& & \quad\mbox{}\cdot g_i\left(\tau;\delta_{u,i}\delta_{v,i}-\sqrt{(1-\delta_{u,i}^2)(1-\delta_{v,i}^2)}\right)
\end{eqnarray}
Note that $\displaystyle 1-\left(\delta_{u,i}\delta_{v,i}+\sqrt{(1-\delta_{u,i}^2)(1-\delta_{v,i}^2)}\right)^2 = \left(\delta_{u,i}\sqrt{1-\delta_{v,i}^2}-\delta_{v,i}\sqrt{1-\delta_{u,i}^2}\right)^2 \geq 0$.  Therefore, the second argument to each invocation of $g_i$ in the previous equation is less than or equal to 1. Since $\tau\in[\tau_\ell,\tau_u]$, 
it follows
that both evaluations of $g_i$ yield nonnegative numbers, and therefore $\Delta\ge 0$.

We next claim that 
\begin{equation}
\Delta=m_i^2n_i^2g_i(\tau;p_i(\tau))^2
\label{eq:Deltaform1}
\end{equation}
for a continuous
$p_i(\tau)\in[\delta_{u,i}\delta_{v,i},1]$ for all $\tau\in[\tau_\ell,\tau_u]$.
In other words, there exists a continuous $p_i(\tau)$ in the range 
$[a,1]$
satisfying the
equation
\begin{equation}
g_i(\tau;p_i(\tau))^2=g_i(\tau;a+c)g_i(\tau;a-c),
\label{eq:tauexp}
\end{equation}
where, for this paragraph, 
$a=\delta_{u,i}\delta_{v,i}$ and 
$c=\sqrt{(1-\delta_{u,i}^2)(1-\delta_{v,i}^2)}$.  This is proved by first
treating $p_i$ as an unknown and expanding \refs{eq:tauexp}.  After simplification,
the result is a quadratic equation for $p_i$.  The facts that $0\le a,c\le 1$
and $a+c\le 1$ allow one to argue that the quadratic equation has a sign
change over the interval $[a,1]$ for all $\tau\in[0,1/(\phi_i\mu_{ii})]$ (hence
for all $\tau\in[\tau_\ell,\tau_u]$). Thus,
the quadratic has a unique root in this interval, which may be taken to be $p_i$; 
it must vary continuously with the coefficients of the quadratic and hence with $\tau$.
The details are left to the reader.
In addition to $\tau$, $p_i(\tau)$ depends
on $\mu_{ii}$, $\phi_i$, $\delta_{u,i}$ and $\delta_{v,i}$.

Thus, by the quadratic formula applied to \refs{eq:muquad}, we can obtain $\zeta_i$ as the larger root
\be\label{eq:zdef}
\zeta_i = \frac{1}{2}(\mbox{trace}(\mb{M}_i(\tau))+ \sqrt{\Delta})=m_in_ig_i(\tau;a_i(\tau)),
\ee
where the second equation comes from adding \refs{eq:trMi} to the square root
of \refs{eq:Deltaform1}
and noting that for any $\tau,a,b$, $(g_i(\tau;a)+g_i(\tau;b))/2=
g_i(\tau;(a+b)/2)$.
Here, we have:
\be
\label{eq:aitaudef}
\displaystyle a_i(\tau)=\frac{1}{2}\left(\delta_{u,i}\delta_{v,i}+p_i(\tau)\right).
\ee
By the earlier bound on $p_i(\tau)$, this implies
$a_i(\tau)\in[\underline{a}_i,\overline{a}_i]$,
where 
\begin{equation}
\underline{a}_i=\delta_{u,i}\delta_{v,i};\quad \overline{a}_i=\frac{1}{2}+\frac{1}{2}\delta_{u,i}\delta_{v,i}.
\label{eq:aibddef}
\end{equation}
Clearly $0\le \underline{a}_i\le \overline{a}_i\le 1$ for all $i$ since 
$\delta_{u,i},\delta_{v,i}\in[0,1]$. Note that tighter bounds are possible by a more careful analysis of $\Delta$. 

Since $\zeta_i=\norm{\mb{P}_i(\tau)}^2$, we can then express $\norm{\mb{P}_i(\tau)}$ as follows:
\begin{equation}
\norm{\mb{P}_i(\tau)}=\sqrt{\zeta_i}=\sqrt{m_in_ig_i(\tau;a_i(\tau))}.
\label{eq:normpieq}
\end{equation}

Next, consider again \refs{eq:gibd1};
the right-hand side is a convex quadratic function of $\tau$ with minimizer at 
$1/(\phi_i(1+\mu_{ii}))$.   It follows from \refs{eq:tauldef1} that
$\tau_\ell\ge 2/\phi_i$ for all $i$.
Thus, for $\tau\in[\tau_\ell,\tau_u]$, we have: 
$$
\tau\geq\frac{2}{\phi_i}\geq \frac{1}{\phi_i(1/2+\mu_{ii})}>\frac{1}{\phi_i(1+\mu_{ii})}.
$$
Thus, the right-hand side of \refs{eq:gibd1} is an
increasing function of $\tau$ for $\tau\in[\tau_\ell,\tau_u]$. 
We then have, for any $\tau\in[\tau_\ell,\tau_u]$
\begin{eqnarray*}
g_i(\tau;a)& \ge& \left(\phi_i(1+\mu_{ii})\left(\frac{1}{\phi_i(1/2+\mu_{ii})}\right)-1\right)^2 \\
&=& \left(\frac{1}{1+2\mu_{ii}}\right)^2.
\end{eqnarray*}
Thus, 
\begin{equation}
\norm{\mb{P}_i(\tau)}\ge \frac{\sqrt{m_in_i}}{1+2\mu_{ii}}
\label{eq:normPlb}
\end{equation}
for any $\tau\in[\tau_\ell,\tau_u]$.
We also have a second lower bound that grows linearly with $\tau$:
\begin{eqnarray}
\sqrt{g_i(\tau;a)}&\ge& \phi_i(1+\mu_{ii})\tau-1 \\
&=&
\phi_i\tau/2+\left(\phi_i(1/2+\mu_{ii})\tau-1\right) \nonumber \\
&\ge &
\phi_i\tau/2, \label{eq:sqrtglblinear}
\end{eqnarray}
where the first inequality follows from \refs{eq:gibd1} and
the other inequality is due to the fact that $\tau\phi_i\ge 2$ as noted
above.
This implies
\begin{eqnarray}
\norm{\mb{P}_i(\tau)}&\ge&\sqrt{m_in_i}\phi_i\tau/2 \nonumber \\
& = & \bar\sigma_i\tau/2. 
\label{eq:normPlblinear}
\end{eqnarray}

Next, we combine this linear lower bound on $\norm{\mb{P}_i(\tau)}$ with an
upper bound on $\norm{\mb{Q}_{ii}}$ in order to be able to take advantage
of \refs{eq:triPQ}.
\begin{claim}
$\norm{\mb{Q}_{ij}}\le (m_in_j)^{\frac{3}{8}}$ with probability 
exponentially close to $1$ as $m_i, n_j\rightarrow \infty$ 
for all $i,j=1,\ldots,k$.
\end{claim}
To establish
the claim observe that $\mb{Q}_{ij}$ 
is random with i.i.d.\ elements that are $b$-subgaussian. 
Thus by Lemma~\ref{lem:2norm}(i),
\begin{equation}
\mathbb{P}\left(\norm{\mb{Q}_{ij}}\geq (m_in_j)^{3/8}\right)\leq 
\exp\left(-\left(\frac{(m_in_j)^{3/4}}{81b^2}-(\log 7)(m_i+n_j)\right)\right),
\label{eq:probq}
\end{equation}
where $u$ is set to be $(m_in_j)^{3/8}$. The right-hand side tends to
zero exponentially fast since $(m_in_j)^{3/4}$ asymptotically dominates $m_i+n_j$
under the assumption that $m_i^{1/4}\le O(n_j^{1/2})$ and
$n_j^{1/4}\le O(m_i^{1/2})$, which was stated as a hypothesis
in the theorem.

Then with a probability exponentially close to $1$, 
the event in \refs{eq:probq} does not happen, hence we assume 
$\norm{\mb{Q}_{ij}}\le (m_in_j)^{\frac{3}{8}}$. 
Focusing on the $i=j$ case for now,
this implies
\begin{equation}
\norm{\mb{Q}_{ii}}\leq\sqrt{m_in_i}/40,
\label{eq:mini40}
\end{equation}
 for large $m_in_i$; since the
theorem applies to the asymptotic range, we assume this inequality holds
true as well.
Combining the inequality \refs{eq:mini40} 
with \refs{eq:normPlblinear} and \refs{eq:triPQ}, we obtain
\begin{eqnarray}
\norm{\mb{P}_i(\tau) + \tau\phi_i \mb{Q}_{ii}} & = & (1+\gamma_i(\tau))
\norm{\mb{P}_i(\tau)} \nonumber \\
& = & (1+\gamma_i(\tau))
\sqrt{m_in_ig_i(\tau;a_i(\tau))}.
\label{eq:1.05}
\end{eqnarray}
In the first line, we have introduced scalar $\gamma_i(\tau)$ to stand for
a quantity in the range $[-1/20,1/20]$ that
varies continuously with $\tau$.  This notation will be used throughout
the remainder of the proof.  The second line follows from
\refs{eq:normpieq}.
Combining \refs{eq:1.05} and \refs{eq:normPlblinear}, we conclude
\begin{equation}
\norm{\mb{P}_i(\tau) + \tau\phi_i \mb{Q}_{ii}} \ge
0.47\bar\sigma_i\tau.
\label{eq:pitauphiqlb}
\end{equation}
Finally, because $\mb{P}_i(\tau) + \tau\phi_i \mb{Q}_{ii}$ is a rescaling
of the right-hand side of
\refs{eq:biicond1} by $\theta$, we conclude
that 
\begin{equation}
\sigma_i \ge
0.47\bar\sigma_i\tau\theta.
\label{eq:sigmaibarsigmailb}
\end{equation}

Applying \refs{eq:1.05} to the formulation of $f(\tau)$ in \refs{eq:ftaudef}, we have, for $\tau\in [\tau_\ell,\tau_u]$:
\begin{eqnarray*}
f(\tau) & = & 
\sum_{i=1}^k\norm{\mb{P}_i(\tau) + \tau\phi_i\mb{Q}_{ii}}^2 - \theta^{-2} \\
& = &
\sum_{i=1}^k 
m_in_i(1+\gamma_i(\tau))^2
g_i(\tau;a_i(\tau)) - \theta^{-2} \\
& = & A(\tau)\tau^2 - 2B(\tau)\tau - C(\tau).
\end{eqnarray*}
The third line is obtained by expanding the quadratic formula 
for $g_i(\tau;a_i(\tau))$, which results in
\begin{eqnarray*}
A(\tau) & = & \sum_{i=1}^k (1+\gamma_i(\tau))^2 \bar{\sigma}_i^2(1+2\mu_{ii}a_i(\tau)+\mu_{ii}^2), \\
B(\tau) & = & \sum_{i=1}^k (1+\gamma_i(\tau))^2\sqrt{m_in_i}\bar{\sigma}_i(a_i(\tau)+\mu_{ii}), \\
C(\tau) & = & \theta^{-2} - \sum_{i=1}^k(1+\gamma_i(\tau))^2m_in_i.
\end{eqnarray*}

We will now prove that there exists $\tau^*\in[\tau_\ell,\tau_u]$ such that $f(\tau^*)=0$ by applying the following lemma, which is a specific form of intermediate theorem for ``pseudo-quadratic'' functions.

\begin{lemma}
\label{lem:inter}
Consider a real-valued function $\hat f(\tau)$ of the form 
$$\hat f(\tau)=A(\tau)\tau^2 -2B(\tau)\tau - C(\tau),$$
where $A(\tau)$, $B(\tau)$, $C(\tau)$ are continuous functions
of $\tau$.  Suppose there are two triples of positive numbers
$(\underline{A},\underline{B},\underline{C})<(\overline{A},\overline{B},\overline{C})$
(where `$<$' is understood element-wise).
Define
\begin{equation}
\tau_\ell'=
\frac{\underline{B}+\sqrt{\underline{B}^2+\underline{A}\cdot\underline{C}}}{\overline{A}},
\label{eq:tauellpdef}
\end{equation}
and
\begin{equation}
\tau_u'=
\frac{\overline{B}+\sqrt{\overline{B}^2+\overline{A}\cdot\overline{C}}}{\underline{A}}.
\label{eq:tauupdef}
\end{equation}
(Clearly $\tau_\ell'<\tau_u'$.)
Suppose further that there is an interval $[\tau_\ell,\tau_u]$ such that
$\tau_\ell\le \tau_\ell'\le \tau_u'\le\tau_u$ and such that for all
$\tau\in[\tau_\ell,\tau_u]$,
$$(\underline{A},\underline{B},\underline{C})\le
(A(\tau),
B(\tau),
C(\tau)) \le 
(\overline{A},\overline{B},\overline{C}).$$
Then there exists a root $\tau^*\in [\tau_\ell',\tau_u']$
(and therefore also in  $[\tau_\ell,\tau_u]$) such that $\hat f(\tau^*)=0$.
\label{lem:intval}
\end{lemma}

\begin{pf}
Some simple algebra shows that 
$\hat f(\tau_\ell')=A(\tau_\ell')(\tau_\ell')^2 -2B(\tau_\ell')(\tau_\ell')-C(\tau_\ell')\le 0$
while
$\hat f(\tau_u')=A(\tau_u')(\tau_u')^2 -2B(\tau_u')\tau_u'-C(\tau_u')\ge 0,$ so there
is a $\tau^*\in[\tau_\ell',\tau_u']$ such that $\hat f(\tau^*)=0$ by the intermediate
value theorem.
\end{pf}

In order to apply Lemma \ref{lem:inter}, we now define the following scalars:
\begin{eqnarray*}
\overline{A} &=&(10/9)\sum_{i=1}^k \bar{\sigma}_i^2(1+2\mu_{ii}\overline{a}_i+\mu_{ii}^2), \\
\underline{A}&=& 0.90\sum_{i=1}^k \bar{\sigma}_i^2(1+2\mu_{ii}\underline{a}_i+\mu_{ii}^2),\\
\overline{B}& =& (10/9)\sum_{i=1}^k \sqrt{m_in_i}\bar{\sigma}_i(\overline{a}_i+\mu_{ii}), \\
\underline{B}&=& 0.90\sum_{i=1}^k \sqrt{m_in_i}\bar{\sigma}_i(\underline{a}_i+\mu_{ii}), \\
\overline{C} &=& (10/9)\left(\theta^{-2} - \sum_{i=1}^km_in_i\right), \\
\underline{C} &=& (9/10)\left(\theta^{-2} - \sum_{i=1}^km_in_i\right).
\end{eqnarray*}
It is obvious that $(0,0)< (\underline{A},\underline{B})< (\overline{A},\overline{B})$.
It follows from \refs{eq:thetarange}, \refs{eq:tauellbd2}, and \refs{eq:cthetarangedef}
below that the
parenthesized quantity in the definitions of $\overline{C},\underline{C}$ is positive,
\crcount{tauellbd}
$$
\theta^{-2}-\sum_{i=1}^km_in_i\geq\left(\frac{1}{4c_{\thethetarange}(\mb{p})^2}-1\right)\sum_{i=1}^km_in_i=1.2^4\left(c_{\thetauellbd}(\mb{p})\right)^2\left(\frac{\rho_m\rho_nk}{k+\rho_m\rho_n-1}\right)\sum_{i=1}^km_in_i>0,
$$
and hence we also have $0<\underline{C}<\overline{C}$.

In addition, given the fact that $a_i(\tau)\in[\underline{a}_i;\overline{a}_i]$ and $\gamma_i(\tau)\in[-1/20;1/20]$ for $\tau\in[\tau_\ell,\tau_u]$, so this establishes
for this interval that
$(\underline{A},\underline{B},\underline{C})\le
(A(\tau),
B(\tau),
C(\tau)) \le 
(\overline{A},\overline{B},\overline{C}).$

We now show that $\tau_\ell'\ge \tau_\ell$ and 
$\tau_u'\le \tau_u$.  We have
$$
\tau_\ell'=\frac{\underline{B}+\sqrt{\underline{B}^2+\underline{A}\cdot\underline{C}}}{\overline{A}}\geq\frac{\sqrt{\underline{A}\cdot\underline{C}}}{\overline{A}}.
$$
Using the facts that $0\leq \overline{a}_i\leq 1$, 
$0\leq\mu_{ii}\leq c_{\themubd}(\mb{p},c_{\thenoisephi})\le 0.08$ (see \refs{eq:c1def} below), 
and $\displaystyle\theta^{-2}-\sum_{i=1}^km_in_i\geq 1.2^4\left(c_{\thetauellbd}(\mb{p})\right)^2\left(\frac{\rho_m\rho_nk}{k+\rho_m\rho_n-1}\right)\sum_{i=1}^km_in_i>0$ as above, we have:
$$
\tau_\ell'\geq \left(\frac{\rho_m\rho_nk}{k+\rho_m\rho_n-1}\right)^{1/2}c_{\thetauellbd}(\mb{p})
\left(\sum_{i=1}^km_in_i\right)^{1/2}\left(\sum_{i=1}^k\bar{\sigma}_i^2\right)^{-1/2}.
$$
Next, observe that 
$$\left(\sum_{i=1}^k\bar{\sigma}_i^2\right)^{-1/2}\ge k^{-1/2}\sigma_1^{-1}$$
while
$$\left(\sum_{i=1}^km_in_i\right)^{1/2}\ge \left(1+\frac{k-1}{\rho_m\rho_n}\right)^{1/2}(m_1n_1)^{1/2}.$$
Since $\phi_1=\sigma_1/\sqrt{m_1n_1}$, we conclude that
$$
\tau_\ell'\geq c_{\thetauellbd}(\mb{p})\phi_1^{-1}=\tau_{\ell},
$$
given the definition of $\tau_{\ell}$ in \refs{eq:tauldef}.
% where
% \begin{equation}
% c_{\thetauellbd}(\mb{p}) = 
% 0.8(c_{\thethetarange}(\mb{p})^{-2}/4-1)^{1/2}k^{-1/2}.
% \label{eq:tauellbd}
% \end{equation}
%Calling upon \refs{eq:tauldef} shows that $\tau_{\ell}'\ge\tau_{\ell}$.

We now consider the condition for $\tau_u'$. We have:
$$
\tau_u'=
\frac{\overline{B}+\sqrt{\overline{B}^2+\overline{A}\cdot\overline{C}}}{\underline{A}}=\frac{\overline{B}}{\underline{A}}+\sqrt{\frac{\overline{B}^2}{\underline{A}^2}+\frac{\overline{A}\cdot\overline{C}}{\underline{A}^2}}.
$$
Using the fact that $0\leq\overline{a}_i\leq 1$, $0\leq\mu_{ii}\leq 0.08$, 
we have 
\begin{eqnarray*}
\frac{\overline{B}}{\underline{A}} & \le & 
\frac{100}{81}\cdot\left(1.08\sum_{i=1}^k\bar{\sigma}_i\sqrt{m_in_i}\right)\left(\sum_{i=1}^k\bar{\sigma}_i^2\right)^{-1} \\
&\le & \frac{4}{3}\cdot\left(\frac{1+(k-1)\sqrt{\rho_m\rho_n}}{1+(k-1)\rho_\sigma^{-2}}\right)\phi_1^{-1},
\end{eqnarray*}
and, using also \refs{eq:thetarange},
\begin{eqnarray*}
\frac{\overline{A}\cdot\overline{C}}{\underline{A}^2} & \le & 
\frac{16}{9}\cdot\left(c_{\thethetarange}(\mb{p})^{-2}-1\right)\left(\sum_{i=1}^km_in_i\right)
\left(\sum_{i=1}^k \bar\sigma_i^2\right)^{-1} \\
&\le & 
\frac{16}{9}\cdot\left(c_{\thethetarange}(\mb{p})^{-2}-1\right)
\left(\frac{1+(k-1)\rho_m\rho_n}{1+(k-1)\rho_\sigma^{-2}}\right)\phi_1^{-2}.
\end{eqnarray*}
Note that $0<c_{\thethetarange}(\mb{p})<1$ given its definition in \refs{eq:cthetarangedef}. Now, combining these terms and 
we conclude that 
\crcount{tauubd}
$$
\tau_u'\le c_{\thetauubd}(\mb{p})\phi_1^{-1}=\tau_u,$$
given the definition of $\tau_u$ in \refs{eq:tauudef} below with $c_{\thetauubd}(\mb{p})$ defined in \refs{eq:tauucondef}.
Thus applying Lemma \ref{lem:inter}, we prove that there exists 
$\tau^*\in[\tau_\ell,
\tau_u]$ such that $f(\tau^*)=0$. This also means the existence of $\lambda^*=\theta\tau^*$. For the remainder of this proof, we will drop the asterisks and simply write these
selected values as $\tau$ and $\lambda$.

{Since $\tnorm{\mb{A}}_{k,2,\theta}^\star=1/\lambda$, the
$\tnorm{\cdot}_{k,2,\theta}^\star$-norm of $\mb{A}$ is already determined at this
step of the proof even though the random variables $\mb{R}$ for the off-diagonal
blocks of $\mb{A}$ are not yet chosen.  (Recall that we are assuming 
for the purpose of this analysis that the random
variables are staged, and that the diagonal-block random variables are chosen
before the off-diagonal blocks.)  It should not be surprising that the norm can
be determined even before all entries are chosen; for many norms such as the
vector $\infty$-norm, it is possible to make small perturbations to many coordinate
entries without affecting the value of the norm.}
\end{subsubsection}
\begin{subsubsection}{Upper bound on $\norm{{W}_{ii}}$}
Now consider the condition 
\refs{eq:bijnormcond}  for block $(i,i)$, $i=1,\ldots,k$.  By 
\refs{eq:normWii}, it suffices to show
\be
\label{eq:biinormcond}
\sigma_2(\mb{P}_i(\tau)+\tau\phi_i\mb{Q}_{ii})
\leq\frac{1}{k+1}
\sigma_1(\mb{P}_j(\tau)+\tau\phi_j\mb{Q}_{jj})
\ee
for all $j=1,\ldots,k$.
In order to analyze $\sigma_2(\mb{P}_i(\tau)+\tau\phi_i\mb{Q}_{ii})$, we start with $\sigma_2(\mb{P}_i(\tau))$. Since $\mb{P}_i(\tau)$ has the rank of at most two, $\bar{\zeta}_i=\sigma_2^2(\mb{P}_i(\tau))$ can be computed as 
the smaller root of the quadratic equation \refs{eq:muquad}, $i=1,\ldots,k$. 
Using the fact that $\zeta_i\bar{\zeta}_i=\mbox{det}(\mb{M}_i)$, we have
\begin{equation}
\bar\zeta_i =
\frac{m_in_i\tau^2\phi_i^2(1-\tau\phi_i\mu_{ii})^2(1-\delta_{u,i}^2)(1-\delta_{v,i}^2)}
{g_i(\tau;a_i(\tau))}
\label{eq:sigma2Pi}
\end{equation}
from \refs{eq:detMformula} and \refs{eq:zdef}.

Now we note that from standard singular value perturbation theory (see, for example, Theorem 7.4.51 from Horn and Johnson \cite{Horn90}) that
\begin{eqnarray*}
\sigma_2(\mb{P}_i+\tau\phi_i\mb{Q}_{ii}) &\le &
\sigma_2(\mb{P}_i)+\tau\phi_i\norm{\mb{Q}_{ii}} \\
&\le & 
\left(\frac{m_in_i\tau^2\phi_i^2(1-\tau\phi_i\mu_{ii})^2(1-\delta_{u,i}^2)(1-\delta_{v,i}^2)}
{g_i(\tau;a_i(\tau))}\right)^{1/2} + \tau\phi_i(m_in_i)^{3/8} \\
&\equiv & T_1+T_2.
\end{eqnarray*}
We handle the two terms separately.  Since we are interested in the asymptotic
case of $m_i,n_i\rightarrow\infty$, we will assume
\begin{equation}
(m_in_i)^{-1/8}\le\frac{1}{10(k+1)\rho_\sigma},
\label{eq:sigmas}
\end{equation}
for all $i,j=1,\ldots,k$.

First, we have:
\begin{eqnarray}
T_1 & = &
\left(\frac{m_in_i\tau^2\phi_i^2(1-\tau\phi_i\mu_{ii})^2(1-\delta_{u,i}^2)(1-\delta_{v,i}^2)}
{g_i(\tau;a_i(\tau))}\right)^{1/2} \nonumber \\
& \le & 
\left(\frac{m_in_i\tau^2\phi_i^2}
{\phi_i^2\tau^2/4}\right)^{1/2} \nonumber \\
&=&2\sqrt{m_in_i} \\
&\le & \frac{2\phi_j\tau\sqrt{m_jn_j}}{6(k+1)} \nonumber\\
& = & \frac{\tau\bar\sigma_j}{3(k+1)} \nonumber\\
& \le & \frac{\sigma_1(\mb{P}_j(\tau) + \tau\phi_j \mb{Q}_{jj})}
{3\cdot 0.47(k+1)}. \label{eq:T1ana}
\end{eqnarray}
The inequality in the second line follows from the fact that $0<1-\tau\phi_i\mu_{ii}\leq 1$ for $\tau\in[\tau_\ell,\tau_u]$ for the
numerator and \refs{eq:sqrtglblinear} for the denominator.  The inequality
in the fourth line follows from 
$\tau\phi_j\ge 6(k+1)\sqrt{\rho_m\rho_n}$, which follows
from \refs{eq:tauldef1}.
The last
line follows from \refs{eq:pitauphiqlb}. 
Next, we have:
\begin{eqnarray*}
T_2 & = & \tau\phi_i(m_in_i)^{3/8} \\
 & = & \tau\bar{\sigma}_i/(m_in_i)^{1/8} \\
 & \le & \tau\bar{\sigma}_j/(10(k+1)) \\
 & \le & \frac{\sigma_1(\mb{P}_j(\tau) + \tau\phi_j \mb{Q}_{jj})}{10\cdot 0.47(k+1)},
\end{eqnarray*}
where the third line follows from \refs{eq:sigmas} and the fourth
again from \refs{eq:pitauphiqlb}.  This inequality and 
\refs{eq:T1ana} together establish
\refs{eq:biinormcond}.
\end{subsubsection}

\begin{subsubsection}{Positivity of $\mb{u}_i$ and $\mb{v}_i$}
The final condition for the $(i,i)$ block is the positivity of
singular vectors. We will show that with high probability, the matrix
$\mb{S}_i=(\mb{P}_i(\tau)+\tau\phi_i\mb{Q}_{ii})^T(\mb{P}_i(\tau)+\tau\phi_i\mb{Q}_{ii})$
is positive, which implies the positivity of the right singular
vector. (At the end of this subsection we consider the left singular vector.)  We
have: $\mb{S}_i=\mb{S}_i^1+\mb{S}_i^2+\mb{S}_i^3+\mb{S}_i^4$, where
$\mb{S}_i^1=(\mb{P}_i(\tau))^T\mb{P}_i(\tau)$,
$\mb{S}_i^2=\tau\phi_i(\mb{P}_i(\tau))^T\mb{Q}_{ii}$,
$\mb{S}_i^3=\tau\phi_i\mb{Q}_{ii}^T\mb{P}_i(\tau)$, and
$\mb{S}_i^4=\tau^2\phi_i^2\mb{Q}_{ii}^T\mb{Q}_{ii}$. Start with
$\mb{S}_i^1$.
Recall $\delta_{u,i}=\mb{e}^T_{m_i}\bar{\mb{u}}_i/\sqrt{m_i}$
and $\bar \sigma_i=\phi_i\sqrt{m_in_i}$.
Then we have:
\begin{align}
S_i^1(l,j)=&\,(1-\tau\phi_i\mu_{ii})^2m_i\left(\sqrt{n_i}\max\{\bar{v}_{i,l},\bar{v}_{i,j}\}\psi_i\left(\sqrt{n_i}\min\{\bar{v}_{i,l},\bar{v}_{i,j}\}\psi_i-\frac{\delta_{u,i}(\bar{v}_{i,l}+\bar{v}_{i,j})}{\max\{\bar{v}_{i,l},\bar{v}_{i,j}\}}\right)+1\right)\nonumber\\
\geq &\,(1-\tau\phi_i\mu_{ii})^2m_i\left[\sqrt{n_i}\max\{\bar{v}_{i,l},\bar{v}_{i,j}\}\psi_i(\xi_{v}\psi_i-2)+1\right], \label{eq:silj}
\end{align}
where we let $\psi_i$ denote $\tau\phi_i/(1-\tau\phi_i\mu_{ii})$ for the
remainder of the analysis of the positivity condition, and where $\xi_v$ was
defined by \refs{eq:xiv_def}.

{From \refs{eq:tauldef4}, $\tau_\ell\ge 4/(\xi_v\phi_i)$.
Since $\tau\in[\tau_\ell,\tau_u]$, we have $\tau\phi_i \ge 4/\xi_v$ (and similarly, $\tau\phi_i \ge 4/\xi_u$)
for all $i=1,\ldots,k$. Thus, since $0<1-\mu_{ii}\tau\phi_i\le 1$, we also
conclude $\psi_i\ge 4/\xi_v$ and hence $\psi_i\xi_v-2\ge \psi_i\xi_v/2$.}
Substituting into
\refs{eq:silj} yields
\begin{eqnarray}
S_i^1(l,j)&\ge &
(1-\tau\phi_i\mu_{ii})^2m_i\left[\sqrt{n_i}\max\{\bar{v}_{i,l},\bar{v}_{i,j}\}
\psi_i^2\xi_{v}/2+1\right] \nonumber \\
&\ge &
(1-\tau\phi_i\mu_{ii})^2m_i(\psi_i^2\xi_{v}^2/2+1), \label{eq:Silj2}
\end{eqnarray}
for $l,j=1,\ldots,n_i$.

{Now considering the matrix $\mb{S}_i^2$, we have:
\begin{eqnarray}
S_i^2(l,j)&=&\tau\phi_i\sum_{s=1}^{m_i}\left(\tau\bar{\sigma}_i\bar{v}_{i,l}\bar{u}_{i,s}-(1-\tau\phi_i\mu_{ii})\right)Q_{ii}(s,j) \nonumber \\
& = &
\tau^2\phi_iT_1 -\tau\phi_i T_2 \label{eq:Si2T1T2}
\end{eqnarray}
where
\begin{eqnarray*}
T_1 &\equiv &\sum_{s=1}^{m_i}(\bar{\sigma}_i\bar{v}_{i,l}\bar{u}_{i,s}+\phi_i\mu_{ii})Q_{ii}(s,j),
\\
T_2 &\equiv & \sum_{s=1}^{m_i} Q_{ii}(s,j).
\end{eqnarray*}
According to Lemma~\ref{lem:indsub}, $T_1$ and $T_2$ are both
subgaussian random variables with parameters 
$b_1\equiv 
b\Vert\bar{\sigma}_i\bar{v}_{i,l}\bar{\mb{u}}_i
+\phi_i\mu_{ii}\mb{e}_{m_i}\Vert$ and $b_2\equiv b\sqrt{m_i}$ respectively.  We can
derive an upper bound on $b_1$:
\begin{eqnarray}
b_1 & = & b\Vert\bar{\sigma}_i\bar{v}_{i,l}\bar{\mb{u}}_i+\phi_i\mu_{ii}\mb{e}_{m_i}\Vert 
\nonumber \\
& \le & 
b\bar{\sigma}_i\bar{v}_{i,l}\Vert\bar{\mb{u}}_i\Vert + 
b\phi_i\mu_{ii}\Vert\mb{e}_{m_i}\Vert \nonumber 
\nonumber \\
& = & 
b\phi_i(m_in_i)^{1/2}\bar{v}_{i,l} + 
b\phi_i\mu_{ii}\sqrt{m_i}\nonumber 
\nonumber \\
& \le & b\phi_i\sqrt{m_i}(\pi_v+\mu_{ii}), \label{eq:b1_bound}
\end{eqnarray}
where the third line used the definition $\phi_i=\sigma_i/(m_in_i)^{1/2}$ and
$\Vert\bar{\mb{u}}_i\Vert=1$, while the fourth line used
\refs{eq:piv_def}.}

{Considering the $T_1$ term first, let us determine the probability that
the negative of
$\tau^2\phi_iT_1$ exceeds $1/6$ times the lower bound given by
\refs{eq:Silj2}:
\begin{eqnarray}
\mathbb{P}\left(\tau^2\phi_iT_1 \le -\frac{(1-\tau\phi_i\mu_{ii})^2m_i\psi_i^2\xi_v^2}{12}\right)
&=&
\mathbb{P}\left(T_1 \le -\frac{m_i\phi_i\xi_v^2}{12}\right)
\nonumber \\
& \le &
\exp\left(-\frac{m_i\xi_v^4}{288 b^2(\pi_v+\mu_{ii})^2}\right),
\label{eq:s2lj1}
\end{eqnarray}
where the first line is obtained by dividing both sides by
$\tau^2\phi_i$ and substituting the definition of $\psi_i$, while the
second line is from
\refs{eq:tsubgau} with 
$t=m_i\phi_i\xi_v^2/12$ and the ``$b$'' of \refs{eq:tsubgau} given
by \refs{eq:b1_bound}.}

{Now let us consider the probability that the negative of
$\tau\phi_i T^2$ exceeds the same quantity:
\begin{eqnarray}
\mathbb{P}\left(\tau\phi_i T_2 \le -\frac{(1-\tau\phi_i\mu_{ii})^2m_i\psi_i^2\xi_v^2}{12}\right)
& = &
\mathbb{P}\left(T_2 \le -\frac{m_i\tau\phi_i\xi_v^2}{12}\right)
\nonumber \\
& \le &
\mathbb{P}\left(T_2 \le -\frac{m_i\xi_v}{3}\right)
\nonumber \\
& \le &
\exp\left(-\frac{m_i\xi_v^2}{18b^2}\right),
\label{eq:s2lj2}
\end{eqnarray}
where, for the first line we again used 
$\psi_i=\tau\phi_i/(1-\tau\phi_i\mu_{ii})$, for the second
$\tau\phi_i\ge 4/\xi_v$ derived above. The third uses
\refs{eq:tsubgau} with 
$t=m_i\xi_v/3$ and the subgaussian parameter given by $b_2$ above.}

{Combining  \refs{eq:Si2T1T2}, \refs{eq:s2lj1}, and \refs{eq:s2lj2} 
via the union bound
yields
\begin{equation}
\mathbb{P}\left(S_i^2(l,j)\leq -\frac{(1-\tau\phi_i\mu_{ii})^2m_i\psi_i^2\xi_{v}^2}{6}\right)
\leq
\exp\left(-\frac{m_i\xi_v^4}{288 b^2(\pi_v+\mu_{ii})^2}\right) 
+
\exp\left(-\frac{m_i\xi_v^2}{18b^2}\right).
\label{eq:s2ij}
\end{equation}}
Note that $\mb{S}_i^3=(\mb{S}_i^2)^T$, which means the analysis is 
the same.

For the matrix $\mb{S}_i^4$, we have 
$S_i^4(l,j)=\tau^2\phi_i^2[(\mb{Q}_{ii}(:,l))^T\mb{Q}_{ii}(:,j)]$, 
where the
square-bracketed factor is 
the inner product of two independent $b$-subgausian random 
vector for all $l\neq j$. (Note that when $l=j$, $S_i^4(l,j)\ge 0$ so
there is nothing to analyze.) {We again bound the probability
that the negative of this term
exceeds $1/3$ times the lower bound given by
\refs{eq:Silj2}:
\begin{eqnarray}
\mathbb{P}\left(S_i^4(l,j)\leq -\frac{(1-\tau\phi_i\mu_{ii})^2m_i\psi_i^2\xi_{v}^2}{6}\right)
& = &
\mathbb{P}\left((\mb{Q}_{ii}(:,l))^T\mb{Q}_{ii}(:,j)\leq -\frac{m_i\xi_{v}^2}{6}\right)
\nonumber \\
& \le &
\exp\left(-m_i\cdot \min\left(\frac{\xi_{v}^4}{576e^2b^4},
\frac{\xi_{v}^2}{24eb^2}\right)\right).
\label{eq:S4ilj2}
\end{eqnarray}
where, for the second line, we applied
Lemma~\ref{lem:subgauinnerprod} with $t=m_i\xi_v^2/6$ and
$n=m_i$.}
Combining \refs{eq:Silj2}, \refs{eq:s2ij}, and \refs{eq:S4ilj2}, we have:
\begin{eqnarray}
\mathbb{P}\left(\min_{l,j}S(l,j)\leq 0\right)&\leq &
n_i(n_i-1)\cdot\left[
\exp\left(\frac{-m_i\xi_{v}^4}{288b^2(\pi_{v}^2+1)}\right)
\right. \nonumber
\\
& & 
\left.+(1/2)
\exp\left(-m_i\cdot \min\left(\frac{\xi_{v}^4}{576e^2b^4},
\frac{\xi_{v}^2}{24eb^2}\right)\right)\right]. \label{eq:Sunion}
\end{eqnarray}

For the left singular vector, define
the matrix,
$$\mb{T}_i=(\mb{P}_i(\tau)+\tau\phi_i\mb{Q}_{ii})(\mb{P}_i(\tau)+\tau\phi_i\mb{Q}_{ii})^T.$$
The analogous analysis 
(i.e., writing $\mb{T}_i=\mb{T}_i^1+\mb{T}_i^2+\mb{T}_i^3+\mb{T}_i^4$ as above
and analyzing the four terms separately)
yields,
\begin{eqnarray}
\mathbb{P}\left(\min_{l,j}T(l,j)\leq 0\right)&\leq &
m_i(m_i-1)\cdot\left[
\exp\left(\frac{-n_i\xi_{u}^4}{288b^2(\pi_{u}^2+1)}\right)
\right. \nonumber
\\
& & 
\left.+(1/2)
\exp\left(-n_i\cdot \min\left(\frac{\xi_{u}^4}{576^2e^2b^4},
\frac{\xi_{u}^2}{24eb^2}\right)\right)\right]. \label{eq:Tunion}
\end{eqnarray}
\end{subsubsection}
\end{subsection}
\begin{subsection}{Analysis for block $(i,j)$, $i\neq j$, $i,j=1,\ldots,k$}
We now consider the off-diagonal $(i,j)$ block, $i\neq j$, $i,j=1,\ldots,k$.
Recall our notation: $\mb{u}_i$, $\mb{v}_i$ stand for the unit-norm dominant left
and right singular vectors respectively
of the right-hand side of \refs{eq:biicond1}, or, equivalently, of
$\mb{P}_i(\tau) + \tau\phi_i\mb{Q}_{ii}$.  

Let us consider the following construction
$$
\mb{V}_{ij}=\tau\left(\frac{\mb{e}_{m_i}(\mb{u}_i^T\mb{R}_{ij})}{\norm{\mb{u}_i}_1}
+
\frac{(\mb{R}_{ij}\mb{v}_j)\mb{e}_{n_j}^T}{\norm{\mb{v}_j}_1}-\frac{\mb{u}_i^T\mb{R}_{ij}\mb{v}_j}{\norm{\mb{u}_i}_1\norm{\mb{v}_j}_1}\mb{e}_{m_i}\mb{e}_{n_j}^T\right).
$$

The matrix $\mb{W}_{ij}=\lambda\mb{R}_{ij}-\theta\mb{V}_{ij}$ clearly
satisfies two orthogonal requirements, \refs{eq:bijcond2} and
\refs{eq:bijcond3}. We now just need to find the conditions so that
$\displaystyle\norm{\mb{W}_{ij}}\leq\frac{1}{k+1}\min_{i=1,\ldots,k}\sigma_i$
and $\norm{\mb{V}_{ij}}_{\infty}\leq 1$. 
\begin{subsubsection}{Upper bound on ${\norm{{\mb{V}}_{ij}}_\infty}$}
We have:
$$
\card{V_{ij}(s,t)}\leq
\tau\left(\card{\frac{\mb{u}_i^T\mb{R}_{ij}(:,t)}{\norm{\mb{u}_i}_1}}
+\card{\frac{\mb{R}_{ij}(s,:)\mb{v}_j}{\norm{\mb{v}_j}_1}}
+\card{\frac{\mb{u}_i^T\mb{R}_{ij}\mb{v}_j}{\norm{\mb{u}_i}_1\norm{\mb{v}_j}_1}}\right).
$$ 
In order to show $\norm{\mb{V}_{ij}}_{\infty}\leq 1$ with high
probability, we will show the sufficient condition that all
probabilities,
\begin{eqnarray}
&\displaystyle\mathbb{P}\left(\tau\card{\frac{\mb{u}_i^T\mb{R}_{ij}(:,t)}
{\norm{\mb{u}_i}_1}}>\frac{1}{3}\right),&\label{eq:rij1} \\
&\displaystyle\mathbb{P}\left(\tau\card{\frac{\mb{R}_{ij}(s,:)\mb{v}_j}
{\norm{\mb{v}_j}_1}}>\frac{1}{3}\right),& \label{eq:rij2} \\
&\displaystyle\mathbb{P}\left(\tau\card{\frac{\mb{u}_i^T\mb{R}_{ij}\mb{v}_j}
{\norm{\mb{u}_i}_1\norm{\mb{v}_j}_1}}>\frac{1}{3}\right),& \label{eq:rij3}
\end{eqnarray}
are exponentially small. 

Since $\tau\in[\tau_\ell,\tau_u]$, we have 
$\tau\sqrt{\phi_i\phi_j} \le 0.3/\mu_{ij}$
by \refs{eq:tauudef3}.
Thus we have:
$$
\mathbb{P}\left(\tau\card{\frac{\mb{R}_{ij}(s,:)\mb{v}_j}{\norm{\mb{v}_j}_1}}>\frac{1}{3}\right)\leq\mathbb{P}\left(\tau\card{\frac{(\mb{R}_{ij}(s,:)-\mu_{ij}\sqrt{\phi_i\phi_j}\mb{e}_{n_j}^T)\mb{v}_j}{\norm{\mb{v}_j}_1}}>\frac{1}{30}\right).
$$
Thus, to analyze \refs{eq:rij2}, it suffices to show that the probability on the
right-hand side of the preceding inequality is exponentially small.
Since $\norm{\mb{v}_j}=1$, 
$((\phi_i\phi_j)^{-1/2}\mb{R}_{ij}(s,:)-\mu_{ij}\mb{e}_{n_j}^T)\mb{v}_j$ 
is a $b$-subgaussian random variable by Lemma~\ref{lem:indsub}. {(Note that $\mb{v}_j$ depends on the 
$(j,j)$ diagonal block of $\mb{A}$, which in turn depends on $\mb{R}_{jj}$
and hence is random.
However, recall also that we have assumed that the random variables in
the block diagonals of $\mb{R}$ are chosen before the off-diagonal blocks,
so that $\mb{v}_j$ may be considered as a deterministic quantity when analyzing
$\mb{R}_{ij}$.)}

By \req{eq:tsubgau}, we have:
\be
\label{eq:vijcond2}
\mathbb{P}\left(\tau\card{\frac{(\mb{R}_{ij}(s,:)-\mu_{ij}\sqrt{\phi_i\phi_j}\mb{e}_{n_j}^T)\mb{v}_j}{\norm{\mb{v}_j}_1}}>\frac{1}{30}\right)
\leq 2\exp\left(-\frac{0.1^2\norm{\mb{v}_j}_1^2}{18b^2\tau^2\phi_i\phi_j}\right).
%\leq 2\exp\left(-\frac{0.1^2\norm{\mb{v}_j}_1^2}{18b^2(0.3/\mu_{ij})^2}\right),
\ee
%where for the second inequality again we applied the fact that $\tau\sqrt{\phi_i\phi_j} \le 0.3/\mu_{ij}$. 
We now must show that the the probability on the right-hand side of \refs{eq:vijcond2} is
exponentially small.  First, we observe that
\crcount{probdenom}
\begin{eqnarray*}
\tau^2\phi_i\phi_j&\le& (\tau_u')^2\phi_i\phi_j \\
&\le &
c_{\thetauubd}(\mb{p})^2\phi_i\phi_j/\phi_1^2 \\
&\le &
 c_{\theprobdenom}(\mb{p}),
\end{eqnarray*}
where $c_{\theprobdenom}(\mb{p}) = c_{\thetauubd}(\mb{p})^2\rho_\sigma^2\rho_m\rho_n$ with $c_{\thetauubd}(\mb{p})$ defined in \refs{eq:tauucondef} below. This follows from the fact that, for any $i,j=1,\ldots,k$,
\begin{equation}
\phi_i/\phi_j=(\bar\sigma_i/\bar\sigma_j)\sqrt{m_j/m_i}\sqrt{n_j/n_i}
\le\rho_\sigma\rho_m^{1/2}\rho_n^{1/2}.
\label{eq:phiratio}
\end{equation}

We now provide a lower bound on $\norm{\mb{v}_j}_1$.
We start with the right
singular vector
$\hat{\mb{v}}(\mb{P}_{j}(\tau))$
of
the matrix
$\mb{P}_{j}(\tau)$.  As noted prior to
\refs{eq:Peig}, this singular vector may be written as
$\hat \alpha_j\bar{\mb{v}}_j+\hat \beta_j\mb{e}_{n_j}$ . 
Let $\mb{v}(\mb{P}_{j}(\tau))$ be the rescaling of
$\hat{\mb{v}}(\mb{P}_{j}(\tau))$ with the scale chosen so that
$\mb{v}(\mb{P}_{j}(\tau))= \alpha_j\bar{\mb{v}}_j+\mb{e}_{n_j}$ 
(i.e., $\beta_j=1$).
Then we can obtain the
value of $\alpha_j$ using the second equation obtained from
\refs{eq:Peig} (see also Lemma 4.5 in \cite{Doan10}), and simplifying
by substituting \refs{eq:gfunc} yields
\begin{equation}
  \alpha_j = \frac{\sqrt{n_j}}{h_i(\tau)} \cdot
\frac{\tau\phi_j[\tau\phi_j -(2a_j(\tau)-\delta_{u,i}\delta_{v,i})h_j(\tau)]}
{\delta_{v,j}h_j(\tau) - \tau\phi_j\delta_{u,j}},
\label{eq:alphaj}
\end{equation}
where $h_j(\tau)=1-\tau\phi_j\mu_{jj}$, which lies in $[0.7,1]$ since 
$\tau\le \tau_u$, and $a_j(\tau)$ is defined as in \refs{eq:aitaudef}.
(Note that the scaling $\beta_j=1$ is valid only if the denominator of the
above fraction is nonzero, which we shall show next.)
Observe that the square-bracketed quantity in the second numerator is 
nonnegative and at least $\tau\phi_j-2$
since $a_j\le 1$ and $\tau\le\tau_u$.

Using the facts that $\delta_u\le \delta_{u,i}\le 1$ and
$\delta_v\le \delta_{v,i}\le 1$
we conclude from \refs{eq:taulsat} that
$\tau\phi_j \ge 2+2\delta_{u,j}/\delta_{v,j}$
and
$\tau\phi_j \ge 2+2\delta_{v,j}/\delta_{u,j}$
 for all $j=1,\ldots,k$ whenever
$\tau\ge\tau_{\ell}$.  

Now, ignoring the additive term of 2 for a moment, 
this assumption
implies that the second denominator
is negative and no more than $\tau\phi_j\delta_{u,j}$ in absolute value. Thus we have:
$$
\alpha_j \le -\frac{\sqrt{n_j}(\tau\phi_j - 2)}{\delta_{u,j}}.
$$
%The fact that $\tau\phi_j\ge 3+ 2\delta_{u,j}/\delta_{v,j}\ge 3$ also
%implies that $\tau\phi_j/3 \le \tau\phi_j-2$, hence we conclude 
%that
%\begin{equation}
%\alpha_j \le -\tau\phi_j\sqrt{n_j}/(3\delta_{u,j}).
%\label{eq:alphajin1}
%\end{equation}
As noted in the previous paragraph $\tau\phi_j-2\ge 2\delta_{u,j}/\delta_{v,j}$, hence
\begin{equation}
\alpha_j \le -2\sqrt{n_j}/\delta_{v,j}.
\label{eq:alphajin}
\end{equation}

Now we write the 1- and 2-norms of $\mb{v}(P_j)$ in terms of $\alpha_j$ and
the other parameters.  Starting with the 1-norm,
\begin{eqnarray*}
\norm{\mb{v}(P_j)}_1 &= & \norm{\alpha_j\bar{\mb{v}}_j + \mb{e}_{n_j}}_1 \\
& \ge & \norm{\alpha_j\bar{\mb{v}}_j}_1 - n_j \\
&= & |\alpha_j|\sqrt{n_j}\delta_{v,j}-n_j \\
& \ge & |\alpha_j|\sqrt{n_j}\delta_{v,j}/2,
\end{eqnarray*}
where, to obtain the last line, we used the fact
that $|\alpha_j|\sqrt{n_j}\delta_{v,j}/2 \ge n_j$, a consequence of
\refs{eq:alphajin}.
Also, $\norm{\mb{v}(P_j)} \le |\alpha_j| + \sqrt{n_j}$ by the triangle inequality.
%\begin{eqnarray*}
%\norm{\mb{v}(P_j)} &= & \norm{\alpha_j\bar{\mb{v}}_j + \mb{e}_{n_j}} \\
%& = & \sqrt{\alpha_j^2+ 2\alpha_j\delta_{v,j}\sqrt{n_j} + n_j} \\
%&\le & \sqrt{\alpha_j^2 + n_j} \\
%& \le & |\alpha_j| + \sqrt{n_j}
%\end{eqnarray*}
%where the third line was obtained since $\alpha_j<0$, while the fourth follows
%from the standard inequality that the 1-norm dominates the 2-norm.  
Thus,
we conclude that
\begin{eqnarray}
\norm{\hat{\mb{v}}(P_j)}_1 & = &
\frac{\norm{\mb{v}(P_j)}_1}{\norm{\mb{v}(P_j)}}  \nonumber \\
& \ge & \frac{|\alpha_j|\sqrt{n_j}\delta_{v,j}/2}{|\alpha_j| + \sqrt{n_j}} \nonumber\\
& = & \frac{\sqrt{n_j}\delta_{v,j}}{2(1 + \sqrt{n_j}/|\alpha_j|)} \nonumber\\
& \ge & \frac{\sqrt{n_j}\delta_{v,j}}{2(1 + \delta_{v,j}/2)} \nonumber\\
& \ge & \frac{\sqrt{n_j}\delta_{v,j}}{3}  \label{eq:vPjlb} 
\end{eqnarray}

Next, we observe by the triangle inequality that
\begin{eqnarray}
\norm{\mb{v}_j}_1 &\ge & 
\norm{\hat{\mb{v}}(\mb{P}_j)}_1 - 
\norm{\hat{\mb{v}}(\mb{P}_j)-\mb{v}_j}_1 \nonumber \\
& \ge & 
\norm{\hat{\mb{v}}(\mb{P}_j)}_1 - 
\sqrt{n_j}\norm{\hat{\mb{v}}(\mb{P}_j)-\mb{v}_j}. \label{eq:vPjdiff}
\end{eqnarray}
We will use Wedin's theorem on perturbation of singular vectors (see Doan and Vavasis \cite{Doan10} and references therein for details) to
analyze the final norm in the above inequality since $\mb{v}_j$ is the
leading singular vector of $\mb{P}_j(\tau)+\tau\phi_j\mb{Q}_{jj}$ while
$\hat{\mb{v}}(\mb{P}_j)$ is the leading singular vector of $\mb{P}_j(\tau)$.

For Wedin's theorem, we choose $\mb{A}=\mb{P}_j(\tau)$, 
$\mb{T} = \tau\phi_j\mb{Q}_{jj}$, and $\mb{B}= \mb{A}+\mb{T}$.  We have: $\Vert \mb{T} \Vert \le \tau\phi_j (m_jn_j)^{3/8}$. In addition,
\begin{eqnarray*}
\sigma_1(\mb{B}) &\ge & \sigma_1(\mb{A}) - \sigma_1(\mb{T}) \\
& \ge & \sqrt{m_jn_jg_j(\tau;a_j)} - \tau\phi_j(m_jn_j)^{3/8},
\end{eqnarray*}
where the second line is obtained from
\refs{eq:normpieq}.
Finally, using \refs{eq:sigma2Pi},
$$\sigma_2(\mb{A})=
\tau\phi_jh_j(\tau)\left(\frac{m_jn_j(1-\delta_{u,j}^2)(1-\delta_{v,j}^2)}
{g_j(\tau;a_j)}\right)^{1/2}.$$
Therefore,
\begin{eqnarray*}
\sin\theta\left(\mb{v}_j,\hat{\mb{v}}(\mb{P}_j(\tau))\right)
& \leq & 
\frac{\tau\phi_j (m_jn_j)^{3/8}}
{\sqrt{m_jn_jg_j(\tau;a_i)} - \tau\phi_j(m_jn_j)^{3/8}
-\tau\phi_jh_j(\tau)\displaystyle\left(\frac{m_jn_j(1-\delta_{u,j}^2)(1-\delta_{v,j}^2)}
{g_j(\tau;a_j)}\right)^{1/2}} \\
& = &
\frac{ (m_jn_j)^{-1/8}}
{\sqrt{g_j(\tau;a_j)}/(\tau\phi_j) - (m_jn_j)^{-1/8}
-h_j(\tau)\displaystyle\left(\frac{(1-\delta_{u,j}^2)(1-\delta_{v,j}^2)}
{g_j(\tau;a_j)}\right)^{1/2}}.
\end{eqnarray*}
Observe that the numerator tends to zero like $(m_jn_j)^{-1/8}$ while the denominator
does not depend on $m_jn_j$ (except for a vanishing term).  Furthermore, the
denominator is positive; this follows from the fact that the
first term in the denominator is at least 0.5 by \refs{eq:sqrtglblinear}
whereas the last term is at most $1/\sqrt{5}$ again by \refs{eq:sqrtglblinear}
and the fact that $\phi_i\tau\ge 5$ thanks to \refs{eq:tauldef1}.

This shows that
\begin{equation}
\label{eq:vPpert}
\norm{\mb{v}_j-\hat{\mb{v}}(\mb{P}_j(\tau))}_2\leq O\left((m_jn_j)^{-\frac{1}{8}}\right).
\end{equation}

Combining \refs{eq:vPjlb}, \refs{eq:vPjdiff}
and \refs{eq:vPpert}, we can then pick a constant less than $1/3$, say $0.3$, 
and claim that
\be
\label{eq:vlb}
\norm{\mb{v}_j}_1\geq 0.3\delta_{v,j}\sqrt{n_j}\geq 0.3\delta_{v}\sqrt{n_j},
\ee
as long as $m_j$, $n_j$ are large.
Combining this bound with \refs{eq:vijcond2}, we can claim that the probability
\refs{eq:rij2} 
 is exponential small:
\be
\label{eq:vijcond2f}
\mathbb{P}\left(\tau\card{\frac{\mb{R}_{ij}(s,:)\mb{v}_j}{\norm{\mb{v}_j}_1}}>\frac{1}{3}\right)\leq 
2\exp\left(-\frac{0.1^20.3^2\delta_{v}^2n_j}{18b^2c_{\theprobdenom}(\mb{p})}\right).
\ee

Similarly, the first probability \refs{eq:rij1}
can also be proved to be exponentially small using the analogous 
lower bound of $\norm{\mb{u}_i}_1$:
\be
\label{eq:ulb}
\norm{\mb{u}_i}_1\geq 0.3\delta_{u}\sqrt{m_i}.
\ee
The bound for the first probability can therefore be written as follows:
\be
\label{eq:vijcond1f}
\mathbb{P}\left(\tau\card{\frac{\mb{u}_i^T\mb{R}_{ij}(:,t)}{\norm{\mb{u}_i}_1}}>\frac{1}{3}\right)
\leq 
2\exp\left(-\frac{0.1^20.3^2\delta_{u}m_i}{18b^2c_{\theprobdenom}(\mb{p})}\right).
\ee
For the third probability \refs{eq:rij3}, we again use the fact that $\tau\sqrt{\phi_i\phi_j}\leq 0.3/\mu_{ij}$ since $\tau\le\tau_u$:
$$
\mathbb{P}\left(\tau\card{\frac{\mb{u}_i^T\mb{R}_{ij}\mb{v}_j}{\norm{\mb{u}_i}_1
\norm{\mb{v}_j}_1}}>\frac{1}{3}\right)
\leq
\mathbb{P}\left(\tau\card{\frac{\mb{u}_i^T(\mb{R}_{ij}-\mu_{ij}\sqrt{\phi_i\phi_j}
\mb{e}_{m_i}\mb{e}_{n_j}^T)\mb{v}_j}{\norm{\mb{u}_i}_1\norm{\mb{v}_j}_1}}>
\frac{1}{30}\right),
$$
where $\mb{u}_i^T(\mb{R}_{ij}/\sqrt{\phi_i\phi_j}-
\mu_{ij}\mb{e}_{m_i}\mb{e}_{n_j}^T)\mb{v}_j$
is a $b$-subgaussian random variable since
$\norm{\mb{u}_i}_2=\norm{\mb{v}_j}_2=1$. We again can bound this
probability using the lower bounds of $\norm{\mb{u}_i}_1$ and
$\norm{\mb{v}_j}_1$ as follows: 
\be
\label{eq:vijcond3f}
\mathbb{P}\left(\tau\card{\frac{\mb{u}_i^T\mb{R}_{ij}\mb{v}_j}{\norm{\mb{u}_i}_1\norm{\mb{v}_j}_1}}>\frac{1}{3}\right)\leq 
2\exp\left(-\frac{0.1^20.3^4\delta_{u}^2\delta_{v}^2m_in_j}{18b^2c_{\theprobdenom}(\mb{p})}\right).
\ee
Combining \refs{eq:vijcond2f}, \refs{eq:vijcond1f}, \refs{eq:vijcond3f}, we obtain the following tail bound:
\begin{eqnarray}
\mathbb{P}\left(\norm{\mb{V}_{ij}}_{\infty}>1\right)&\leq& 
2\exp\left(-\frac{0.1^20.3^2\delta_{v}^2n_j}{18b^2c_{\theprobdenom}(\mb{p})}\right)  \nonumber\\
& & \,\mbox{} + 
2\exp\left(-\frac{0.1^20.3^2\delta_{u}^2m_j}{18b^2c_{\theprobdenom}(\mb{p})}\right) \nonumber\\
& & \,\mbox{} + 
2\exp\left(-\frac{0.1^20.3^4\delta_{u}^2\delta_{v}^2m_in_j}{18b^2c_{\theprobdenom}(\mb{p})}\right).\label{eq:vijcondf}
\end{eqnarray} 
\end{subsubsection}
\begin{subsubsection}{Upper bound on $\norm{{\mb{W}}_{ij}}$}
The second constraint for this type of block is $\displaystyle\norm{\mb{W}_{ij}}\leq\frac{1}{k+1}\min_{i=1,\ldots,k}\sigma_i$. Using the fact that $\mb{Q}_{ij}=\mb{R}_{ij}/\sqrt{\phi_i\phi_j}-\mu_{ij}\mb{e}_{m_i}\mb{e}_{n_j}^T$, we have:
$$
\mb{W}_{ij}=\displaystyle\tau\theta\sqrt{\phi_i\phi_j}\left(\mb{Q}_{ij}-\frac{\mb{e}_{m_i}\mb{u}_i^T\mb{Q}_{ij}}{\norm{\mb{u}_i}_1} - \frac{\mb{Q}_{ij}\mb{v}_j\mb{e}_{n_j}^T}{\norm{\mb{v}_j}_1}+\frac{\mb{u}_i^T\mb{Q}_{ij}\mb{v}_j}{\norm{\mb{u}_i}_1\norm{\mb{v}_j}_1}\mb{e}_{m_i}\mb{e}_{n_j}^T\right).
$$

We will establish that 
$\displaystyle\norm{\mb{W}_{ij}}\leq\frac{1}{k+1}\min_{i=1,\ldots,k}\sigma_i$ 
by showing that 
\begin{eqnarray}
\tau\theta\sqrt{\phi_i\phi_j}\norm{\frac{\mb{Q}_{ij}}{2}-\frac{\mb{e}_{m_i}\mb{u}_i^T\mb{Q}_{ij}}{\norm{\mb{u}_i}_1}}&\le& \frac{1}{3(k+1)}\min_{i=1,\ldots,k}\sigma_i,
\label{eq:qij1}  \\
\tau\theta\sqrt{\phi_i\phi_j}\norm{\frac{\mb{Q}_{ij}}{2}-\frac{\mb{Q}_{ij}\mb{v}_j\mb{e}_{n_j}^T}{\norm{\mb{v}_j}_1}}&\le&\frac{1}{3(k+1)}\min_{i=1,\ldots,k}\sigma_i, 
\label{eq:qij2}  \\
\tau\theta\sqrt{\phi_i\phi_jm_in_j}\card{\frac{\mb{u}_i^T\mb{Q}_{ij}\mb{v}_j}{\norm{\mb{u}_i}_1\norm{\mb{v}_j}_1}}&\le&\frac{1}{3(k+1)}\min_{i=1,\ldots,k}\sigma_i.
\label{eq:qij3}
\end{eqnarray}
Given that $m_i,n_i\rightarrow\infty$ for all $i=1,\ldots,k$, we make the
following assumption:
\crcount{cdelta}
\begin{equation}
(m_jn_i)^{-1/8}\le
\frac{.47}{3(k+1)(\bar{\rho}_m\bar{\rho}_n)^{1/8}\rho_{\sigma}c_{\thecdelta}(\mb{p})},
\label{eq:sigmas2}
\end{equation}
for all $i,j=1,\ldots,k$,
where
we introduce
\begin{equation}
c_{\thecdelta}(\mb{p}) = 
\max\left\{\frac{1}{2}+\frac{1}{0.3\delta_u},
\frac{1}{2}+\frac{1}{0.3\delta_v},
\frac{1}{0.3^2\delta_u\delta_v}\right\}.
\label{eq:cthecdelta}
\end{equation}
Now, inequality \refs{eq:qij1} is derived as follows:
\begin{eqnarray*}
\norm{\frac{\mb{Q}_{ij}}{2}-\frac{\mb{e}_{m_i}\mb{u}_i^T\mb{Q}_{ij}}{\norm{\mb{u}_i}_1}}
& \le &
\norm{\mb{Q}_{ij}}\cdot\norm{\frac{\mb{I}}{2}- 
\frac{\mb{e}_{m_i}\mb{u}_i^T}{\norm{\mb{u}_i}_1}} \\
& \le & 
(m_in_j)^{3/8}\cdot\left(\frac{1}{2} +\frac{\sqrt{m_i}}{\norm{\mb{u}_i}_1}\right) \\
& \le &
(m_in_j)^{3/8}\cdot\left(\frac{1}{2} +\frac{1}{0.3\delta_{u,i}}\right).
\end{eqnarray*}
The first line uses submultiplicativity of the 2-norm
since we have:
$$
\mb{Q}_{ij}/2-\mb{e}_{m_i}\mb{u}_i^T\mb{Q}_{ij}/\norm{\mb{u}_i}_1
=(\mb{I}/2- 
\mb{e}_{m_i}\mb{u}_i^T/\norm{\mb{u}_i}_1)\mb{Q}_{ij}.
$$
The second uses the triangle
inequality, and the third uses \refs{eq:ulb}.
Multiply by the scalar $\tau\theta\sqrt{\phi_i\phi_j}$ and let $l=1,\ldots,k$
be arbitrary:
\begin{eqnarray*}
\tau\theta\sqrt{\phi_i\phi_j}\norm{\frac{\mb{Q}_{ij}}{2}-\frac{\mb{e}_{m_i}\mb{u}_i^T\mb{Q}_{ij}}{\norm{\mb{u}_i}_1}}
& \le & 
\tau\theta\sqrt{\phi_i\phi_j}
(m_in_j)^{3/8}\cdot\left(\frac{1}{2} +\frac{1}{0.3\delta_{u,i}}\right) \\
& = &
\tau\theta
\sqrt{\bar\sigma_i\bar\sigma_j}
\frac{m_i^{1/8}n_j^{1/8}}{m_j^{1/4}n_i^{1/4}}
\cdot\left(\frac{1}{2} +\frac{1}{0.3\delta_{u,i}}\right) \\
& \le &
\frac{.47\tau\theta}{3(k+1)}\bar\sigma_l  \\
&\le &
\frac{\sigma_l}{3(k+1)}.
\end{eqnarray*}
The third line follows from \refs{eq:sigmas2} and
the last from \refs{eq:sigmaibarsigmailb}.
Inequality 
\refs{eq:qij2} is established using the same argument.
Finally, \refs{eq:qij3} is established by a similar
argument starting from the inequality
$|\mb{u}_i^T\mb{Q}_{ij}\mb{v}_j|\le\norm{\mb{Q}_{ij}}\leq (m_in_j)^{3/8}$.
\end{subsubsection}
\end{subsection}
\begin{subsection}{Analysis for block $(k+1,j)$, $j=1,\ldots,k$}
We now consider the $(k+1,j)$ block. Similar to the above approach, we will construct the following matrix $\mb{V}_{k+1,j}$:
$$
\mb{V}_{k+1,j}=\tau\frac{\bar{\mb{R}}_{k+1,j}\mb{v}_j\mb{e}_{n_j}^T}{\norm{\mb{v}_j}_1}.
$$
\begin{subsubsection}{Upper bound on $\norm{{\mb{V}}_{k+1,j}}_\infty$}
The condition $\norm{\mb{V}_{k+1,j}}_{\infty}\leq 1$ can be dealt with using the same approach as before. We have:
$$
V_{k+1,j}(s,t)=\tau\frac{\bar{\mb{R}}_{k+1,j}(s,:)\mb{v}_j}{{\norm{\mb{v}_j}_1}}.
$$
Since $\tau\le\tau_u$, we can conclude from \refs{eq:tauudef4} that
$\tau\le 0.9/(\mu_{ij}\sqrt{\phi_i\phi_j})$ for all $i=k+1,\ldots,k_0$. Thus, we have
$$
\mathbb{P}\left(\tau\card{\frac{\bar{\mb{R}}_{k+1,j}(s,:)\mb{v}_j}{\norm{\mb{v}_j}_1}}>1\right)
\leq
\mathbb{P}\left(\tau\card{\frac{(\bar{\mb{R}}_{k+1,j}(s,:)-\mu_{i(s),j}\sqrt{\phi_{i(s)}\phi_j}\mb{e}_{n_j}^T)\mb{v}_j}{\norm{\mb{v}_j}_1}}>0.1\right),
$$
where $i(s)$ is the corresponding original block (row) index for the 
$s$th row of $\bar{\mb{R}}_{k+1,j}$. Since $\norm{\mb{v}_j}=1$, 
$(\bar{\mb{R}}_{k+1,j}(s,:)/\sqrt{\phi_{i(s)}\phi_j}-\mu_{i(s),j}\mb{e}_{n_j}^T)\mb{v}_j$ 
is a 
$b$-subgaussian random variable. Thus, by \req{eq:tsubgau}, we have:
$$
\mathbb{P}\left(\tau\card{\frac{(\bar{\mb{R}}_{k+1,j}(s,:)-\mu_{i(s),j}\sqrt{\phi_{i(s)}\phi_j}\mb{e}_{n_j}^T)\mb{v}_j}{\norm{\mb{v}_j}_1}}>0.1\right)\leq  
2\exp\left(-\frac{0.1^2\norm{\mb{v}_j}_1^2}{2b^2\tau^2\phi_{i(s)}\phi_j}\right).
$$
To show this is exponentially small, we first analyze the denominator.  
We start by noting that
\crcount{probdenomtwo}
\begin{eqnarray*}
\tau^2\phi_{i(s)}\phi_j & \le & (\tau_u')^2\phi_{i(s)}\phi_j \\
& \le & c_{\thetauubd}(\mb{p})^2\phi_{i(s)}\phi_j/\phi_1^{2} \\
&\le & c_{\thetauubd}(\mb{p})^2c_{\thenoisephi}\phi_j/\phi_1 \\
& \le & c_{\thetauubd}(\mb{p})^2c_{\thenoisephi}\rho_{\sigma}(\rho_m\rho_n)^{1/2} \\
& \equiv & c_{\theprobdenomtwo}(\mb{p}).
\end{eqnarray*}
The second line was obtained from
\refs{eq:tauucondef} and the third from \refs{eq:noisephi}, and the last
line introduces another constant.
Combining with \refs{eq:vlb} for the numerator, we obtain the following tail bound:
\be
\label{eq:vk1jcondf}
\mathbb{P}\left(\norm{\mb{V}_{k+1,j}}_{\infty}>1\right)\leq 2\exp\left(-\frac{0.1^2
0.3^2\delta_{v}^2n_j}{2b^2c_{\theprobdenomtwo}(\mb{p})}\right).
\ee
\end{subsubsection}
\begin{subsubsection}{Upper bound on $\norm{{\mb{W}}_{k+1,j}}$}
Now consider $\mb{W}_{k+1,j}$.  It is clear that $\mb{W}_{k+1,j}\mb{v}_j=\mb{0}$. In addition, we have:
\begin{equation}
\mb{W}_{k+1,j}=\tau\theta\mb{\Phi}
\bar{\mb{Q}}_{k+1,j}\left(\mb{I} 
-\frac{\mb{v}_j\mb{e}_{n_j}^T}{\norm{\mb{v}_j}_1}\right),
\label{eq:wk+1form}
\end{equation}
where $\bar{\mb{Q}}_{k+1,j}\in\R^{\bar{m}_{k+1}\times n_j}$ is a $b$-subgaussian matrix 
that is a concatenation of $\mb{Q}_{lj}$, $l=k+1,\ldots,k_0$ and
$$\mb{\Phi} = 
\left(
\begin{array}{ccc}
\sqrt{\phi_{k+1}\phi_j}\mb{I}_{m_{k+1}} & &  \\
& \ddots & \\
& & \sqrt{\phi_{k_0}\phi_j}\mb{I}_{m_{k_0}}
\end{array}
\right).$$
By the same argument as before,
\begin{equation}
\left\Vert \mb{I} 
-\frac{\mb{v}_j\mb{e}_{n_j}^T}{\norm{\mb{v}_j}_1}\right\Vert \le 
1 + \frac{1}{0.3\delta_{v,j}}\leq c_{\thecdelta}(\mb{p})+1/2.
\label{eq:Iveineq}
\end{equation}
where $c_{\thecdelta}(\mb{p})$ was defined by \refs{eq:cthecdelta}.
Also,
\begin{equation}
\Vert \mb{\Phi}\Vert = \sqrt{\bar\phi_{k+1}\phi_j},
\label{eq:normPhi}
\end{equation}
where $\displaystyle\bar\phi_{k+1}=\max_{i=k+1,\ldots,k_0}\phi_i$.
Now suppose 
\crcount{qkponej}
\begin{equation}
\Vert \bar{\mb{Q}}_{k+1,j}\Vert \le c_{\theqkponej}(\mb{p})(m_jn_j)^{1/2}/\sqrt{c_{\thenoisephi}},
\label{eq:barqbd}
\end{equation}
where 
\begin{equation}
c_{\theqkponej}(\mb{p}) = \frac{0.47}{\rho_\sigma(c_{\thecdelta}(\mb{p})+1/2)(k+1)}
\label{eq:qkponej}
\end{equation}
and $c_{\thenoisephi}$ was defined in \refs{eq:noisephi}.
(Below we will argue that \refs{eq:barqbd} happens with high probability.) 

Using the hypothesis \refs{eq:barqbd},
\begin{eqnarray*}
\norm{\mb{W}_{k+1,j}} & \le &
\tau\theta \norm{\mb{\Phi}}\cdot \norm{\mb{I} 
-\frac{\mb{v}_j\mb{e}_{n_j}^T}{\norm{\mb{v}_j}_1}}\cdot
\norm{\bar{\mb{Q}}_{k+1,j}} \\
& \le &
\tau\theta\sqrt{\bar\phi_{k+1}\phi_j}(c_{\thecdelta}(\mb{p})+1/2)
c_{\theqkponej}(\mb{p})(m_jn_j)^{1/2}/\sqrt{c_{\thenoisephi}} \\
&\le &
\tau\theta\phi_j(c_{\thecdelta}(\mb{p})+1/2)c_{\theqkponej}(\mb{p})(m_jn_j)^{1/2} \\
&= &
\tau\theta\bar{\sigma_j}\frac{0.47}{\rho_\sigma(k+1)} \\
&\le &
\frac{0.47\tau\theta}{k+1}\bar{\sigma}_k \\
& \le &
\frac{1}{k+1}\cdot\min_{i=1,\ldots,k}\sigma_k.
\end{eqnarray*}
The first line follows from
\refs{eq:wk+1form}, the second from \refs{eq:normPhi},
\refs{eq:Iveineq}, and \refs{eq:barqbd}.
The third and fifth follow from \refs{eq:noisephi} and \refs{eq:rhosig_def}
respectively, and the last from \refs{eq:sigmaibarsigmailb}.

Now we show that the hypothesis
\refs{eq:barqbd} holds with high probability using
Lemma~\ref{lem:2norm}.  As mentioned above, $\bar{m}_{k+1}$ denotes the number of
rows of $\bar{\mb{Q}}_{k+1,j}$, i.e., $m_{k+1}+\cdots+m_{k_0}$.
\begin{eqnarray*}
\mathbb{P}\left(\norm{\bar{\mb{Q}}_{k+1,j}}>
c_{\theqkponej}(\mb{p})(m_jn_j)^{1/2}/\sqrt{c_{\thenoisephi}}
\right)
&\leq&
\exp\left(-\frac{\displaystyle 8c_{\theqkponej}(\mb{p})^2}{81b^2c_{\thenoisephi}}m_jn_j+(\log 7)(\bar{m}_{k+1} + n_j)\right) \\
&=&\exp\left(-\frac{\displaystyle 4c_{\theqkponej}(\mb{p})^2}{81b^2c_{\thenoisephi}}m_jn_j+(\log 7)\bar{m}_{k+1}\right)\\
&&\,\mbox{}\cdot\exp\left(-\frac{\displaystyle 4c_{\theqkponej}(\mb{p})^2c_{\thenoisephi}}{81b^2}m_jn_j+(\log 7)n_j\right).
\end{eqnarray*}
The second exponent  in the second line tends to $-\infty$ linearly with $m_j$; the
first exponent also tends to $-\infty$ linearly provided that 
\begin{equation}
\frac{\bar{m}_{k+1}}{m_jn_j}\le K <
\frac{4c_{\theqkponej}(\mb{p})^2}{81b^2c_{\thenoisephi}(\log 7)},
\label{eq:mkmjnj}
\end{equation}
where $K$ is some constant (independent of $m_i,n_i$ for any $i$),
which holds under the assumption \refs{eq:mnp1bd}.
\end{subsubsection}
The analysis of $(i,k+1)$ block is similar for $i=1,\ldots,k$.
\end{subsection}
\begin{subsection}{Analysis for block $(k+1,k+1)$}
\begin{subsubsection}{Upper bound on $\norm{{\mb{V}}_{k+1,k+1}}_{\infty}$}
For the last block $(k+1,k+1)$, we will simply construct
$\mb{V}_{k+1,k+1}\in\mathbb{R}^{\bar{m}_{k+1}\times\bar{n}_{k+1}}$
from $(k_0-k)^2$ sub-blocks
$\mb{V}_{st}^{(k+1)}\in\mathbb{R}^{m_s\times n_t}$,
$$
\mb{V}_{st}^{(k+1)}=\tau\mu_{st}\sqrt{\phi_s\phi_t}\mb{e}_{m_s}\mb{e}_{n_t}^T, 
\quad s,t=k+1,\ldots,k_0.
$$
Since $\tau\le \tau_u$, by \refs{eq:tauudef5}
we have: $\tau\leq 0.9/(\mu_{st}\sqrt{\phi_s\phi_t})$ for all $s,t=k+1,\ldots,k_0$. 
Thus, $\norm{\mb{V}_{k+1,k+1}}_{\infty}\leq 1$.
\end{subsubsection}
\begin{subsubsection}{Upper bound on $\norm{{\mb{W}}_{k+1,k+1}}$}
We have, $\mb{W}_{k+1,k+1}$ is composed of
blocks: $\mb{W}_{s,t}^{(k+1)}=\tau\theta\left(\sqrt{\phi_s\phi_t}\bar{\mb{Q}}_{st}+\bar{\mb{B}}_{s,t}\right)$, where $\bar{\mb{Q}}_{s,t}\in\mathbb{R}^{m_s\times n_t}$.  We will write 
the sum as:
$$\mb{W}_{k+1,k+1}=\tau\theta(\mb{\Phi}_2\bar{\mb{Q}}_{k+1,k+1}\mb{\Phi_3}+\bar{\mb{B}}_{k+1,k+1})$$
where $\bar{\mb{Q}}_{k+1,k+1}$ contains entries chosen from a 
$b$-subgaussian distribution, and
$$
\mb{\Phi}_{2} = \left(
\begin{array}{ccc}
\sqrt{\phi_{k+1}}\mb{I}_{m_{k+1}} & & \\
& \ddots & \\
& & 
\sqrt{\phi_{k_0}}\mb{I}_{m_{k_0}} 
\end{array}
\right),
$$
and
$$
\mb{\Phi}_{3} = \left(
\begin{array}{ccc}
\sqrt{\phi_{k+1}}\mb{I}_{n_{k+1}} & & \\
& \ddots & \\
& & 
\sqrt{\phi_{k_0}}\mb{I}_{n_{k_0}} 
\end{array}
\right).
$$

We have: $\displaystyle\norm{\bar{\mb{B}}_{k+1,k+1}}=\max_{l=k+1,\ldots,k_0}\bar{\sigma}_l=\bar{\sigma}_{k+1}$ and $\norm{\mb{\Phi}_2} = \norm{\mb{\Phi}_3}=(\bar\phi_{k+1})^{1/2}$, where
$\bar\phi_{k+1}$ was defined as in \refs{eq:normPhi}.
Thus
\begin{equation}
\norm{\mb{W}_{k+1,k+1}}\leq\tau\theta\bar\phi_{k+1}\norm{\bar{\mb{Q}}_{k+1,k+1}}+\tau\theta\bar{\sigma}_{k+1}.
\label{eq:wk+1k+1}
\end{equation}

Applying the assumption \refs{eq:sigmanoisebd} to the second term of
\refs{eq:wk+1k+1}, we have:
\begin{eqnarray*}
\tau\theta\bar{\sigma}_{k+1} & \le & 
\frac{0.47\tau\theta\bar\sigma_k}{2(k+1)} \\
&\le & \frac{1}{2(k+1)}\cdot\min_{i=1,\ldots,k}\sigma_i.
\end{eqnarray*}
The second line follows from \refs{eq:sigmaibarsigmailb}.

Now turning to the first term, let us suppose that
\begin{equation}
\norm{\bar{\mb{Q}}_{k+1,k+1}} \le \frac{0.23\sqrt{m_ln_l}}{(k+1)c_{\thenoisephi}},
\label{eq:qk+1k+1}
\end{equation}
where $c_{\thenoisephi}$ is from \refs{eq:noisephi} and $l$ is
the index of the min of $\sigma_1,\ldots,\sigma_k$.
(Below we will argue that this holds with probability exponentially close
to 1.) Then
\begin{eqnarray*}
\tau\theta\bar\phi_{k+1}\norm{\bar{\mb{Q}}_{k+1,k+1}}
& \le & \frac{0.23\tau\theta\phi_{k+1}\sqrt{m_ln_l}}{(k+1)c_{\thenoisephi}} \\
& \le & \frac{0.23\tau\theta\phi_l\sqrt{m_ln_l}}{k+1} \\
& = & \frac{0.46\tau\theta\bar{\sigma}_l}{2(k+1)} \\
& \le & \frac{1}{2(k+1)}\cdot\min_{i=1,\ldots,k}\sigma_i.
\end{eqnarray*}
The second line  uses \refs{eq:qk+1k+1} and
the last line uses \refs{eq:sigmaibarsigmailb} and the choice of $l$.
Thus, we have analyzed both of the terms of \refs{eq:wk+1k+1} and established $\displaystyle\norm{\mb{W}_{k+1,k+1}}\le \frac{1}{k+1}\min_{i=1,\ldots,k}\sigma_i$ as required.

We now analyze the probability that \refs{eq:qk+1k+1} fails.
According to Lemma~\ref{lem:2norm}
$$
\mathbb{P}\left(\norm{\bar{\mb{Q}}_{k+1,k+1}}>
\frac{0.23\sqrt{m_ln_l}}{(k+1)c_{\thenoisephi}}\right)
\leq
\exp\left(-\frac{8\cdot 0.23^2m_ln_l}{81b^2(k+1)c_{\thenoisephi}}+(\log 7)(\bar{m}_{k+1} + \bar{n}_{k+1})\right). 
$$
This quantity tends to zero exponentially fast as long as
\begin{equation}
\frac{\bar{m}_{k+1}+\bar{n}_{k+1}}{\min_{i=1,\ldots,k}m_in_i}\le 
K < \frac{8\cdot 0.23^2}{81b^2(k+1)c_{\thenoisephi}(\log 7)},
\label{eq:mk1nk1}
\end{equation}
where $K$ is some constant (independent of the matrix size),
which holds under the assumption \refs{eq:mnp1bd}.
\end{subsubsection}
\end{subsection}
\begin{subsection}{Definitions of the scalars}
The definitions of the scalars appearing in the theorem and the proof
can now be provided based on the inequalities developed during the proof.

We start by defining $\tau_\ell$ as follows:
\begin{equation}
\tau_\ell = c_{\thetauellbd}(\mb{p})\phi_1^{-1},
\label{eq:tauldef} 
\end{equation}
where $c_{\thetauellbd}(\mb{p})$ is defined as
\begin{eqnarray}
c_{\thetauellbd}(\mb{p}) & = & \rho_\sigma\sqrt{\rho_m\rho_n}\max\left\{6(k+1)\sqrt{\rho_m\rho_n},\frac{4}{\xi_u},\frac{4}{\xi_v},2 + \frac{2}{\delta_u}, 2 + \frac{2}{\delta_v}\right\}. \label{eq:tauellbd2}
\end{eqnarray} 
Applying inequality \refs{eq:phiratio}, the following inequalities that have already been used in the preceding analysis indeed hold:
\begin{eqnarray}
\tau_\ell & \ge & 6(k+1)\sqrt{\rho_m\rho_n}\max_{i=1,\ldots,k}\phi_i^{-1}, \label{eq:tauldef1}\\
\tau_\ell& \ge & \max\left\{\frac{4}{\xi_u},\frac{4}{\xi_u}\right\}\max_{i=1,\ldots,k}\phi_i^{-1}, \label{eq:tauldef4} \\
\tau_\ell& \ge & \left(2+\max\left\{\frac{2}{\delta_u},\frac{2}{\delta_v}\right\}\right)\max_{i=1,\ldots,k}\phi_i^{-1}. \label{eq:taulsat} 
\end{eqnarray}
The constant $c_{\thethetarange}(\mb{p})$ is then defined as follows:
\begin{equation}
c_{\thethetarange}(\mb{p}) = \frac{1}{2}\left(1.2^4\left(c_{\thetauellbd}(\mb{p})\right)^2\left(\frac{k\rho_m\rho_n}{k+\rho_m\rho_n-1}\right)+1\right)^{-1/2}.
\label{eq:cthetarangedef}
\end{equation}

Next, we define
\begin{equation}
\tau_u = c_{\thetauubd}(\mb{p})\phi_1^{-1},
\label{eq:tauudef}
\end{equation}
where
\begin{equation}
c_{\thetauubd}(\mb{p}) = \frac{4}{3}\left(\frac{1+(k-1)\sqrt{\rho_m\rho_n}}{1+(k-1)\rho_\sigma^{-2}}+\sqrt{\left(\frac{1+(k-1)\sqrt{\rho_m\rho_n}}{1+(k-1)\rho_\sigma^{-2}}\right)^2+\frac{(1+(k-1)\rho_m\rho_n)\left(c_{\thethetarange}(\mb{p})^{-2}-1\right)}{1+(k-1)\rho_\sigma^{-2}}}\right).
\label{eq:tauucondef}
\end{equation}
Note that $\displaystyle c_{\thetauellbd}(\mb{p})=\frac{25}{36}\cdot\sqrt{\frac{k+\rho_m\rho_n-1}{k\rho_m\rho_n}}\cdot\sqrt{c_{\thethetarange}(\mb{p})^{-2}/4-1}$ from \refs{eq:cthetarangedef}, which implies $c_{\thetauellbd}(\mb{p})<c_{\thetauubd}(\mb{p})$ or $\tau_\ell<\tau_u$.

We now define
\begin{eqnarray}
c_{\themubd}(\mb{p},c_{\thenoisephi}) & = & 
\frac{1}{c_{\thetauubd}(\mb{p})}\cdot\min\left\{\frac{0.3}{\rho_{\sigma}\sqrt{\rho_m\rho_n}},\frac{0.9}{\left(c_{\thenoisephi}\rho_{\sigma}\sqrt{\rho_m\rho_n}\right)^{1/2}},\frac{0.9}{\sqrt{c_{\thenoisephi}}}\right\}.
\label{eq:c1def}
\end{eqnarray}
Clearly, since $c_{\thetauubd}(\mb{p})\geq (4/3)\left(1+\left(c_{\thethetarange}(\mb{p})\right)^{-1}\right)\geq 4$ and $\rho_{\sigma},\rho_m,\rho_n\geq 1$, we have: $c_{\themubd}(\mb{p},c_{\thenoisephi})\leq 0.075<0.08$.

Now, using \refs{eq:phiratio} and the upper bound $\phi_i/\phi_j\le c_{\thenoisephi}$ for all
$i=1,\ldots,k+1$ and all $j=1,\ldots,k$ (a restatement of
\refs{eq:noisephi}), the following inequalities indeed hold:
\begin{eqnarray}
\tau_u &\le& \min_{i=1,\ldots,k}\frac{0.3}{\mu_{ii}\phi_i}, \label{eq:tauudef1} \\
\tau_u &\le& \min_{i,j=1,\ldots,k}\frac{0.3}{\mu_{ij}\sqrt{\phi_i\phi_j}}, \label{eq:tauudef3} \\
\tau_u &\le& \min_{i=k+1,\ldots,k_0; j=1,\ldots,k}\frac{0.9}{\max\left\{\mu_{ij},\mu_{ji}\right\}\sqrt{\phi_i\phi_j}}, \label{eq:tauudef4} \\
\tau_u &\le& \min_{i,j=k+1,\ldots,k_0}\frac{0.9}{\mu_{ij}\sqrt{\phi_i\phi_j}}.
\label{eq:tauudef5} 
\end{eqnarray}

The last scalar to define is
$c_{\theminiNoise}(\mb{p},c_{\thenoisephi},b)$. We define it as follows:
\begin{equation}
c_{\theminiNoise}(\mb{p},c_{\thenoisephi},b) =
\min\left(
\frac{4c_{\theqkponej}(\mb{p})^2}{81b^2c_{\thenoisephi}(\log 7)},
\frac{8\cdot 0.23^2}{81b^2(k+1)c_{\thenoisephi}(\log 7)}\right),
\label{eq:c0def}
\end{equation}
where $c_{\theqkponej}(\mb{p})$ was defined by
\refs{eq:qkponej}.
\end{subsection}
\section{Numerical Examples}
\label{sec:ex}
\subsection{Biclique example}
We consider a simple example that involves a bipartite graph $G=(U,V,E)$ with two non-overlapping bicliques given by $U_1\times V_1$ and $U_2\times V_2$, where $U_1\cap U_2=\emptyset$ and $V_1\cap V_2=\emptyset$. The remaining edges in $E$ are inserted at random with probability $p$. The $U$-to-$V$ adjacency matrix can be written in the form $\mb{A}=\mb{B}+\mb{R}$, where $\mb{B}$ is a block diagonal matrix with $k_0=3$ diagonal blocks, the last of which is a block of all zeros while the other two of which are blocks of all ones. If $U_1\cup U_2=U$ and $V_1\cup V_2=V$, we can consider $\mb{B}$ with just $k_0=2$ diagonal blocks. We also assume that $\card{U_1}=\card{U_2}=1/2\card{U}=m/2$ and $\card{V_1}=\card{V_2}=1/2\card{V}=n/2$. We would like to find these $k=2$ planted bicliques within the graph $G$ under the presence of random noise simultaneously. 

\begin{figure}[htp]
\centering
\includegraphics[width=0.7\textwidth]{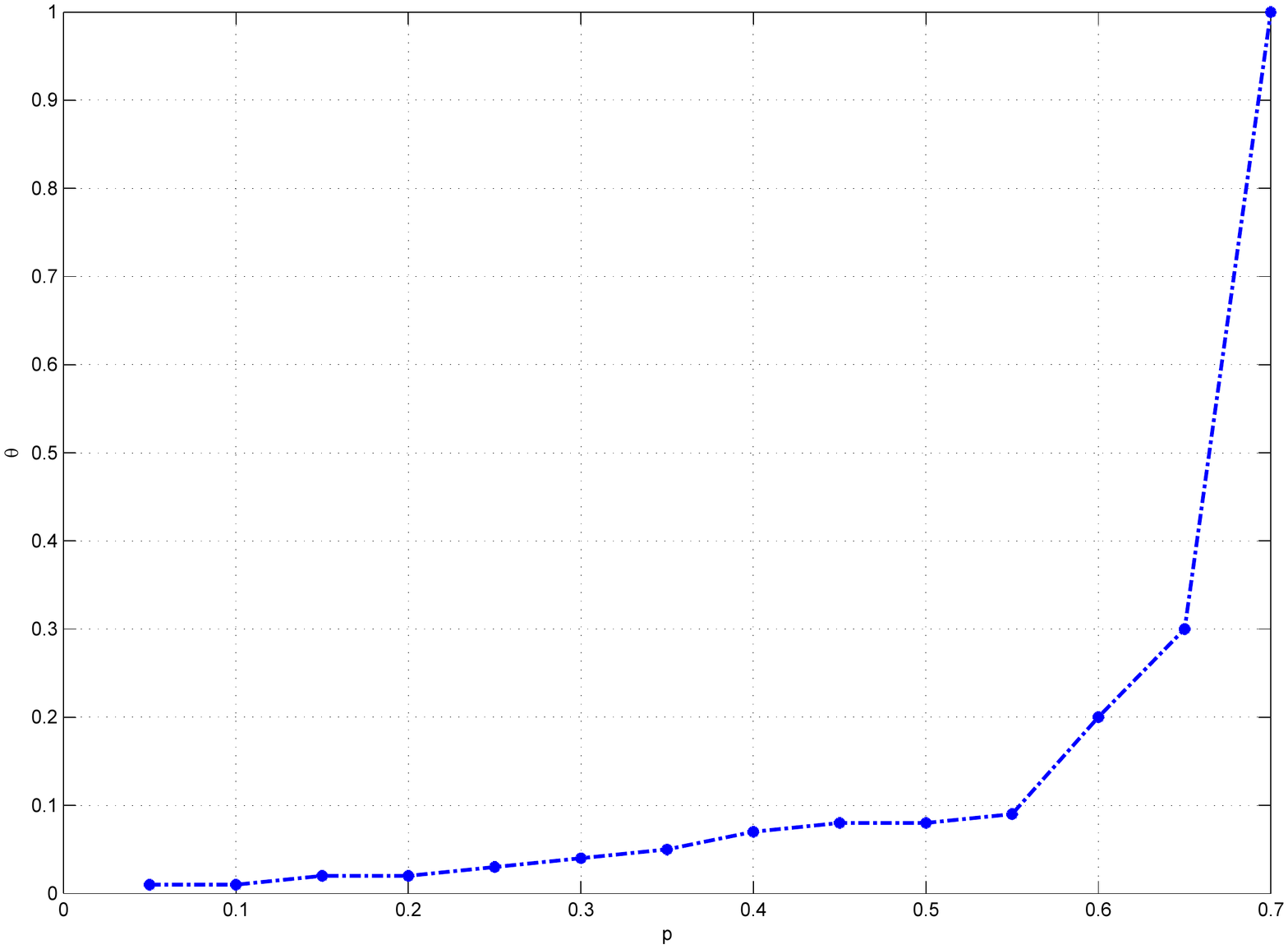}
\caption{Minimum values of $\theta$ to recover two planted bicliques for different values of $p$}\label{fig:mintheta}
\end{figure}

\begin{figure}[htp]
\centering
\includegraphics[width=0.7\textwidth]{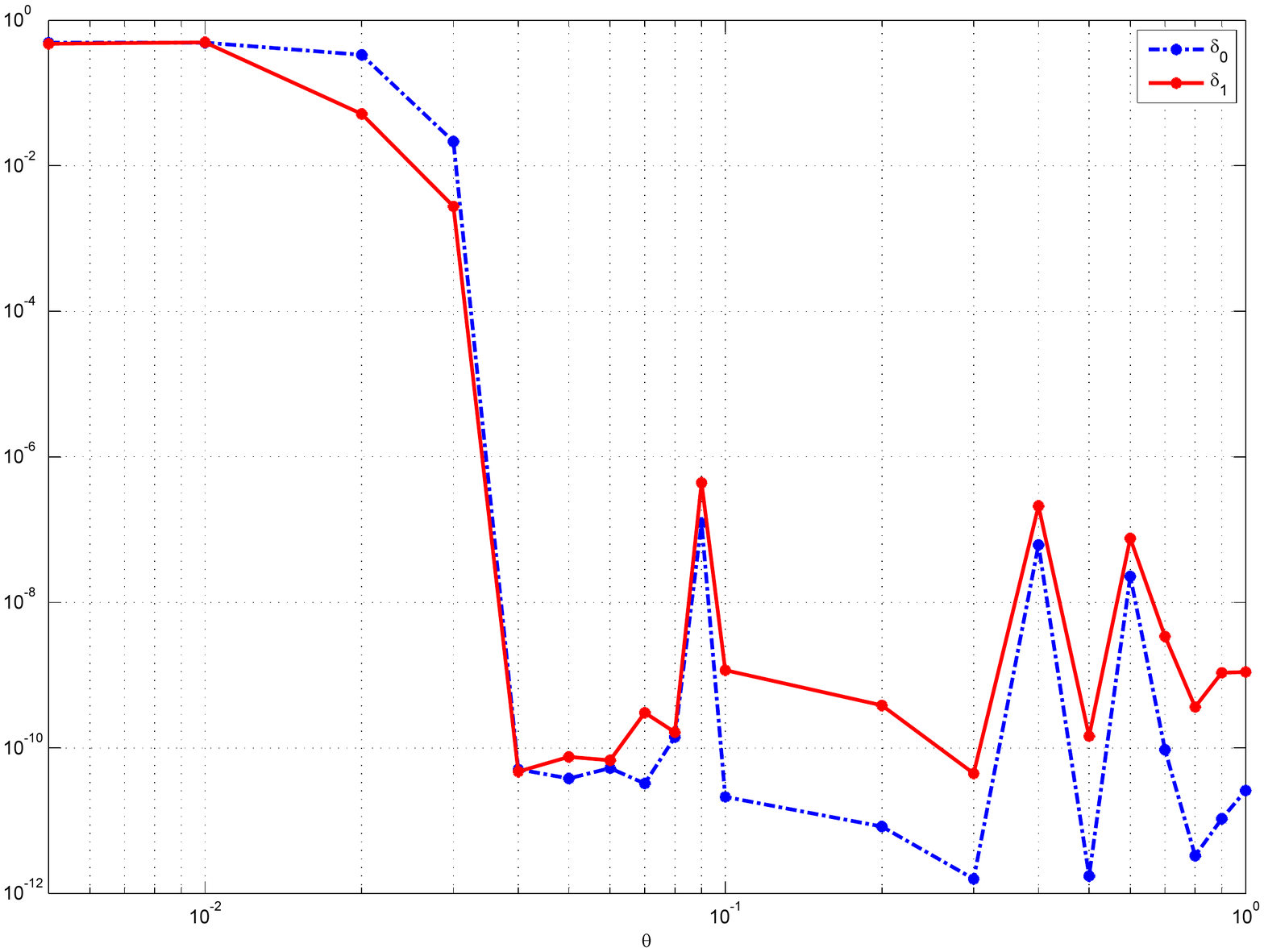}
\caption{Maximum differences between entries in diagonal blocks and off-diagonal blocks for $p=0.30$}\label{fig:differences}
\end{figure}

For this example, $\bar{\mb{u}}_i=\mb{e}_{m_i}/\sqrt{m_i}$ and $\bar{\mb{v}}_i=\mb{e}_{n_i}/\sqrt{n_i}$ for $i=1,2$. In addition, $\bar\sigma_i=\sqrt{m_in_i}$, $i=1,2$, which means $\phi_1=\phi_2=1$. We can then choose $\rho_u=\rho_v=1$, $\xi_u=\xi_v=1$, $\pi_u=\pi_v=1$, $\rho_m=\rho_n=1$, and $\rho_\sigma=1$. Under the random setting described above, $\mu_{ij}=p$ for all $i\neq j=1,2$. Given that $k=k_0$, we can set $c_0=0$ and there is no need to consider the conditions related to noise blocks. With $\bar{\mb{u}}_i=\mb{e}_{m_i}/\sqrt{m_i}$ and $\bar{\mb{v}}_i=\mb{e}_{n_i}/\sqrt{n_i}$ for $i=1,2$, the analysis is simpler and we only need $c_{\thetauellbd}(\mb{p})=2$ since \refs{eq:tauldef1} and \refs{eq:tauldef4} are not needed while \refs{eq:taulsat} can be relaxed to $\displaystyle\tau_\ell\ge 2\max_{i=1,\ldots,k}\phi_i^{-1}$. The constant $c_{\thethetarange}(\mb{p})$ has a better approximation:
$$
c_{\thethetarange}(\mb{p})=\frac{1}{2}\left(\frac{36}{25}c_{\thetauellbd}(\mb{p})-1\right)^{-1}=\frac{25}{94}\approx 0.266.
$$
We then can compute $c_{\themubd}(\mb{p},c_{\thenoisephi})$ as follows:
$$
c_{\themubd}(\mb{p},c_{\thenoisephi})=0.3/c_{\thetauubd}(\mb{p})=(0.9/4)\left(1+\left(c_{\thethetarange}(\mb{p})\right)^{-1}\right)^{-1}\approx 0.047,
$$ 
which means with $p\leq 0.047$, we are able to recover two planted cliques using the proposed convex formulation in \refs{eq:nprob} with $0.376\cdot(mn)^{-1/2}\leq\theta\leq 0.752\cdot(mn)^{-1/2}$ with high probability. The results are quite restricted given the way how we construct the dual solutions solely based on matrices of all ones. Having said that, these conditions are theoretical sufficient conditions. Practically, the convex formulation \refs{eq:nprob} with a wider range of $\theta$ can recover planted bicliques under the presence of more random noise, i.e., higher probability $p$. The numerical computation is performed with CVX \cite{cvx} for the biclique example discussed here with $m=n=50$. We test the problem with $10$ values of $p$ ranging from $0.05$ to $0.95$. For each value of $p$, we construct a random matrix $\mb{A}$ and solve \refs{eq:nprob} with $20$ different values of $\theta$ ranging from $0.005$ to $1.0$. The solution $\bX$ is scaled so that the maximum value of its entries is $1$. We compare $\bX$ and $\mb{B}$ by taking the maximum differences between their entries in diagonal blocks of $\mb{B}$, $\delta_1$, and that of off-diagonal blocks, $\delta_0$. For this example, we are not able to recover two planted bicliques, i.e., the block diagonal structure of the matrix $\mb{B}$, with $\theta=0.005$ for any $p$ given large values for $\delta_0$ and $\delta_1$. It is due to the fact that for smaller values of $\theta$, the objective of achieving better rank-$2$ approximation is more prominent than the objective of achieving the sparse structure. {In addition, we cannot recover the two bicliques for $p\geq 0.75$}. Figure \ref{fig:mintheta} shows the minimum values $\theta_{\min}(p)$ of $\theta$ with which \refs{eq:nprob} can be used to recover the planted bicliques when there is a significant reduction in the values of $\delta_0$ and $\delta_1$. The graph indicates that we need larger $\theta$ for the settings with more random noise. Figure \ref{fig:differences} plots these differences (in log scale) for $p=0.30$ and we can see that $\delta_0$ and $\delta_1$ change from $10^{-2}$ to $10^{-10}$ between $\theta=0.03$ and $\theta=0.04$. When the planted bicliques can be recovered, all of these values are in the order of $10^{-6}$ or less, which indicates the recovery ability of our proposed formulation for this example under the presence of noise. Note that for this special example of binary data, the range of the values of $\theta$ with which two planted bicliques can be recovered is usually large enough to cover the whole remaining interval $[\theta_{\min}(p),1]$ considered in this experiment.
%Figure \ref{fig:zerobound} and \ref{fig:onebound} plot the actual values of these differences, $\delta_0$ and $\delta_1$, for different values of $\theta$ and $p$ when the planted bicliques can be recovered. We can see that all of these values are in the order of $10^{-6}$ or less, which indicates the recovery ability of our proposed formulation for this example under the presence of noise.
 
\begin{figure}[htp]
\centering
\includegraphics[width=0.7\textwidth]{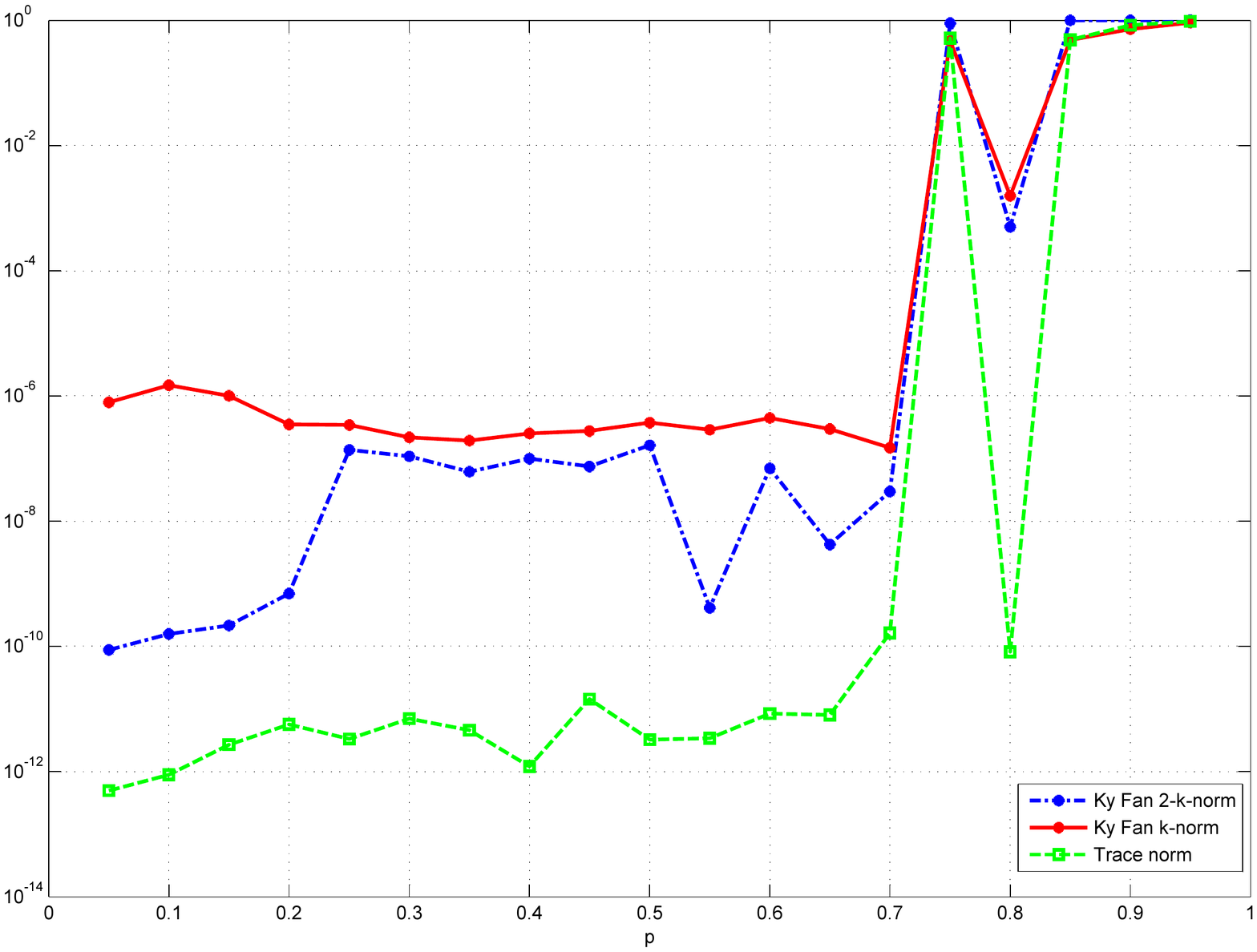}
\caption{Maximum differences between entries for three models}\label{fig:differentnorms}
\end{figure}

\begin{comment}
\begin{figure}[htp]
\centering
\includegraphics[width=0.7\textwidth]{zerobound}
\caption{Maximum difference between entries in off-diagonal blocks for different values of $p$ and $\theta$}\label{fig:zerobound}
\end{figure} 

\begin{figure}[htp]
\centering
\includegraphics[width=0.7\textwidth]{onebound}
\caption{Maximum difference between entries in diagonal blocks for different values of $p$ and $\theta$}\label{fig:onebound}
\end{figure}
\end{comment}
{Under the setting of this experiment, two blocks have the same size, i.e., $m_1n_1=m_2n_2=mn/4$, which means $\bar{\sigma}_1=\bar{\sigma}_2$. As mentioned previously, if we replace the Ky Fan $2$-$k$-norm in \refs{eq:nprob} by the Ky Fan $k$-norm, it is likely that we can still retrieve the information of singular vectors, which is enough for this experiment. We now run the Ky Fan $k$-norm formulation with different levels of noise by varying $p$ from $0.05$ to $0.95$. Similarly, we also test the trace norm formulation proposed by Ames \cite{Ames13} under the Bernoulli model with $\alpha=1$ and $\beta=p$ given this is a biclique instance. Figure \ref{fig:differentnorms} show the plots of $\max\{\delta_0,\delta_1\}$ obtained from the three different models. It shows that all of three models can handle noisy instances with $p\leq 0.7$ with the trace norm model achieving the best result in terms of accuracy. It is due to the fact that if the trace norm model is successful, it returns the (unique) exact solution. In the next examples, we will demonstrate that if singular values are needed as parts of the recovery result, both Ky Fan $k$-norm and the trace norm model are not able to deliver.
}

{\subsection{Examples with synthetic gene expression data}
In this section, we apply our formulation for synthetic gene expression data sets studied in Preli\'c et al. \cite{prelic2006systematic}. Under this setting, biclusters are \emph{transcription modules}, which are defined by a set of genes ${\cal G}_i$ and a set of experimental conditions ${\cal C}_i$. Preli\'c et al. \cite{prelic2006systematic} provide two types of biclusters, constant clusters with binary gene expression matrices, which are similar to data inputs in the bicliqe problem, and additive clusters with integer gene expression matrices. We will focus on additive clusters in this section. Following Preli\'c et al. \cite{prelic2006systematic}, we will examine the effects of noise with $k=10$ non-overlapping transcription modules, each of which consists of $10$ genes and $5$ experimental conditions. The resulting gene expression matrices $\mb{E}$ are $100\times 50$ matrices with element values range from $0$ to $100$. Within the implanted biclusters, the values are at least $50$ while the background values, i.e., outside the biclusters, are less than $50$. {Furthermore, average gene expression values are different from one implanted bicluster to another and within each bicluster, the values are also different from one another}. We add random normal noise, $r_{ij}\sim N(0,(50\sigma)^2)$, where $\sigma$ is the noise level, $0\leq\sigma\leq 0.1$, to the gene expression values while maintaining their non-negativity, i.e., $e_{ij}\leftarrow \max\{e_{ij}+r_{ij},0\}$. More details of how to construct these gene expression matrices can be found in Preli\'c et al. \cite{prelic2006systematic}. 
 
In order to compare different biclustering methods, Preli\'c et al. \cite{prelic2006systematic} defined a match score of two biclusters ${\cal B}=({\cal G}_i,{\cal C}_i)_{i=1,\ldots,k}$ and ${\cal B}'=({\cal G}'_i,{\cal C}'_i)_{i=1,\ldots,k}$ as
\be
\label{eq:score}
S^*_G({\cal B},{\cal B}')=\frac{1}{k}\sum_{i=1}^k\max_{j=1,\ldots,k}\frac{\card{{\cal G}_i\cap{\cal G}'_j}}{\card{{\cal G}_i\cup{\cal G}'_j}}.
\ee
Clearly, $S^*_G({\cal B},{\cal B}')\in [0,1]$ and $S^*_G({\cal B},{\cal B}')=1$ if $\cal B$ and ${\cal B}'$ are the same. The match score is not symmetric and given the implanted bicluster ${\cal B}^*$, each biclustering method with the resulting bicluster ${\cal B}$ is measured by two measures, the \emph{average bicluster relevance}, $S^*_G({\cal B},{\cal B}^*)$, and the \emph{average module recovery}, $S^*_G({\cal B}^*,{\cal B})$. According to Preli\'c et al. \cite{prelic2006systematic}, we can also define a similar match score $S^*_C$ for experimental conditions. Having said that, to be consistent with the comparative study discussed in Preli\'c et al. \cite{prelic2006systematic}, we will focus only on $S^*_G$ match scores. In addition, for these gene expression applications, we also believe that it is of greater importance to correctly determine the clustering of the genes rather than of the experimental conditions. Now, for each noise level between $0$ and $0.1$, we will generate $10$ noisy gene expression matrices and as in Preli\'c et al. \cite{prelic2006systematic}, the two performance measures will be averaged over these $10$ instances. Similar to the biclique example, we solve \refs{eq:nprob} with $20$ different values of $\theta$ ranging from $0.005$ to $1.0$. For each run, the resulting matrix is scaled to best approximate the (noisy) input matrix, i.e., to minimize $\norm{\alpha\mb{X}^*-\mb{E}}$, and element values are rounded down to zeros according to an appropriate threshold. The threshold is determined when there is a significant ratio (usually in the order of $10^3$) between two consecutive sorted element values of the resulting matrix. The final computational issue is how to select the appropriate value for the parameter $\theta$. Theoretically, there is a range of $\theta$ in which the recovery holds. For example, when all data blocks are square matrices of size $n$, $\theta$ is required to be in the order of $O(1/(n\sqrt{k}))$. Having said that, it is difficult to find correct constants in practice. For this particular example, we follow the heuristic used in Doan et al. \cite{doan2013proximal}, which finds the balance between the magnitude of the resulting matrix measured by the norm of its $k$-approximation and the approximation averaging effect measured by the norm of the residual. Figure \ref{fig:lcurve} shows the plot of these two measures for our first run without noise ($\sigma=0$) and an appropriate value of $\theta$ can be selected from the distinct middle range. We pick $\theta=0.07$, which is in the middle of that range. Sorted element values of the resulting matrix is plot in Figure \ref{fig:values} and we can see a significant transition (with a ratio of more than $10^4$) between large and small values. The threshold for zero rounding can be set to be $5\times 10^{-4}$ in this case knowing that all larger element values are larger than $5$. 

\begin{figure}[htp]
\centering
\includegraphics[width=0.7\textwidth]{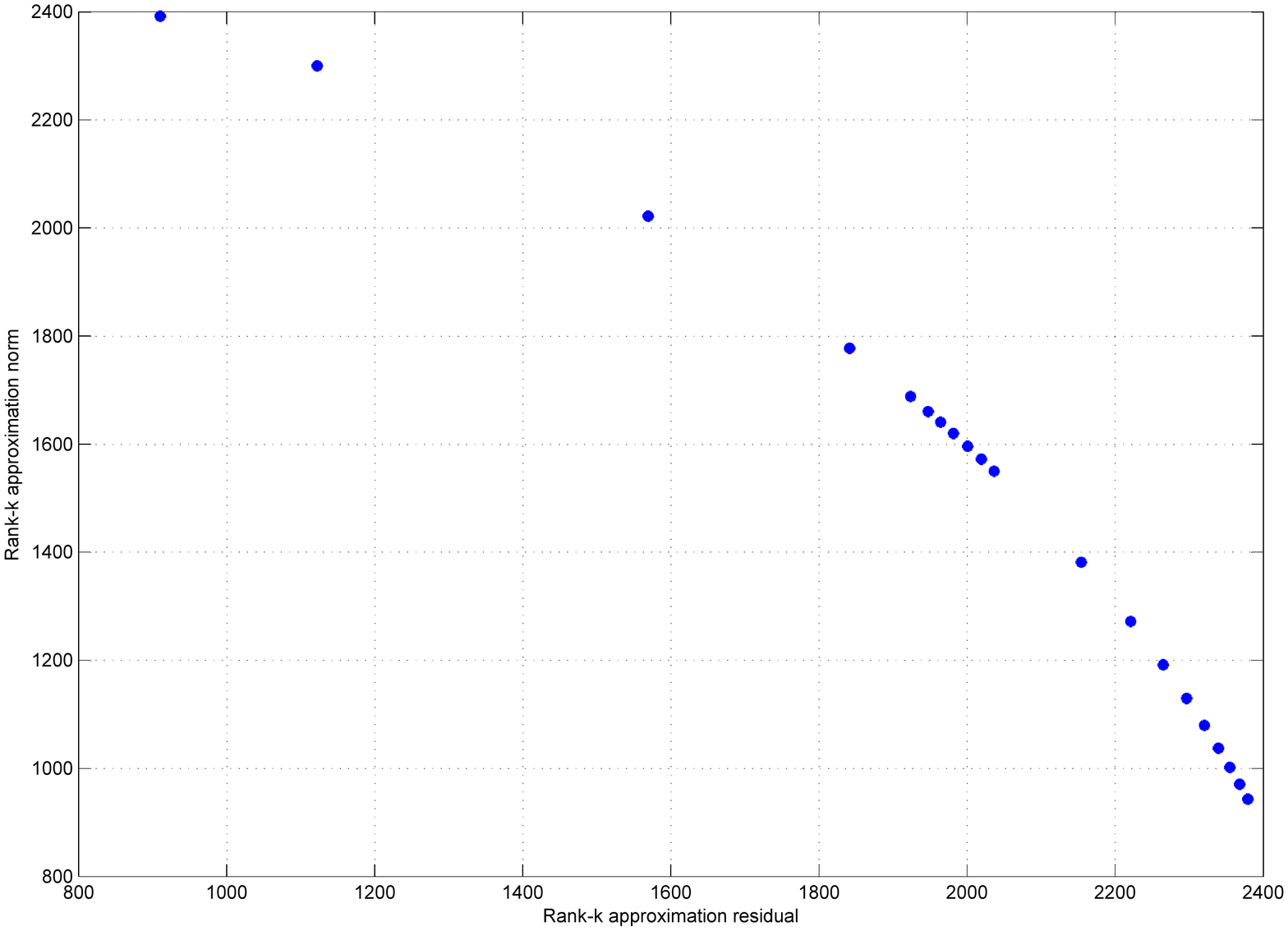}
\caption{Approximation averaging effect vs. magnitude of resulting blocks}\label{fig:lcurve}
\end{figure} 

\begin{figure}[htp]
\centering
\includegraphics[width=0.7\textwidth]{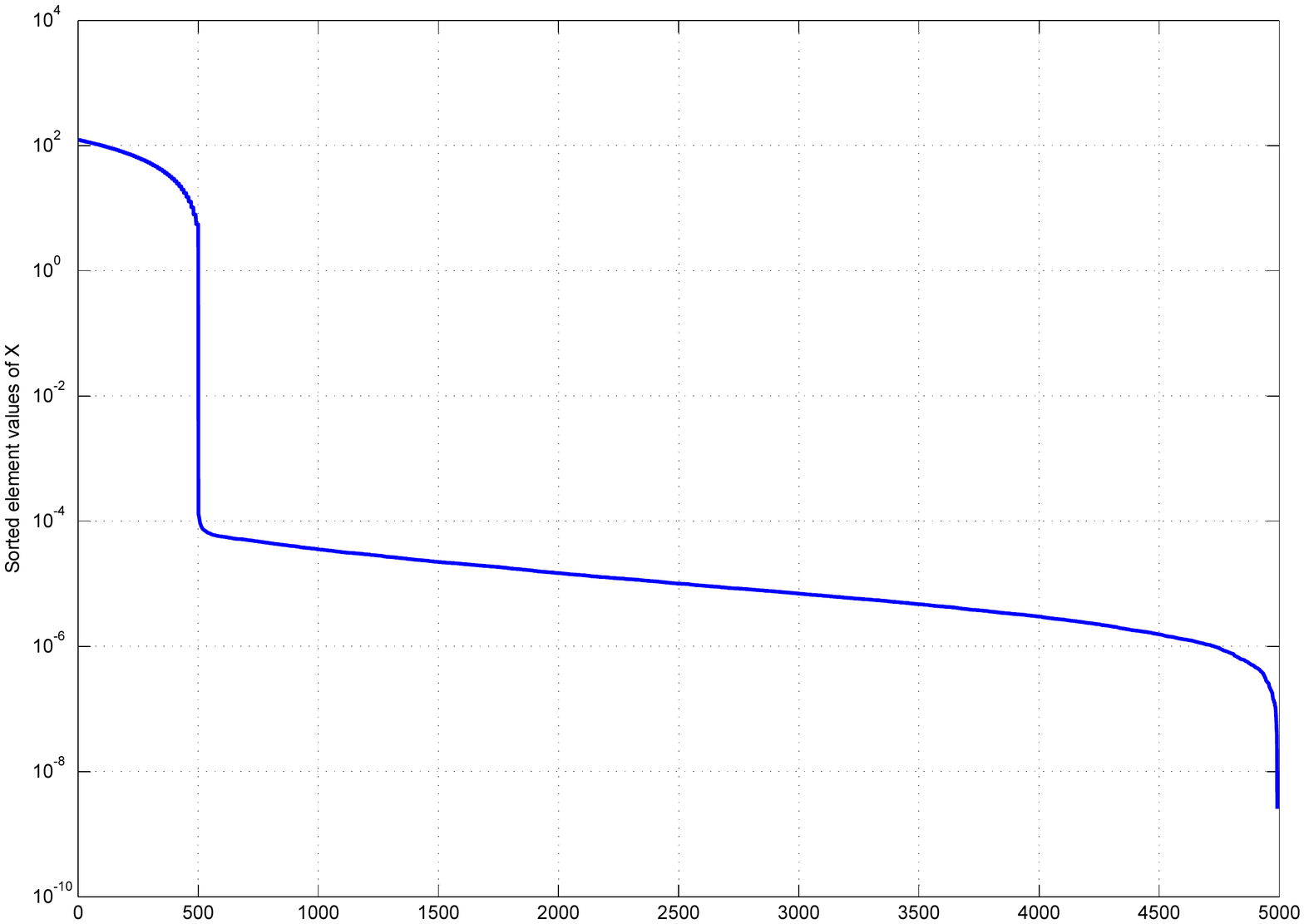}
\caption{Distinction between large and small element values of the resulting matrix}\label{fig:values}
\end{figure} 

The recovered transcription modules are displayed in Figure \ref{fig:recovery} alongside the display of the original gene expression data. It clearly shows that all $10$ transcription module are recovered exactly, which means both performance measure, average bicluster relevance and average module recovery, achieve the maximum value of $1$. {In addition, differences in gene expression levels between different implanted biclusters are present in the recovered transcription modules as in the original gene expression data. We also try to run the Ky Fan $k$-norm formulation and the trace norm model proposed by Ames \cite{Ames13} for the original gene expression data. Figure \ref{fig:newrecovery} shows the recovered transcription modules from the two models. Even though the recovered modules are correct, there is no significant difference in gene expression levels from one implanted bicluster to another as in the original gene expression data in the results of these two models. Furthermore, the trace norm model, which is developed for biclique problems, provides a single gene expression level within each implanted bicluster and this level is the same for all implanted biclusters. It shows that these two models cannot recover the information of singular values as expected.
}
\begin{figure}[htp]
\centering
\includegraphics[width=0.8\textwidth]{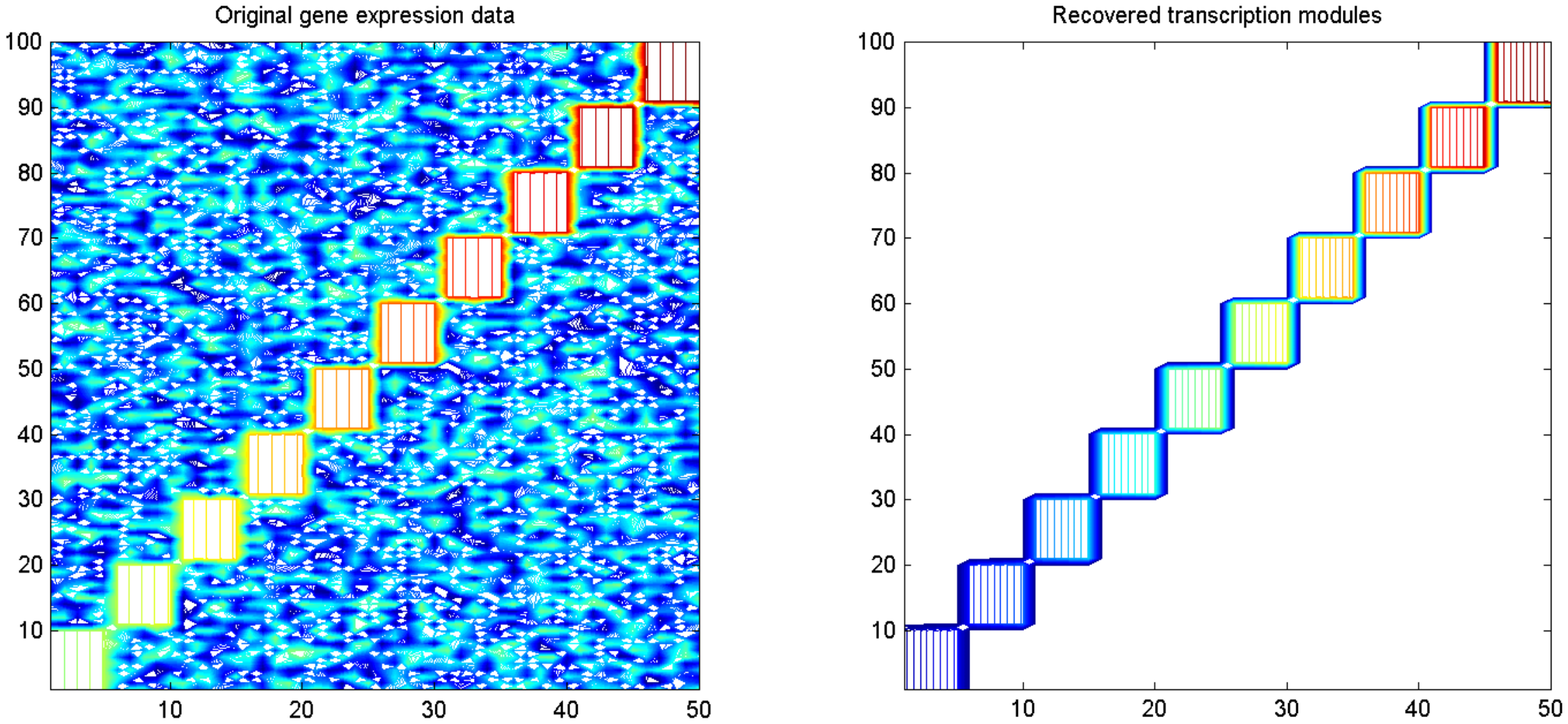}
\caption{Original gene expression data vs. recovered transcription module}\label{fig:recovery}
\end{figure} 

\begin{figure}[htp]
\centering
\includegraphics[width=0.8\textwidth]{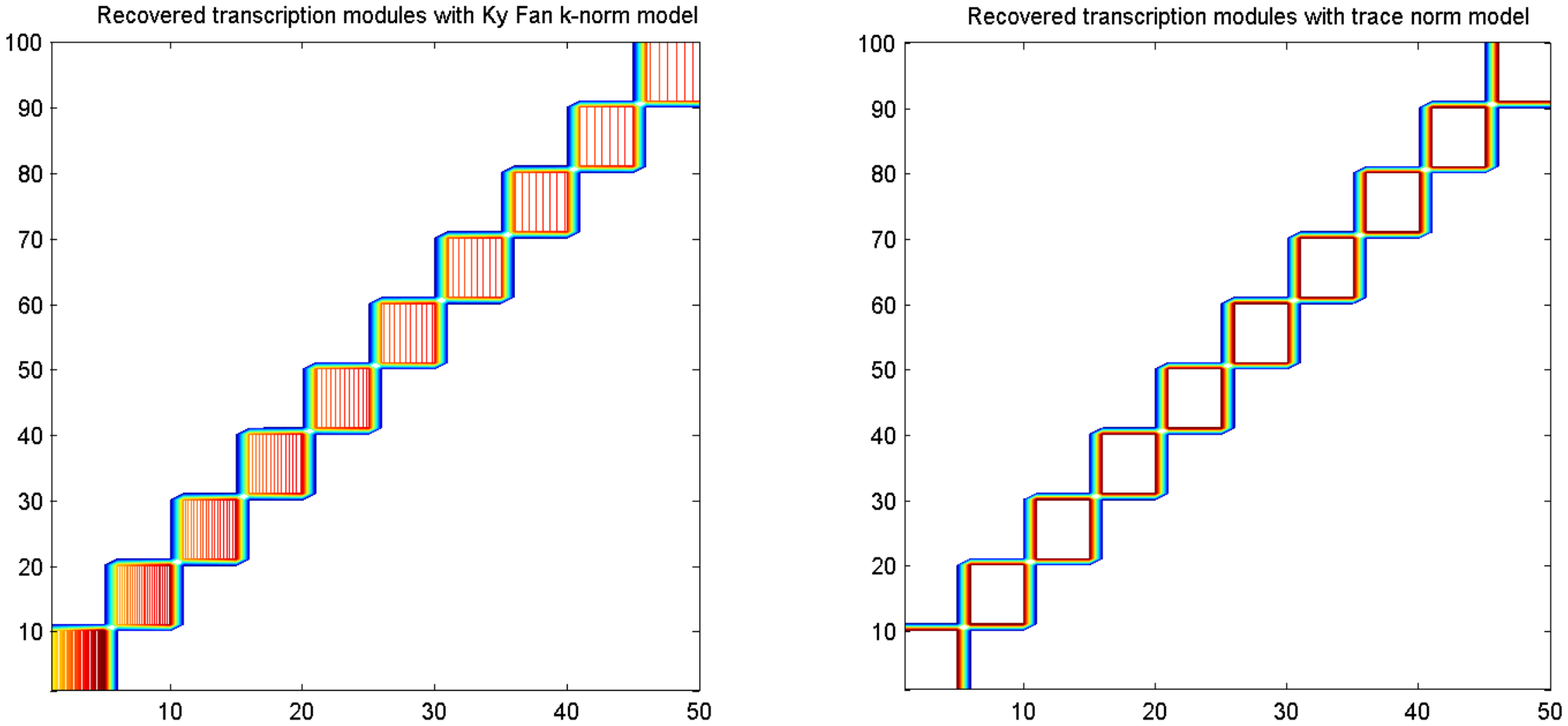}
\caption{Recovered transcription modules from two different models}\label{fig:newrecovery}
\end{figure}

The effect of noise is captured in Figure \ref{fig:matchscores}. Both measures, average bicluster relevance and average module recovery, are the same in these instances and they are very close to 1 with the minimum value is larger than 0.99. As compared to results reported in Preli\'c et al. \cite[Figs. 3(a),3(b)]{prelic2006systematic}, for this particular numerical example, our proposed method is comparable to (if not better) the best algorithms such as BiMax, ISA, and Samba. 

\begin{figure}[htp]
\centering
\includegraphics[width=0.7\textwidth]{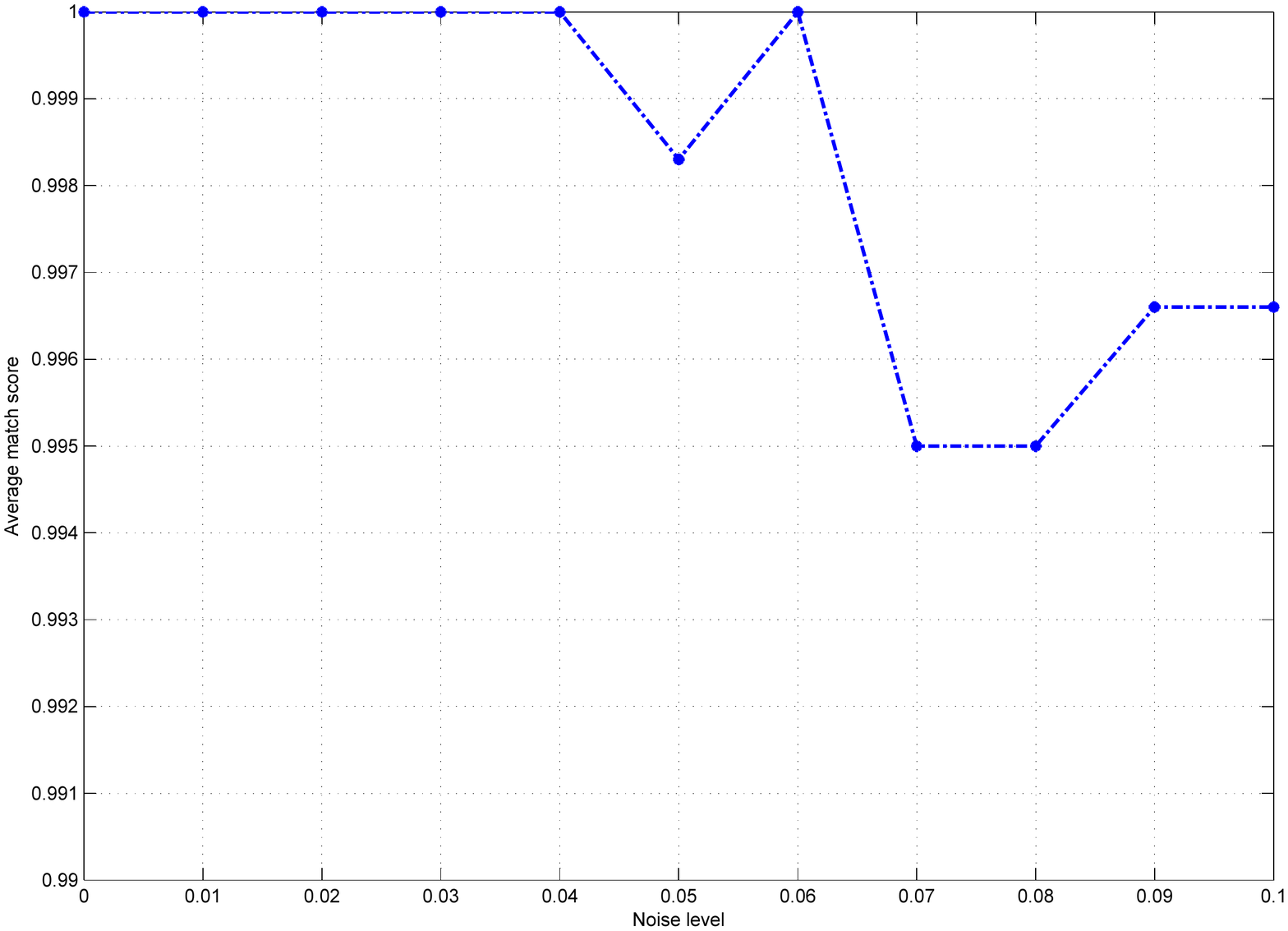}
\caption{Match scores with different noise levels}\label{fig:matchscores}
\end{figure}

When the noise level goes higher, not all of $10$ modules can be recovered given the fact that the noisy background data can be misunderstood for actual expression data. Figure \ref{fig:noisedata} shows an example of noisy gene expression data at the noise level of $\sigma=0.3$. We run the proposed formulation with $k=10$ and recover $6$ largest modules, which are not all perfect. We solve the problem again with $k=6$ instead and achieve much better results. The results are shown in Figure \ref{fig:diffk}.

\begin{figure}[htp]
\centering
\includegraphics[width=0.8\textwidth]{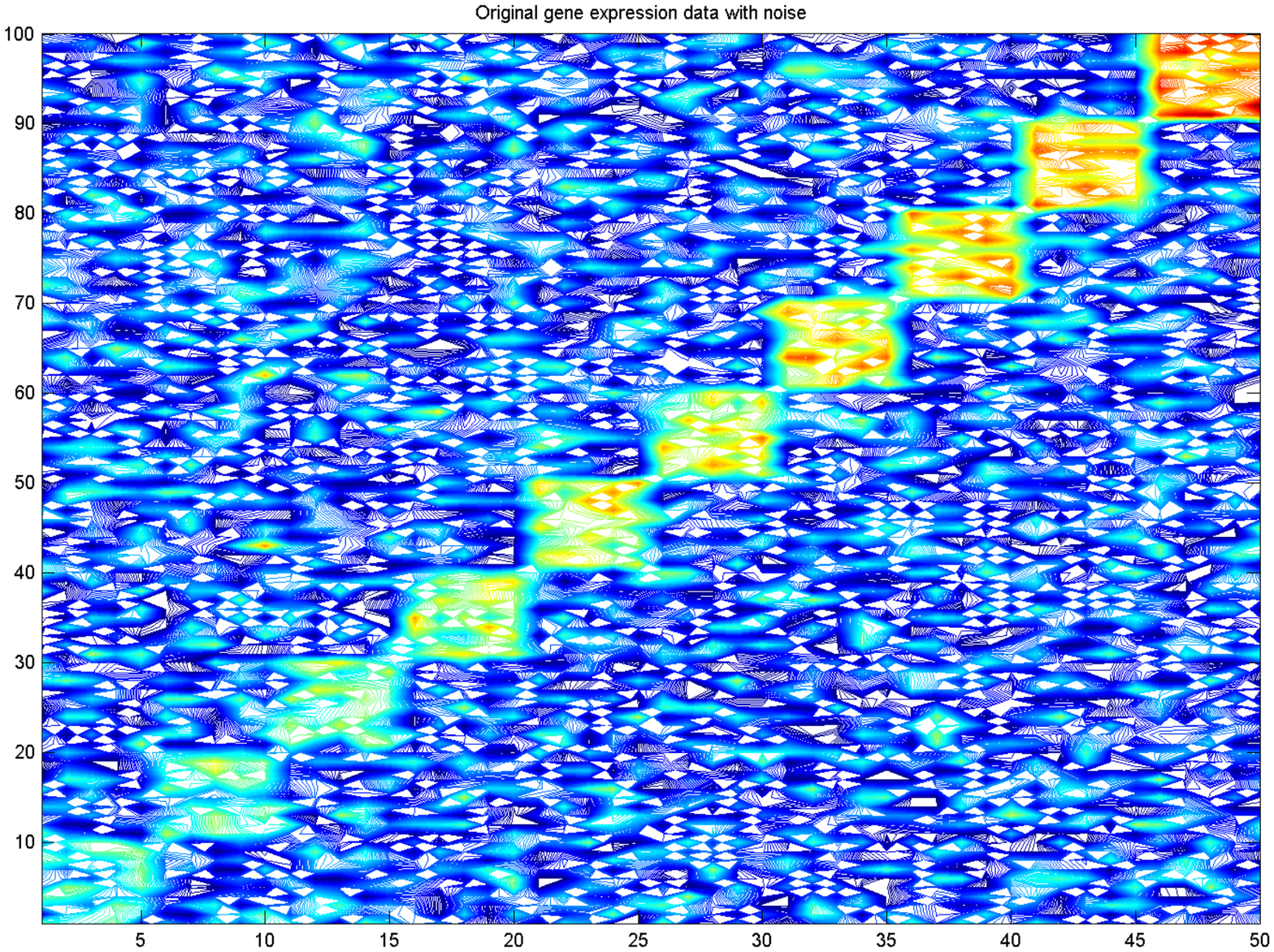}
\caption{A noisy gene expression data matrix with $\sigma=0.3$}\label{fig:noisedata}
\end{figure} 

\begin{figure}[htp]
\centering
\includegraphics[width=0.8\textwidth]{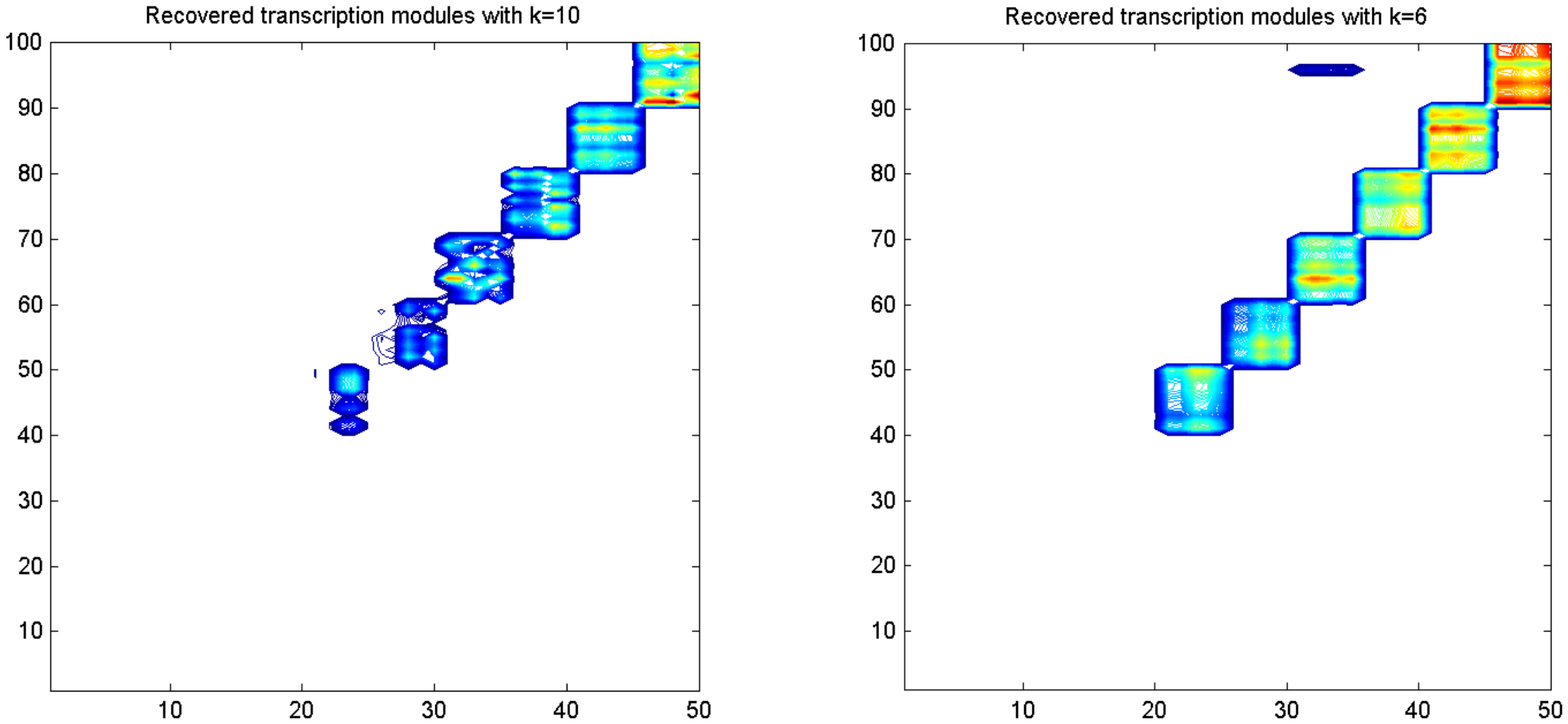}
\caption{Recovery modules obtained with different $k$}\label{fig:diffk}
\end{figure}

We conclude this section with a remark regarding algorithms used to solve the optimization problem \refs{eq:nprob}. For the numerical examples discussed in this section, we solve its equivalent semidefinite optimization formulation \refs{eq:sdpprob} that involves semidefinite constraints for matrices of size $(m+n)\times (m+n)$. For instances with $m=50$ and $n=100$, the computational time in 64-bit Matlab 2013b with the CVX solver on our machine (3.50 GHz CPU and 16.0 GB RAM) is approximately 130 seconds. Clearly, for larger instances, we would need to develop appropriate first-order algorithms for the problem. A similar algorithmic framework as the one in Doan et al. \cite{doan2013proximal} developed for the nuclear norm formulation could be an interesting topic for future research. 
}
\section{Conclusions}
We have shown that a convex optimization problem with Ky Fan $2$-$k$-norm and $\ell_1$-norm can recover the $k$ largest blocks of nonnegative block diagonal matrices under the presence of noise under certain conditions. This is an extension of the work in \cite{Doan10} and it could be used in biclustering applications.

\section*{Acknowledgements}
We would like to thank two referees for their helpful comments and suggestions.

\renewcommand{\baselinestretch}{1.00}
\small
\bibliographystyle{plain}
\bibliography{LCK}
\end{document}